\documentclass[11pt]{article}
\usepackage{latexsym,amsmath,amscd,amssymb,amsthm,graphics}
\usepackage{subcaption}
\usepackage{enumerate,todonotes}
\usepackage[margin=2cm]{geometry}
\usepackage{graphicx}
\usepackage[square,numbers]{natbib}
\usepackage{hyperref}
\hypersetup{colorlinks}
\usepackage{cleveref}
\usepackage{url}
\usepackage{framed, comment}
\usepackage[all]{xy}
\usepackage{epstopdf}
\usepackage{float}
\usepackage{stmaryrd}

\newcommand{\wt}[1]{\widetilde{#1}}

\newcommand\addtag{\refstepcounter{equation}\tag{\theequation}}

\newtheorem{theorem}{Theorem}[section]
\newtheorem*{theorem*}{Theorem}
\newtheorem{definition}[theorem]{Definition}

\newtheorem{remark}[theorem]{Remark}

\theoremstyle{method}
\newtheorem{method}{Method}[section]

\AddToHook{cmd/appendix/before}{

}

\def\bx{\boldsymbol{x}}

\makeatother\def\b{\boldsymbol}

\newcommand{\cred}[1]{\textcolor{red}{#1}}

\begin{document}

\title{Lévy areas, Wong Zakai anomalies in diffusive limits of Deterministic Lagrangian Multi-Time Dynamics}

\author{Theo Diamantakis$^1$ James Woodfield$^1$}
\addtocounter{footnote}{1} 
\footnotetext{Department of Mathematics, Imperial College, London SW7 2AZ, UK. 
\addtocounter{footnote}{1} }

\date{\today}

\maketitle

\makeatother


\begin{abstract}
Stochastic modelling necessitates an interpretation of noise. In this paper, we describe the loss of deterministically stable behaviour in a fundamental fluid mechanics problem, conditional to whether noise is introduced in the sense of Itô, Stratonovich or a limit of Wong-Zakai type. We examine this comparison in the wider context of discretising stochastic differential equations with and without the Lévy area. We study Stochastic Advection by Lie Transport and its derivation from homogenisation theory, which introduces drift corrections of the same class naturally. From a numerical viewpoint, we demonstrate performing higher order discretisations with the use of a Lévy area can also lead to the loss of conserved area and angle quantities. Such behaviour is not physically expected in the Stratonovich model.
From the viewpoint of homogenisation, the qualitative properties of the Wong-Zakai anomaly are physically motivated as arising due to correlations from a fast and mean scale fluid decomposition.



\end{abstract}

\newpage

\section{Introduction}

Data driven stochastic parametrisations for
geophysical fluids dynamics are an effective and computationally tractable approach to accounting for unresolved scales \cite{cotter2017}. Faster and smaller scales are crucial to accurate forecasting, but direct simulation is neither practical nor necessary to exhibit such effects in solution behaviour. SALT- (Stochastic Advection by Lie Transport) is a stochastic parametrisation \cite{SALTalgo} based on the particle-trajectory mapping in fluid mechanics that conserves quantities known as Casimirs \cite{holm2015variational}, examples include helicity and enstrophy.


The SALT ansatz describing the stochastic fluid velocity is given by
\begin{equation}
\mathrm{d}x_t = u_t(x_t) \mathrm{d}t + \sum_{k=1}^m\xi_k(x_t) \circ \mathrm{d} W^k_t\, ,\label{saltansatz}\end{equation}

where $x_t$ denotes a Lagrangian fluid particle position, $u_t, \{\xi_k\}_{k=1}^m$ are vector fields that model the resolved fluid velocity and unresolved effects respectively and $W_t = \{W_t^k\}^m_{k=1}$ is an $m$ dimensional standard Brownian motion. A variety of subgrid effects may be parametrised, including eddy currents, wave interactions or even phytoplankton blooms. In a rapidly changing global climate one expects these oceanic processes to certainly vary. One can interpret this as a change in the calibration vector fields $\xi_k$, but also a change in the underlying noise $\circ \mathrm{d}W$ from correlations to even a weakening of the Brownian assumption (i.e coloured noise or a fractional Brownian motion). In particular, the effects of such a coloured noise have been shown to cause bifurcations in mechanical systems of SALT type in the recent work of \cite{DHP2023}, and in the real world, much attention has been drawn to the possibility of a collapse in the Atlantic meridional overturning circulation (AMOC) modelled in \cite{Ditlevsen2023} as a bifurcation of a stochastic dynamical system. 

In this paper we investigate an instance of this noise driven bifurcating term, the Wong-Zakai anomaly. Referring to the Wong-Zakai theorem \cite{wong1965convergence}, one has that for a given smooth convergent approximation of Brownian motion $W^n(t) \rightarrow W(t)$, we recover convergence of the solutions controlled differential equations driven by $W^n$. Given $X^n_t, X_t$ solving the following equations below, we have $X^n_t \rightarrow X_t$, \[\dot{X}^n_t = \sum_{i=1}^mA_i(X^n_t)\dot{W}_i^n(t),\]  \[\mathrm{d}X_t = \sum_{i=1}^mA_i(X_t) \circ \mathrm{d}W^i_t,\] provided the fixed vector fields $\{A_i\}^m_{i=1}$ are sufficiently smooth \emph{and commute}\footnote{If the vector fields don't commute, one can attain Stratonovich convergence when approximations to the Weiner process $W$ are more specific i.e. piecewise linear and mollifier approximations to $W$ are used.}.

It was shown by Ikeda-Watanabe~\cite[Chapter 7, Theorem 7.2]{ikeda2014stochastic} the possibility of approximations $W^n$ other than piecewise linear or mollifications of $W$ (for an example of such a construction see \cite{sussmann1991limits}) converging uniformly in time, $L^2$ in probability such that when $\{A_i\}^m_{i=1}$ fail to commute as vector fields, the limit process $X^n(t) \rightarrow Y_t \neq X_t$ solves the equation \[\mathrm{d}Y_t = \sum_{1 \leq i \leq j \leq m}s_{ij}[A_i, A_j](Y_t) \mathrm{d}t +  \sum_{i=1}^mA_i(Y_t) \circ \mathrm{d}W^i_t.\] Where $[\cdot, \cdot]$ denotes the commutator of vector fields and $s_{ij}$ a skew symmetric constant matrix particular to the approximation constructed in \cite{ikeda1977class}. The additional term not present in the smooth differential equation is what we refer to as a Wong-Zakai anomaly. Our usage of ``anomaly" follows the terminology of rough paths literature referring to $s_{ij}$ as an ``area anomaly". In this context, $s_{ij}$ obstructs the Stratonovich convergence seen in \cite{wong1965convergence} and when $m=2$ can be interpreted as changing the signed area enclosed by a planar Brownian motion and the chord between $0$ and $W_t$. Although our language emphasises the Stratonovich integral as the ``default" choice, as is common in the stochastic fluid transport communities, we must remark a non-vanishing Wong-Zakai anomaly is the true generic case for Brownian motions constructed through functional central limit theorems, as shown in \cite{gottwald2022timereversibility, engel2024nonlinear} with the Stratonovich integral appearing only exceptionally.  

The functional central limit theorem, and more generally homogenisation theory, constructs a smooth approximation of SDE where the approximation to Brownian motion is given by fast scale dynamics; consequently it exhibits the Wong-Zakai anomaly. In the context of geometric mechanics, Cotter, Gottwald, Holm \cite{cotter2017} a posteriori derive the SALT stochastic velocity ansatz \eqref{saltansatz} introduced in \cite{holm2015variational} from a coupled fast slow system using homogenisation theory, they hint that Wong-Zakai anomaly additional drift terms may arise in the more general setting of homogenisation theory (see equations 4.9, 4.10 in chapter 4 in \cite{cotter2017}). In this paper we adapt this derivation in the more geometric language of \cite{Holm2019a, Holm2019b} and show that similar assumptions of the type seen in of Cotter, Gottwald and Holm lead to exact expressions for the drift and diffusion coefficients as averaged flow velocities and spacio-temporal eigenvectors/values perhaps meaningful to the SALT algorithm in \cite{SALTalgo}. The additional drift from homogenisation theory, recognised in \cite{cotter2017}[Eq. 4.9], will be shown to split into a Wong-Zakai anomaly and Itô correction term, motivating the preferred Stratonovich interpretation of noise for this model. 

As remarked previously, Wong-Zakai anomaly drift terms are linked to stochastic analysis and approximation of SDE. For Brownian motion, the phenomenon of perturbing the Lévy stochastic area is a common occurrence in (cross correlated) coloured noise approximations \cite{DHP2023, Friz_2015} (a special case of homogenisation with a fast scaled Ornstein-Uhlenbeck process). This is considered more abstractly in the rough path context where although $W^n(t) \rightarrow W_t$ may be proven on the level of a Brownian path, we fail to obtain convergence of the canonically defined $\int^t_0 W^n(s) \mathrm{d}W^n(s)$ to the analogous iterated Stratonovich integrals $\int_0^t W_s \circ \mathrm{d}W_s =: \mathbb{W}_t^{\text{Strat}}$. Instead, one has convergence to $\mathbb{W}_t^{\text{Strat}}$ perturbed by an antisymmetric term $\mathsf{s}t \in \mathfrak{s}\mathfrak{o}(m)$. Under both of these considerations it has been shown that a Wong-Zakai anomaly can arise in the corresponding equation. The appearance of a noise induced drift term, equivalent to what we denote the Wong Zakai anomaly, has also been observed independently by Sussmann and McShane in limits of equations driven by non standard approximations of Brownian motion \cite{sussmann1991limits,mcshane1972equations,ikeda2014stochastic}.

The effect this anomalous drift term has on stochastic fluid dynamics needs to be understood. Already in \cite{Resseguier_Mémin_Heitz_Chapron_2017} additional drifts recovered from dynamic mode decomposition have shown to have very notable practical modelling advantages. Also, further scale-selected drift velocities of different inhomogeneities are considered in \cite{bauer2020deciphering} where the resulting numerical solutions of quasigeostrophic flow have strong evidence of changing the coherent large-scale structures.

Also, arising in a similar manner but an entirely different phenomena altogether, is the effect of approximating and using Lévy's stochastic area in the context of higher order (in this case Stratonovich) numerical integrators (see the supplementary appendix \ref{Sec:Lévy Area contribution in higher order Stratonovich Integrators} and Eq. \eqref{EQ: Stratonovich Taylor Expansion}). In this context iterated bracket of vector fields also appear, but require multiplication by Lévy's stochastic area for a Weiner process. Unlike the Wong-Zakai anomaly occurring at strong order $0$, the contribution of Lévy's stochastic area occurs at strong order $1$. It is an ongoing discussion as to what extent the simulation of the Lévy area is necessary or effective in stochastic geophysical simulations \cite{boulvard2023diagnostic}.


Point vortex dynamics have been a classical playground for fluid dynamics \cite{aref2007point}, bearing similar underlying geometric structure to the infinite dimensional fluid dynamics problem but in a lower dimensional space under a singular momentum map \cite{MARSDEN1983}. There is a growing interest and precedence in numerically investigating classical fluid instabilities under the effect of transport noise, in \cite{flandoli2022effect} the Kelvin Helmholtz instability was investigated, and a delay in the onset of instability could be achieved. In this paper we investigate both the Wong-Zakai anomaly drift and the Lévy area appearing in higher order numerical schemes, and its effect on point vortex dynamics. Having derived a Stratonovich noise with Wong-Zakai anomaly, we examine the effect of this noise induced drift on stable configurations of stochastic point vortex solutions of the Euler equations, one of the most fundamental non steady noisy fluid models, and show the loss of stable equilibrium configurations. 

\paragraph{Outline and objectives of this paper.} This paper studies the SALT model as homogenised limit of Lagrangian multi-time dynamics and the effect of the Wong-Zakai anomaly it can exhibit. It is organised in several sections.
\begin{itemize}
    \item In section \ref{sec:2}, we state the general SALT variational principle for stochastic Euler-Poincaré equations on an abstract Lie algebra. We define key notations and in section \ref{sec:homogenisation} interpret the stochastic drift and diffusion coefficients as average components of a mean and fast map through homogenisation theory, and show that a Wong-Zakai anomaly and Itô correction term naturally arises. 
    \item In section \ref{sec:continuumeuler} we specialise the ideas of \cref{sec:2} to the case of diffeomorphism groups and examine the Wong-Zakai anomaly's geometric properties. We further specialise to the case of point vortex solutions of the Euler equations in \ref{sec:SPV2D} using the stochastic velocity with Wong-Zakai anomaly ansatz motivated via homogenisation in \ref{sec:homogenisation}. We examine the effects of the Wong-Zakai anomaly on the conserved quantities inherent to the deterministic system and contrast its effects to the well known Itô-Stratonovich correction in \ref{Example:Stratonovich Point Vortex System}, \ref{Example:Itô Point Vortex System} and \ref{Example:Type I }. The choice of point vortices provide us with a solid benchmark and we examine the numerical and dynamical aspects of the Lévy area and Wong-Zakai anomalies. 
    \item In \cref{Numerical investigation} we consider the implications of the Lévy area in the context of integrators of stochastic systems. In section \ref{sec: numerics}, we use an Additive Runge Kutta method of the point vortex system consistent with several interpretations of noise. We also discuss the required modification in order to treat the Wong-Zakai anomaly drift(emerging from peturbed Lévy areas) and also investigate the implications of approximating and using Lévys stochastic area for a Weiner process in the context of higher order numerical methods. \Cref{sec:results}, contains the results of the point vortex experiments, where particular interest was focused on preserving the area between three point vortices under specific rotational and translational vector fields. 
    \item In the supplementary appendix \ref{sec:appendixh} we compile a selection of results from homogenisation theory used in this paper. We also comment on the links with rough paths, the Wong-Zakai theorem and the extensions that can be pursued in the context of this paper, outlined in the conclusion \ref{sec: conclusion}. A derivation of the Wong-Zakai anomaly through a numerical viewpoint is presented in Appendix \ref{Sec:Lévy Area contribution in higher order Stratonovich Integrators}.
\end{itemize}

\section{Stochastic models of ideal fluids}\label{sec:2}

In this section we shall outline the abstract formulation of geometric mechanics with added stochasticity.  

\subsection{Background for SALT} Initially formulated in \cite{holm2015variational}, \emph{Stochastic Advection by Lie Transport} (SALT) allows the inclusion of noise into the transport velocity of the underlying deterministic fluid. The logic in this approach is twofold; from the modelling perspective this allows parametrisation of subgrid scale effects using calibrated ``data vector fields". From the theoretical viewpoint the advantage of SALT noise is the preservation of deterministic Poisson and Hamiltonian structures, an important consequence is SALT obeys a stochastic Kelvin circulation theorem \cite{Luesink2021}. 

In a recent paper of Crisan, Holm, Leahy and Nilssen \cite{crisan2022variational} the SALT approach to geometric mechanics has been further generalised to use rough paths as noise, yielding rough differential equations (RDE). The use of rough paths allows the incorporation of non-Markovian random perturbations to the Lagrangian fluid trajectories which allow for memory effects, as an example. 

In this paper we are largely interested Brownian motion as driving noise, motivated as the limit of homogenised fast chaotic dynamics. Brownian motion is one particular choice of rough path, but also possesses it's own standalone theory of stochastic differential equations (SDE) that are generalised by RDE. We shall adopt a ``rough paths informed" SDE approach, meaning the following.

\begin{itemize}
    \item We shall interpret SDE in sense of either the Stratonovich, Itô or Stratonovich with Wong-Zakai anomalies \cite{sussmann1991limits}. This shall save having to define more advanced constructions from rough paths and improve accessibility for the reader.
    \item On occasion, we may lift the Brownian path $W_t \in \mathbb{R}^m$ to it's enhancement as a rough path that specifies the notion of integration. For example:

    \[\boldsymbol{W}^{\text{Strat}}_{s,t} = (W_t - W_s, \mathbb{W}^{\text{Strat}}_{s,t}), \quad \left(\mathbb{W}^{\text{Strat}}_{s,t}\right)^{ij} := \int_s^t (W_r - W_s)^i \circ \mathrm{d}W^j_r \in \mathbb{R}^m \otimes \mathbb{R}^m.\]

    Formally, we interpret RDEs driven by a Brownian rough path with the appropriate SDE using the symbolic identification $``\mathrm{d}\boldsymbol{W}^{\text{Strat}} = \circ \mathrm{d}W"$. Changes to the signature $\mathbb{W}^{\text{Strat}}$ are reflected in how the SDE must be interpreted. For example, $\mathbb{W}^{\text{Strat}}_{s,t} -\frac{1}{2}I_{m \times m}(t-s) = \mathbb{W}^{\text{Itô}}_{s,t}$, where:

    \[ \left(\mathbb{W}^{\text{Itô}}_{s,t}\right)^{ij} := \int_s^t (W_r - W_s)^i \mathrm{d}W^j_r, \quad ``\mathrm{d}\boldsymbol{W}^{\text{Itô}} = \mathrm{d}W".\]

    This corresponds to the quadratic covariation correction terms between Itô and Stratonovich integrals. A rough path $\boldsymbol{Z}_{s,t} = (Z_t - Z_s, \mathbb{Z}_{s,t})$ who's notion of rough integration obeys ``standard" integration by parts with no additional terms is known as geometric, and must satisfy the identity $\operatorname{Sym}(\mathbb{Z}_{s,t}) = \frac{1}{2} Z_{s,t} \otimes Z_{s,t}$. This holds true for Stratonovich enchanced Brownian motion but not for Itô. The antisymmetric part is known as the Lévy stochastic area when $Z_t = W_t$ and unlike the symmetric part is unconstrained by $W_t$. The Lévy stochastic area for Itô and Stratonovich Brownian motion coincide. 

    \item When the choice of enhancement of the path $W_t$ is yet to be made we shall simply write $*\mathrm{d}W$ to denote that the rough driver is $(W_t - W_s, \mathbb{W}_{s,t})$ with $\mathbb{W}_{s,t} \in \{\mathbb{W}^{\text{Strat}}_{s,t}, \mathbb{W}^{\text{Itô}}_{s,t},  \mathbb{W}^{\text{Strat}}_{s,t} + \mathsf{s}(t-s)\}$ allowing various choices of enhancement as a rough path. The deterministic, skew-symmetric $m \times m$ matrix $\mathsf{s}$ is a perturbation of the Lévy stochastic area for Stratonovich enchanced Brownian motion, which incurs a correction term \emph{not} given by a quadratic covariation, we refer to as the Wong-Zakai anomaly (L\'evy area correction in other works). A key point to note is that due to skewsymmetry of $\mathsf{s}$ we still have a geometric rough path where the integration by parts formula matches the smooth case\footnote{The rough integral defined using the enhancement $\mathbb{W}^{\text{Strat}}_{s,t} + \mathsf{s}(t-s)$ can be in fact shown to coincide in value with Stratonovich integration. Due to smoothness and antisymmetry of $F(t-s) = \mathsf{s}(t-s)$, any extra contribution is contracted with a symmetric matrix of quadratic covariations \cite[Section 5.2]{Hairer2020} and a distinction only arises when driving a RDE.}. 
\end{itemize}
The signature is particularly useful when we are considering smooth approximations to Brownian motion $B^\varepsilon$ with $\mathrm{d}B^\varepsilon := \dot{B}^\varepsilon \mathrm{d}t$, as they provide a way of tracking which stochastic integral we will converge to. With our conventions established we are in position to write down a stochastic variational principle.

We consider a Lie group $G$ with Lie algebra $\mathfrak{g}$ and vector space $V$, which allows defining the semidirect product Lie group $G \ltimes V$ via a (left) group action. This action is represented abstractly as concatenation $$\left(g_{1}, v_{1}\right) \cdot \left(g_{2}, v_{2}\right):=\left(g_{1} g_{2}, v_{1}+g_{1} v_{2}\right), \quad vg^{-1} := gv.$$ 

Examples include matrix vector multiplication (a linear representation), or the pushforward of tensor quantities via a diffeomorphism. We take $g_t \in G, u_t \in \mathfrak{g}, \mu_t \in \mathfrak{g}^*$ as arbitrary continuous in time elements, and $\xi_i \in \mathfrak{g}, a_0 \in V$ fixed. The noise/geometric rough path is Brownian motion $B_t$ which is $\mathbb{R}^m$ valued and posses a perturbation to it's Lévy area by a (possibly zero) $\mathfrak{s}\mathfrak{o}(m)$ matrix $\mathsf{s}$. 

The reduced Lagrangian $\ell: \mathfrak{g} \times V \rightarrow \mathbb{R}$ is thought as arising from a parametric Lagragian $L_{a_0} :TG \rightarrow \mathbb{R}$ that is right invariant at all parameter values, that is $L_{a_0}(g, v) = L_{a_0g^{-1}}(e,vg^{-1}) := \ell(vg^{-1}, a_0g^{-1})$. For a detailed discussion of this procedure see \cite{WeinsteinBook2002}. Assuming this reduction has been carried out the following constrained variational principle produces the stochastic Euler-Poincar\'e equations as introduced in \cite{holm2015variational}. The Hamilton-Pontryagin approach was originally derived in \cite{Gay-Balmaz2018} and can also be seen as a special case of the rough Hamilton-Pontryagin variational principle in \cite{crisan2022variational} and \cite{DHP2023}.

The action in the variational principle is, 
\begin{equation}S = 
\int_{0}^{T} \ell(u_t,a_0g_t^{-1}) \mathrm{d}t +\Big\langle \mu_t , \mathrm{d}g_tg_t^{-1} - u_t \mathrm{d}t -\frac{1}{2}\sum_{1\leq i,j \leq m}\mathsf{s}^{ij}[\xi_i, \xi_j]\mathrm{d}t - \sum_{k=1}^m \xi_k \circ \mathrm{d}W^{k}_t \Big\rangle \, ,\label{abstractVarP}\end{equation}

The meaning of these integrals are to be defined as follows, 

\[\int \left\langle \mu_t, \mathrm{d}g_tg_t^{-1} \right\rangle := \int \left\langle \mu_t, \alpha_t \right\rangle \mathrm{d}t + \int \sum_{k=1}^m \left\langle \mu_t, \beta^k_t \right\rangle \circ \mathrm{d}W^k_t \,,\]

\[ \int \left\langle \mu_t, u_t \mathrm{d}t + \frac{1}{2}\mathsf{s}^{ij}[\xi_i, \xi_j] \mathrm{d}t \right\rangle := \int \left\langle \mu_t, u_t + \frac{1}{2}\sum_{1 \leq i,j \leq m}\mathsf{s}^{ij}[\xi_i, \xi_j] \right\rangle \mathrm{d}t \,,\] 

\[\int \langle \mu_t, \xi \circ \mathrm{d}W_t \rangle := \int \langle \mu_t, \xi \rangle \circ \mathrm{d}W_t := \int \sum_{k=1}^m\langle \mu_t, \xi_k \rangle \circ \mathrm{d}W^k_t \,.\]

The vector fields $\alpha_t, \beta^k_t$ are the drift and diffusion coefficients of the SDE defined by the time derivative of the stochastic flow $\mathrm{d}g_t$ on the tangent bundle of the Lie group $TG$. 

\begin{equation}\mathrm{d}g_t = \alpha_t(g_t) \mathrm{d}t + \beta_t(g_t) \circ \mathrm{d}W_t \,.\label{diffeovariation}\end{equation}

The variation of this action $\delta S = 0$ results in the following stochastic semidirect product Euler-Poincar\'e equations,
\begin{equation}
\begin{aligned}
&\mathrm{d} \frac{\delta \ell}{\delta u_t} + \operatorname{ad}^*_{u_t}\frac{\delta \ell}{\delta u_t} \mathrm{d}t + \frac{1}{2}\sum_{1 \leq i,j \leq m} \operatorname{ad}^*_{[\xi_i, \xi_j]}\frac{\delta \ell}{\delta u_t} \mathrm{d}t + \sum_{k=1}^m \operatorname{ad}^*_{\xi_k}\frac{\delta \ell}{\delta u_t}\circ \mathrm{d}W^{i}_t = \frac{\delta \ell}{\delta a_t} \diamond a_t \mathrm{d}t \,, \\
&\mathrm{d}a_t + \left(u_t + \frac{1}{2}\sum_{1\leq i,j \leq m}\mathsf{s}^{ij}[\xi_i, \xi_j]\right)a_t \mathrm{d}t + \sum_{k=1}^m \xi_k a_t\circ \mathrm{d}W^{k}_t\,.
\end{aligned}\label{abstractEP}\end{equation}

The operator $\operatorname{ad}^*$ appearing in the Euler-Poincaré equation is the coadjoint operator, defined as the dual of the map $\operatorname{ad}_\xi \eta :=[\xi,\eta]$ which is the bracket associated to $\mathfrak{g}$.

\begin{align}
\begin{split}
\mathrm{ad}^{*}: \mathfrak{g} \times \mathfrak{g}^{*} & \longrightarrow \mathfrak{g}^{*}, \\
(\xi, \mu) & \mapsto \operatorname{ad}_{\xi}^{*}(\mu),\\
 \left\langle\operatorname{ad}_{\xi}^{*} \mu, \eta\right\rangle & =\left\langle\mu, \operatorname{ad}_{\xi} \eta\right\rangle, \quad \text { for all } \eta, \xi \in \mathfrak{g} .
\end{split}
\label{def: ad-star}
\end{align}

The advected quantity terms appearing in the right hand side of \eqref{abstractEP} make use of the diamond operator and the Lie algebra action on $V$. Differentiating the group multiplication on $G \ltimes V$ determines a $\mathfrak{g}$ action on $V$; the diamond operator is then defined as map that defines a dual of this action from the vector space pairing to the Lie algebra pairing as follows,

$$\langle v \diamond a, \xi \rangle_{\mathfrak{g}^*\times\mathfrak{g}} = \langle v, \,\xi a \rangle_{V^*\times V} \quad{\text { for all } \xi \in \mathfrak{g}, a \in V, v \in V^*.}$$

Specific instances of these operators are to be discussed in following \cref{sec:EPeuler}. 

In order to make sense of the Euler-Poincar\'e equations as arising from a variational principle, one requires the noise to be given by a geometric rough path in order to make use of standard product and chain rules of calculus as shown in \cite{crisan2022variational}. This rules out a choice of Itô noise, which is non geometric\footnote{One can derive an Itô SDE by writing in the Stratonovich interpretation and always including the relevant correction term, see \cite{Holm2020a}.}. The most general lift of Brownian motion to a geometric rough path that we shall consider is the perturbation of the Stratonovich signature $\mathbb{B}^{\text{Strat}}$ by the deterministic constant matrix $\mathsf{s} \in \mathfrak{s}\mathfrak{o}(m)$. This remains geometric as a consequence of skew symmetry, since:

\[\operatorname{Sym}(\mathbb{B}^{\text{Strat}}_{t,s} + \mathsf{s}(t-s)) = \frac{1}{2}(B_{t} -B_s) \otimes (B_{t} - B_s).\]

One can interpret the $\mathsf{s} \neq 0$ case as a rough differential equation driven by the geometric rough path $(B_t - B_s, \mathbb{B}^{\text{Strat}}_{s,t} + \mathsf{s}(t-s))$ or a Stratonovich equation with drift. In keeping with the ``rough paths informed" SDE approach we will favour the interpretation of a Stratonovich SDE with drift. In \cref{sec:homogenisation} we will explain how nonzero $\mathsf{s}$ can appear naturally in the context of deriving stochastic velocities through homogenisation theory.

The Lagrange multiplier $\mu_t \in \mathfrak{g}^*$ constrains $\alpha_t = u_t + \frac{1}{2}\mathsf{s}^{ij}[\xi_i, \xi_j], \beta_t = \xi$ in \eqref{diffeovariation} to produce the stochastic velocity relation:
\begin{equation}
\mathrm{d}g_tg_t^{-1} = u_t \mathrm{d}t + \frac{1}{2}\mathsf{s}^{ij}[\xi_i, \xi_j] \mathrm{d}t + \xi \circ \mathrm{d}W_t.\label{lagpaths}\end{equation}

The replacement of the deterministic ODE reconstruction equation by a SDE with additional Lie algebra elements $\xi_i \in \mathfrak{g}$ can be motivated through homogenisation theory first derived in \cite{cotter2017}. This will be discussed in section \ref{sec:homogenisation} along with the modifications caused by perturbations to the L\'evy area from the homogenisation procedure.

\subsubsection{Example: Stochastic Incompressible Euler Equations}\label{sec:EPeuler}

A multitude of stochastic fluid models can be derived as a specific case of the Euler-Poincaré equations defined in \eqref{abstractEP}, such as the Euler-Boussinesq equations \cite{Pan2023}, the Camassa Holm equation \cite{Tyranowski2016} (of which relates to the Burgers, KdV and Hunter-Saxton equations) and (isentropic or incompressible) Euler equations \cite{crisan2022variational, Crisan2019, takao2020stochastic} via particular choices of Lagrangian $\ell$ and group $G$. We shall demonstrate the incompressible Euler case, first for generic Riemannian manifolds and the familiar case in flat Euclidean space.

The appropriate (formal) Lie group to consider is $G = H^s(M, M)$, the volume preserving Sobolev diffeomorphisms with $s$ weak derivatives in $L^2$ on a manifold $M$ \cite{Marsden1984}. Let us make explicit some of constructions we have stated above for the specific case of this group. For diffeomorphism groups the operation is given by composition and our diffeomorphisms act on a vector space $V$ of tensor fields via pushforward, i.e. $a_t := g_t a_0 := (g_t)_* a_0 = (g_t^{-1})^* a_0$\footnote{This is also sometimes notated $ag^{-1}$, and for the case of $a$ given by a $0$-form (a scalar) $(g^{-1})^* a$ is indeed interpreted as the composition $a \circ g^{-1}$, however this algebraic notation can be unclear for more general tensor pullback actions (which contain transformation terms of the basis elements), or outside the context of a matrix Lie group acting on vectors via right multiplication.  }.

It follows from these definitions that the Lie algebra is given by divergence free (incompressible) vector fields on $M$, with the Lie algebra bracket given by the usual Jacobi-Lie vector field commutator. The Lie algebra action on $V$ is given by taking the Lie derivative of tensor element $a_t \in V$ with respect to a vector field (such as $u_t$) acting on $a_t$. The Lie derivative is also the same operator used for the coadjoint action $\operatorname{ad}^*$. It is crucial to note however that although both lines in \eqref{abstractEP} will be written in terms of Lie derivatives, $\frac{\delta \ell}{\delta u_t}$ and $a_t$ belong to different spaces and thus will not translate to the same vector calculus operations should one choose to calculate in local coordinates\footnote{Certain exceptions can occur where these may coincide (such as Cartesian coordinates with access to the Euclidean metric and choice of a kinetic energy Lagrangian), but this should not be expected to occur in general.}. The diamond operator can be calculated by integration by parts depending on what exact type of tensor field $a_t$ is chosen to be, since both the $\langle \cdot , \cdot \rangle_{\mathfrak{g}^* \times \mathfrak{g}}, \langle \cdot , \cdot \rangle_{V^* \times V} $ pairings are $L^2$ type pairings involving suitably contracted geometric objects.

Obtaining the incompressible Euler equations require a kinetic energy Lagrangian with the choice of group described above, however, it is common to impose incompressibility of the fluid via a stochastic constraint allowing us to set $G$ as larger group of diffeomorphisms (not necessarily volume preserving) with $s$ weak derivatives. This involves making use of a single advected quantity, the mass density, or volume element, $D_t = (g_t)_* \left(D_0\operatorname{vol}_{\boldsymbol{g}}\right)$, to constrain (the components of) $D_t$ to equal one. Our action integral will take the form:
\begin{equation}
    \int_0^T\left(\ell(u_t, D_t) \mathrm{d}t+\left\langle p_0, D_t-1\right\rangle \mathrm{d} t+\left\langle p,D_t-1\right\rangle \circ \mathrm{d} W_t\right) +\Big\langle \mu , \mathrm{d}g_tg_t^{-1} - u_t \mathrm{d}t - \frac{1}{2}\mathsf{s}^{ij}[\xi_i, \xi_j] \mathrm{d}t - \xi \circ \mathrm{d}W_t \Big\rangle
\end{equation}

It is noted in \cite{street2021semi} that for stochastic fluid dynamics one requires the pressure Lagrange multiplier to contain both deterministic and stochastic constraints $p_0 = p_0(x,t)$ and $p = (p_k(x,t))_{k=1}^{ m}$ for each specified $\xi_k$. The kinetic energy Lagrangian $\ell$ on the Riemannian manifold $(M, \boldsymbol{g})$ is given as \cite{arnold1966geometrie}

$$\ell(u_t, D_t) = \frac{1}{2}\int_M  \boldsymbol{g}(u_t,u_t) D_t\,.$$
The stochastic Euler-Poincaré equations \eqref{abstractEP} in the specific case of this group and Lagrangian are given to be:
$$\begin{aligned}
& (\mathrm{d} + \pounds_{u_t \mathrm{d}t + \frac{1}{2}\mathsf{s}^{ij} [\xi_i, \xi_j]\mathrm{d}t + \xi \circ \mathrm{d}W_t} ) (u^\flat_t \otimes D_t) = \left( \frac{1}{2}\bold{d} \boldsymbol{g}(u_t,u_t) \mathrm{d}t + \bold{d} p_0 \mathrm{d}t + \bold{d} p \circ \mathrm{d}W_t \right) \otimes D_t\\
&(\mathrm{d} + \pounds_{u_t \mathrm{d}t + \frac{1}{2}\mathsf{s}^{ij} [\xi_i, \xi_j]\mathrm{d}t + \xi \circ \mathrm{d}W_t} )D_t = 0, \quad D_t \equiv 1\cdot \operatorname{vol}_{\boldsymbol{g}}
\end{aligned}$$

Here, the variational derivative of the Lagrangian with respect to vector field $u_t \in \mathfrak{X}(M)$ is a 1-form density $\frac{\delta \ell}{\delta u_t} = u^\flat_t \otimes D_t = \mu \in \mathfrak{X}^*$, with $\flat = \flat_{\boldsymbol{g}}$ denoting the metric induced musical isomorphism from vectors to covectors and a bold $\bold{d}$ denotes exterior derivative. It is customary to consider only the 1-form component as is written since $D_t$ passes through the operator $(\mathrm{d} + \pounds_{u_t})$ due to advection, this in turn allows us to formally ``divide" the left and right hand side by $D_t$ terms in this equality of 1-form densities to obtain an equality of 1-forms. Likewise, the diamond term $\frac{\delta \ell}{\delta D_t}$ is a semimartingale  $\frac{1}{2}\boldsymbol{g}(u,u)\mathrm{d}t + p_0 \mathrm{d}t + p_k \circ \mathrm{d}W^k \in V^* = H^s(M, \mathbb{R})$\footnote{This semimartingale may also be denoted by $\frac{\delta \ell}{\delta D_t} \circ \mathrm{d}S_t$, note that we have $s$-times weak differentiability in space.} with $\frac{\delta \ell}{\delta D} \diamond D = \bold{d}\frac{\delta \ell}{\delta D} \otimes D $.

When $M = \mathbb{R}^d$ or $\mathbb{T}^d$ with $\boldsymbol{g}$ given by the Euclidean metric we may drop the overt distinction between 1-forms and vector fields, discard the density part through the ``division by $D_t$" procedure mentioned above, and obtain the following stochastic Euler system (with additional noise induced drift):

\begin{equation}
\begin{aligned}
& \mathrm{d} u_t + u_t \cdot \nabla u_t\mathrm{d}t + \frac{1}{2}\mathsf{s}^{ij} \left([\xi_i, \xi_j] \cdot \nabla u_t + (\nabla[\xi_i, \xi_j])^T u_t \right) \mathrm{d}t + \left(\xi_k \cdot \nabla u_t + (\nabla \xi_k)^T u_t\right) \circ \mathrm{d}W^k_t = \nabla p_0 \mathrm{d}t + \nabla p_k \circ \mathrm{d}W^k_t\\[8pt]
&\nabla \cdot u_t = \nabla \cdot \xi_i = 0.
\end{aligned}\label{EulerwithLevy}\end{equation}

Where we defined the operator $(\nabla v)^T u := \sum_j u_j \nabla v_j$ and made use of $(\nabla u)^T u = \frac{1}{2} \nabla (u \cdot u)$. An analysis of the well-posedness properties of the vorticity formulation (see \eqref{vorticity}) with no Wong-Zakai anomaly can be found in \cite{Lang2023} and the rough case in \cite{CRISAN2022109632}. We also note that the pressure constraint and advection equation initially only imply incompressibility of the sum total stochastic velocity field:
$$\nabla \cdot (u_t \mathrm{d}t + \frac{1}{2}\mathsf{s}^{ij}[\xi_i, \xi_j] \mathrm{d}t + \xi_k \circ \mathrm{d}W^k_t) = 0$$

The Doob-Meyer decomposition implies that each $\xi_k$ and $u_t + \frac{1}{2}\mathsf{s}^{ij}[\xi_i, \xi_j]$ individually are divergence free. Since each $\xi_k$ is forced to be incompressible, and incompressible vector fields form a Lie-subalgebra of $\mathfrak{X} = \operatorname{Vect}(M)$, it also follows that $\nabla \cdot [\xi_i, \xi_j] = 0$, so we recover the usual requirement that $\nabla \cdot u_t = 0$.

\subsection{Homogenisation of deterministic geometric fluid dynamics}\label{sec:homogenisation}

In section \ref{sec:2} we stated a SALT variational principle with an extra drift $\frac{1}{2}\mathsf{s}^{ij}[\xi_i, \xi_j] \mathrm{d}t$. In this section we shall motivate the additional drift and stochastic reconstruction equation used in the variational principle \eqref{abstractVarP} as a limit of (possibly deterministic) multi-time dynamics, through homogenisation. Stochastic flows have been considered by Mikulevicius and Rozovskii \cite{doi:10.1137/S0036141002409167}, and models such as Location Uncertainty (LU), akin to SALT, also postulate a stochastic velocity field with a drift and noise term \cite{npg-30-237-2023}. The modelling assuming that the noise in these models ought to be the limit of a fast fluctuation flow will produce analogous drift terms in these models.

In the SALT model, homogenisation theory was first applied to the reconstruction equation $\dot{g}_tg^{-1}_t = u_t$ in Cotter, Gottwald and Holm \cite{cotter2017} under a composition of two maps $g_{t, t/\varepsilon} = (\operatorname{id} + \zeta_{t/\varepsilon}) \circ \overline{g}_t$, a mean flow and near identity fast map. 
Weak convergence to a stochastic (in time) diffeomorphism is shown under some chaoticity assumptions, which are natural modelling assumptions in the context of an ocean or atmospheric system. 

A key remark made in \cite{cotter2017}[Eq. 4.9] is the multiplicative nature of the noise arising from this decomposition induces drift correction terms depending on the fast map $\zeta_{t/\varepsilon}$. We shall obtain the result of \cite{cotter2017} using some of the more geometric perspectives seen in \cite{Holm2019a} \cite{Holm2019b} and explicitly identify the noise induced correction term noted in \cite{cotter2017} to be a combination of the Itô-Stratonovich correction term and a Wong-Zakai anomaly. The presence of Wong-Zakai anomalies in the derivation of SALT is particular of interest, as the effect of perturbations to the velocity of finite dimensional physical systems such as the SALT rigid body has been studied in \cite{DHP2023}. In this paper the setting of fluid dynamics motivates carrying over these ideas to the setting of infinite dimensional Lie groups. The coloured noise approximations and their limits containing Wong-Zakai anomalies seen in \cite{DHP2023} are in fact a special case of the generalisation in this paper. One may interpret the coloured noise limit as stochastic averaging of an Ornstein Uhlenbeck process, a form of (stochastic) homogenisation \cite{Pavliotis2008}. Arguments involving limiting time scales to relate deterministic fluid models to a stochastic equation have also been considered in \cite{w12102950}, notably, the Wong-Zakai theorem is also invoked. We anticipate the relevance of these drift corrections for stochastic fluid models other than SALT. 

The composition of maps theory and equivalent decompositions of Eulerian fluid velocity have been used in the context of wave mean flow interaction (WMFI) \cite{HHS2023} as well in the Generalised Lagrangian Mean (GLM) theory seen in \cite{Holm2002}.
We follow the language of this approach, where $G$ is taken as a diffeomorphism group and we consider curves $g_t$ through $G$ such that $g^{\varepsilon}_{t} = \widetilde{g}_{t/ \varepsilon} \circ \overline{g_t}$. $\widetilde{g}_{t/\varepsilon} \in G$, $\overline{g}_t \in G$ denote the fast and slow map respectively.
For simplicity we take $M = \mathbb{R}^d$. The full velocity $u \in \mathfrak{X}(\mathbb{R}^d)$ is defined as the vector field satisfying $\frac{\partial}{\partial t} g_t(X) = u_t(g_t(X))$ where $X = g_0(X) \in \mathbb{R}^d$ is the initial particle label. 

We follow a calculation seen in \cite{Holm2019a}, \cite{Holm2019b}.  Taking the time derivative of this decomposition gives, by the chain rule:
\begin{align}
u^\varepsilon_t =: \dot{g}^\varepsilon_t \circ (g_t^{ \varepsilon})^{-1} &= \frac{1}{\varepsilon}\dot{\widetilde{g}}_{t / \varepsilon} \circ \overline{g}_t \circ ( \widetilde{g}_{t / \varepsilon} \circ  \overline{g}_t)^{-1} + \left(\frac{\partial \widetilde{g}_{t / \varepsilon}}{\partial \overline{g}_t} \cdot \dot{\overline{g}}_t\right) \circ (\widetilde{g}_{t / \varepsilon} \circ \overline{g}_t)^{-1}\\
&= \frac{1}{\varepsilon}\dot{\widetilde{g}}_{t / \varepsilon} \circ \widetilde{g}^{-1}_{t / \varepsilon} + \widetilde{g}_{t / \varepsilon \, *} ( \dot{\bar{g}}_t \circ \bar{g}^{-1}_t) =: \frac{1}{\varepsilon}\dot{\widetilde{g}}_{t / \varepsilon} \widetilde{g}^{-1}_{t / \varepsilon} + \operatorname{Ad}_{\widetilde{g}_{t/\varepsilon}}\overline{u}_t \end{align}




Evaluating at a Lagrangian coordinate $x_t = g_t(X)$ turns the above relation into an ordinary differential equation.

\begin{equation}
u^\varepsilon_t(x^\varepsilon_t) = \dot{g}_t{g}^{-1}_t(x^\varepsilon_t) = \dot{x}^\varepsilon_t =  \frac{1}{\varepsilon}\dot{\widetilde{g}}_{t / \varepsilon}\widetilde{g}^{-1}_{t/\varepsilon}({x}^\varepsilon_t) + \operatorname{Ad}_{\widetilde{g}_{t/\varepsilon}}\overline{u}_t (x^\varepsilon_t)\label{velocitydecomp}\end{equation}

Our aim is to apply deterministic homogenisation to take $\varepsilon \rightarrow 0$. The scale separation appearing in \cref{velocitydecomp} is appropriate for homogenisation, however the $\varepsilon$ dependencies in $\wt{g}_{t/\varepsilon}$ preclude the use of this theory as is currently written. We introduce a fast variable $\lambda_t^\varepsilon \in \mathbb{R}^m$ and specify the form of the velocity associated to $\wt{g}$ such that \cref{velocitydecomp} is an equation only in time and the slow and fast variables $(x^\varepsilon, \lambda^\varepsilon)$.

A skew product ODE follows from two further assumptions on the structure of the terms in the equation for $u^\varepsilon_t$. 

\begin{enumerate}
    \item 

    The fast vector field $\dot{\widetilde{g}}_{t / \varepsilon}\widetilde{g}^{-1}_{t/\varepsilon}(\cdot)$ is assumed to have the following eigenvalue eigenvector decomposition seen in \cite{cotter2017}
    \begin{align}
    \dot{\widetilde{g}}_{t / \varepsilon}\widetilde{g}^{-1}_{t/\varepsilon}(\cdot) = \sum_{k=1}^m \lambda^{\varepsilon, k}_{t} v_k(\cdot ) \,\,.
    \end{align} 
    Where $v_k: \mathbb{R}^d \rightarrow \mathbb{R}^d$ denotes the $k$-th eigenvector velocity field, assumed spatially dependent and temporally independent. The fast in time spatially independent $k$-th eigenvalue $\lambda_{t}^{\varepsilon,k}:\mathbb{R} \rightarrow \mathbb{R}$, is specified by a system of ODEs. More specifically, it is assumed that $\lambda_{t}^{\varepsilon}$ (the vector who's $k$'th component is $\lambda_t^{\varepsilon,k}$) solves the following initial value problem  
    \begin{align}
        \dot{\lambda}^\varepsilon_t = \frac{1}{\varepsilon^2}h(\lambda^\varepsilon_t), \quad \lambda_0^{\varepsilon}=\lambda_0,\label{eq: ivp for lambda}
    \end{align}
    with initial condition $\lambda_0$,
    and function $h:\mathbb{R}^{m} \rightarrow \mathbb{R}^m$. The system of ODE's is assumed well posed, and we denote the induced flow map $\phi_t:\mathbb{R}^{m}\rightarrow\mathbb{R}^m,$ describing the solution of \cref{eq: ivp for lambda}.
    
    

    
    It is assumed that there exists a unique ergodic invariant measure $\nu$ to which the initial value $\lambda_0$ is distributed according to and that under this measure, the eigenvalues $\lambda^{k, \varepsilon}_t$ are $\nu$-mean zero (for all times $t$ as a consequence of invariance):

    \[\langle \lambda^{k,\varepsilon}_t \rangle_{\nu(\mathrm{d}\lambda)} := \int \lambda_t^{k, \varepsilon} \nu(\mathrm{d}\lambda_0) = \int \lambda^k_0 \nu(\mathrm{d}\lambda_0)=  0\]

    This ensures the map $f_0(x_t^\varepsilon, \lambda^\varepsilon_t) = \lambda^{k,\varepsilon}_t v_k(x^\varepsilon_t)$  satisfies the centering condition $\langle f_0 \rangle_{\nu(\mathrm{d}\lambda)} = 0$. 

    \item We also assume the pushed forward mean velocity can be expressed as a function of $\lambda_t^\varepsilon$ and $x^\varepsilon_t$, that is, there exists an $f_1$ such that:
\[\operatorname{Ad}_{\widetilde{g}_{t/\varepsilon}}\overline{u}_t(x^\varepsilon_t) = f_1(x^\varepsilon_t,\lambda^\varepsilon_t,t).\]
\end{enumerate}




\begin{remark} The second assumption enforces that the mean velocity field, when pushed forward by the fast map at time $t$ has no fast $t/\varepsilon$ temporal dependence outside of $\lambda^\varepsilon_t$. The assumption is trivially satisfied if the flows $\widetilde{g}, \overline{g}$ commute at all times, as this implies $\operatorname{Ad}_{\widetilde{g}_{t/\varepsilon}}\overline{u}_t(x) = \overline{u}_t(x) = f_1(x,t)$. We leave for future works the further weakening of this assumption. More general $\varepsilon$ dependence is permitted in \cite[Theorem 5.5]{10.1214/21-AIHP1203} and alternative approaches via rough paths stability can circumvent the above assumption.

\end{remark} 



With the assumptions 1. and 2. equation \eqref{velocitydecomp} matches the form of equation 1.1 seen in \cite{kelly2017chaos} and 11.2.1 in \cite{Pavliotis2008}, (also stated in equation \eqref{xeps} in the supplementary appendix \ref{sec:appendixh}):

\begin{equation}
\dot{x}^\varepsilon_t = \frac{1}{\varepsilon}f_0(x^\varepsilon_t, \lambda^\varepsilon_t)  + f_1(x^\varepsilon_t, \lambda^\varepsilon_t, t)\label{skewfluidsystem}\end{equation}
\begin{equation}
\dot{\lambda}^\varepsilon_t = \frac{1}{\varepsilon^2}h(\lambda^\varepsilon_t), \quad \lambda_0^{\varepsilon}=\lambda_0.
\label{fastevals}\end{equation}

$$\begin{aligned}
f_1\left(x^\varepsilon_t, \lambda_t^{\varepsilon}, t\right) & =:\frac{\partial \widetilde{g}_{t / \varepsilon}}{\partial \bar{g}_t}\overline{u}_t(\widetilde{g}^{-1}_{t/\varepsilon} x^\varepsilon_t) \\
f_0(x^\varepsilon_t, \lambda_t^{\varepsilon}) &=: \sum_{k=1}^m v_k\left(x^\varepsilon_t\right) \lambda_t^{\varepsilon, k}, 
\end{aligned}$$

with $f_0$ satisfying the centering condition with respect to the measure $\nu$. It is in the above form where Theorem 1.1 in \cite{kelly2017chaos} and the methods of \cite[Section 11.2]{Pavliotis2008} apply when the unspecified dynamics for $\lambda^\varepsilon_t$ are assumed to satisfy an iterated weak invariance principle (WIP). For precise details on the specific class of dynamical systems that satisfy such a WIP we refer the reader to \cite{kelly2017chaos, 284645bd-f971-3dfe-b374-dab0f8c0c677}. 

We note that the $t$-dependence of $f_1$ is of no consequence to matching the form \eqref{xeps}, as the phase space of the slow dynamics $\dot{x}^\varepsilon_t$ may be expanded to $\dot{z}^\varepsilon_t = (\dot{x}^\varepsilon_t, \dot{T}_t)$ with the trivial equation $\dot{T}_t = 1$ appended. The trivial variable $T_t = t + T_0$ is $\mathcal{O}(1)$, passes through $\nu$-averaging unchanged and does not interact with the $\mathcal{O}(1/\varepsilon)$ term $f_0$. This ensures the inclusion of $T_t$ will not affect the centering condition or introduce any difference in the formulae from the $t$-independent $f_1$ case.

As such, homogenisation theory \cite{Pavliotis2008} and the iterated weak invariance principle (Theorem 1.1 \cite{kelly2017chaos}, see also \ref{KM2017}) will therefore apply to equations \eqref{skewfluidsystem}, \eqref{fastevals} resulting in convergence in distribution to a stochastic differential equation as $\varepsilon \rightarrow 0$. Thus the solution $x^\varepsilon_t$ of the ODE \eqref{skewfluidsystem} converges in distribution to a stochastic process $Z_t$.

Where $Z_t$ is the solution of an (Itô) stochastic differential equation, whose coefficients may be explicitly deduced when the following formulae converge:


\[\addtag\begin{aligned}\label{homogenisedfluid}
    \mathrm{d}Z_t &= U_t(Z_t)\mathrm{d}t + \sigma(Z_t) \mathrm{d}W_t\\
    U_t(Z_t) &:= \left\langle \frac{\partial \widetilde{g}_{t / \varepsilon}}{\partial \bar{g}_t}\overline{u}_t(\widetilde{g}^{-1}_{t/\varepsilon} (\cdot))\right\rangle_{\nu(\mathrm{d}\lambda)}(Z_t) + \left(\int_0^\infty \left\langle \phi_s(\lambda_0)^i \lambda^j_0\nabla v_i(\cdot) v_j(\cdot)\right\rangle_{\nu(\mathrm{d}\lambda)} \mathrm{d}s\right)(Z_t)\\
    \sigma \sigma^T(Z_t) &:= \left(\int_0^\infty \Big\langle \lambda^i_0\phi_s(\lambda_0)^j  v_i(\cdot) \otimes  v_j(\cdot) + \phi_s(\lambda_0)^i \lambda^j_0 v_i(\cdot) \otimes  v_j(\cdot) \Big\rangle_{\nu(\mathrm{d}\lambda)} \mathrm{d}s\right)(Z_t)
\end{aligned}\]

The coefficients in \eqref{homogenisedfluid} are calculated at $\varepsilon = 1$, using the flow map of the eigenvalue equation for $\lambda^\varepsilon_t$. Because of the invariant distribution assumption, $\lambda_0 \overset{d}{=} \lambda^\varepsilon_t := \phi_t(\lambda_0)$ we denote the integration of these random variables\footnote{Note that moments of $\lambda^\varepsilon_t$ are $t$-independent as a consequence of invariance, but mixed expressions of the type $\langle F(\lambda_0, \phi_t(\lambda_0)) \rangle_{\nu (\mathrm{d}\lambda)}$ need not be, see \cite{Pavliotis2008}[Ex. 11.7.1].} shorthand by $\nu(\mathrm{d}\lambda)$. 

The equation for $Z_t$ defines an Itô diffusion, it follows (up to regularity and completeness) from Theorem 9.2 in \cite{Kunita_1984} that this equation defines a stochastic flow $\Phi_t$ satisfying $\Phi_t(Z_0) = Z_t$ which can be identified as a stochastic diffeomorphism, one solving the SDE: 

$$\mathrm{d}\Phi_t \Phi^{-1}_t(x) = U_t(x)\mathrm{d}t + \sigma(x)\mathrm{d}W.$$

It remains to show that this is the form of the velocity constraint shown in the variational principle \eqref{abstractEP}, in fact, we will show the velocity is of the form seen in \cite{DHP2023} where an additional drift term may arise depending on the properties of the fast $\lambda$ dynamics. As remarked in \cite{cotter2017}, the resulting expression for $U$ contains both the average of the mean flow and an additional drift from the fluctuating displacement vector field identified in Cotter, Gottwald and Holm as $\int_0^\infty\langle \nabla_x \partial_t{\zeta}(\phi_t(\lambda),\overline{g}_t(x)) \, \zeta(\lambda, \overline{g}_t(x)) \rangle_{\nu(\mathrm{d}\lambda)} \mathrm{d}t$. We shall show that our analogous term (exchanging $\partial_t \zeta_{t/\varepsilon} \circ \overline{g}_t$ with $\dot{\widetilde{g}}_{t/\varepsilon}\widetilde{g}_{t/\varepsilon}^{-1} = \lambda^{k,\varepsilon}_t v_k$) representing the additional drift incurred by the fast map in \eqref{homogenisedfluid} is in fact a combination of the Itô-Stratonovich correction, which always appears (allowing us to change away from Itô noise), and the Wong-Zakai anomaly that appears conditionally. The remaining anomalous drift, after absorbing the Itô Stratonovich correction will lead to a stochastic equation for velocity vector field of the form:

\begin{equation}\mathrm{d}\Phi_t \Phi^{-1}_t(x) = \overline{U}_t(x)\mathrm{d}t +\frac{1}{2}\mathsf{s}^{ij}[\sigma_i, \sigma_j](x)\mathrm{d}t + \sigma(x) \circ \mathrm{d}W\label{WZvelocityeq}\end{equation}

Where we have denoted the drift with no correction term $\overline{U}_t = \left\langle \frac{\partial \tilde{g}_{t / \varepsilon}}{\partial \bar{g}_t}\overline{u}_t(\widetilde{g}^{-1}_{t/\varepsilon} (\cdot))\right\rangle_{\nu(\mathrm{d}\lambda)}$, this is precisely the velocity considered in \cite{DHP2023} when considering a coloured noise approximation of the SDE. Indeed, one can think of this approximation as a particular case of stochastic averaging when using the Ornstein-Uhlenbeck process as the fast $\lambda$ dynamics. When SALT is motivated from a homogenisation theory perspective this further supports investigation of this additional drift present from the model's derivation.  

As mentioned above, convergence to the limit SDE \eqref{WZvelocityeq} is valid under a large class \cite{284645bd-f971-3dfe-b374-dab0f8c0c677} of dynamical systems modelling the fast modes $\lambda$. Particular Green-Kubo type representations of the noise coefficient $\sigma$ were shown to hold in \cite{kelly2017chaos} under conditions that the $\lambda$ dynamics satisfy moment estimates. As discussed in \cite{engel2024nonlinear}, these representations agree with the classical formulae seen in \cite{Pavliotis2008} when the assumptions of both theories are met. A thorough analysis of the ergodic theory and convergence criterion for the validity of such representations is given in \cite{10.1214/21-AIHP1202, 10.1214/21-AIHP1203}. 

To relate $\sigma$ in terms of the eigenvectors $v_k$, we first, as in \cite{kelly2017chaos}, define a bilinear operator $\mathfrak{B}(v,w) = \int_0^\infty \langle v w \circ \phi_s \rangle_{\nu (\mathrm{d}\lambda )} \mathrm{d}s$. The convergence of $\mathfrak{B}$ holds under the assumptions of moment estimates of the type \cite{10.1214/21-AIHP1202, 10.1214/21-AIHP1203} which henceforth we assume are satified.

It follows we may write the correction term in \eqref{homogenisedfluid} in terms of $\mathfrak{B}$,


\begin{equation}
\left(\int_0^\infty \left\langle \phi_s(\lambda_0)^i \lambda^j_0\nabla v_i(\cdot) v_j(\cdot)\right\rangle_{\nu(\mathrm{d}\lambda)} \mathrm{d}s\right)^l = \mathfrak{B}(\lambda^j_0v^\alpha_j, \lambda^i_0\partial_\alpha v^l_i) = \mathfrak{B}(\lambda^j_0, \lambda^i_0)\partial_\alpha v^l_i v^\alpha_j := \left( E^{ij}\nabla v_i \cdot  v_j \right)^l \, .\label{StratWZform}\end{equation}

Equation \eqref{StratWZform} contains the Itô-Stratonovich correction. To verify this, one must explicitly calculate the diffusion coefficient $\sigma$ in \eqref{homogenisedfluid}, again in terms of $\mathfrak{B}$.



$$\begin{aligned}
\left(\sigma(x)\sigma(x)^T\right)^{\alpha \beta} &= \left(\int_0^\infty \Big\langle \lambda^i_0 \phi_s(\lambda_0)^j  v_i(\cdot) \otimes  v_j(\cdot) + \phi_s(\lambda_0)^i \lambda^j_0 v_i(\cdot) \otimes  v_j(\cdot) \Big\rangle_{\nu(\mathrm{d}\lambda)} \mathrm{d}s\right)^{\alpha \beta}\\
&= \mathfrak{B}(\lambda^i_0 v^\alpha_i, \lambda^j_0 v^\beta_j) + \mathfrak{B}(\lambda^j_0 v^\alpha_j, \lambda^i_0 v^\beta_i) = \mathfrak{B}(\lambda^i_0, \lambda^j_0)v^\alpha_iv^\beta_j + \mathfrak{B}(\lambda^j_0 , \lambda^i_0 )v^\alpha_jv^\beta_i \\ 
&= \left( E^{ji}v_i(x) \otimes v_j(x) + E^{ij}v_i(x)\otimes v_j(x) \right)^{\alpha \beta}\end{aligned}$$

We decompose the matrix $E$ into it's symmetric and antisymmetric parts. Let $E = M + \mathsf{s}^\prime$ with $M$ the symmetric part possessing a Cholesky decomposition\footnote{It can be shown that $\mathfrak{B}$ is positive semidefinite, and thus so is $M$, allowing the use of the Cholesky decompostion.} $M = DD^T$. 
It follows that 
\[ \sigma(x)\sigma(x)^T = 2M^{ij}v_i(x) \otimes v_j(x) = 2D^{ik}D^{jk}v_i(x) \otimes v_j(x) = \sqrt{2}D^{ik}v_i(x) \sqrt{2}D^{jk}v^T_j(x) \]

With the last equality following from the definition of the outer product\footnote{Note the indices $D^{jk}$ are not swapped in this transpose, as this is not a matrix vector product, rather it is a sum over a basis.} $x \otimes y = xy^T$.  It follows from the above calculation that $\sigma^{\alpha \beta}(x) = \sqrt{2}D^{i\beta}v^{\alpha}_i(x)$, where again we have denoted $v^{\alpha}_i(x)$ as the $\alpha$ component of the $i$-th vector. To verify this claim note that 
\[(\sigma \sigma^T)^{\alpha \beta} = \sigma^{\alpha k} \sigma^{\beta k} = \sqrt{2}D^{ik}v^{\alpha}_i \sqrt{2}D^{jk}v^{\beta}_j = 2 M^{ij}v^\alpha_iv^\beta_j = (2M^{ij}v_i \otimes v_j)^{\alpha \beta}.\]

It follows that \eqref{StratWZform} can be written as:

\begin{equation}
E^{ij}\nabla v_i \cdot  v_j = M^{ij}\nabla v_i \cdot  v_j + (\mathsf{s}^\prime)^{ij}\nabla v_i \cdot  v_j 
\label{Lévyplusstrat}\end{equation}

The symmetric part of \eqref{Lévyplusstrat} is identified as the Itô-Stratonovich correction term with the covariance matrix $D$. To see this, observe that 

\[(M^{ij}\nabla v_i \cdot  v_j)^l = D^{ik}D^{jk}\frac{\partial v^l_i}{\partial x^\alpha}v^\alpha_j = \frac{1}{2}\frac{\partial \sigma^{lk}}{\partial x^\alpha}\sigma^{\alpha k}\]

One recovers the higher dimension Itô-Stratonovich correction term \cite{Pavliotis_2014}. For the remaining term, observe that because of the antisymmetry of the $i,j$ indices due to the matrix $\mathsf{s}^\prime$ we can write the gradient term as a Lie bracket. Furthermore, one can write $\mathsf{s}^\prime$ in terms of the basis defined by $D$, as $\mathsf{s}^\prime = D \mathsf{s} D^T$. 
\[(\mathsf{s}^\prime)^{ij}\nabla v_i \cdot  v_j = D^{i \alpha}\mathsf{s}^{\alpha \beta}D^{j \beta} [v_i, v_j] = \frac{1}{2}\mathsf{s}^{\alpha \beta}[\sigma^{(\cdot)}_\alpha, \sigma^{(\cdot)}_\beta]\]

This is the form of the Wong-Zakai anomaly defined as the commutator of the columns of the diffusion matrix $\sigma$ and contracted with the linear perturbation to the Lévy area $\mathsf{s}$, a skew symmetric matrix. As a result, the velocity equation \eqref{WZvelocityeq} is obtained by decomposition of the full correction term, and changing to Stratonovich noise. 

\subsection{Rough paths theory intepretation of covariance and L\'evy area} We shall now describe how we interpret the matrices $D, \mathsf{s}$ and how they may arise in a physical context. 

The Green-Kubo type formula employed to calculate $\sigma$ can be though of as the autocorrelation tensor of the fast map $f_0 = \sum_{k=1}^m v_k(x_t) \lambda_t^{\varepsilon, k}$, this can be made precise in the sense of the Ergodic theorem \cite[Section 11.6]{Pavliotis2008} (see also the supplementary appendix \ref{sec:homog_formulas}) or the operator $\mathfrak{B}$ in Theorem \ref{KM2017} attributed to Kelly, Melbourne \cite{kelly2017chaos}.

From this viewpoint the matrix $M$ is the variance-covariance matrix of the resulting Brownian motion when one applies the weak invariance principle to calculate $\sigma$. The bilinear form $\mathcal{B}$ in this context is the correlations between the fast scale eigenvalues of the fast vector field. 

Note, that we can absorb the correlation matrix $D$ as part of the definition of the resulting diffusion tensor $\sigma$, using a standard Brownian motion that satisfies $\mathbb{E}[W^i_tW^j_t] = \delta^{ij}t$, this is the approach taken above. Alternatively we can define a correlated Brownian motion $B_t$ satisfying $\mathbb{E}[B^i_tB^j_t] = M^{ij}t$, with the noise coefficient $\sigma$ not containing any factors of $D$. In this context we define $\xi_i = \sqrt{2}v_i$, where we have the vector fields in the SALT Euler-Poincaré equations originating from the eigenvectors of the underlying fast flow (modulo factors of $\sqrt{2}$)\footnote{One may rescale either the definition of $\mathfrak{B}$ or the matrices $E = M + \mathsf{s}^\prime$ (c.f. \cite[Remark 1.4]{kelly2017chaos}) to eliminate these factors of $\sqrt{2}$. We shall follow the convention that $\mathfrak{B}$ agrees with formulae for the drift \& diffusion coefficients found in \cite{Pavliotis2008}.}.

It is important to check both these interpretations result in the same equations. Recall the fast dynamics are chosen such that they obey the iterated weak invariance principle \eqref{KM2014} \cite{kelly2016}[Thm. 9.1], that is, the approximation $W^\varepsilon_t := \varepsilon \int_0^{t/\varepsilon^2}\phi_s(\lambda) \mathrm{d}s$ has a canonical lift to a rough path via iterated integrals, such that they converge to the limiting rough path as follows:
\[(W^\varepsilon_t - W^\varepsilon_s , \mathbb{W}^\varepsilon) \rightarrow (DW_t - DW_s, D\otimes D(\mathbb{W}_{s,t}^{\text{Strat}}) +\mathsf{s}^\prime(t-s))\] 

Where we denote $D \otimes D ( \mathbb{W}_{s,t}^{\text{Strat}}) := D\mathbb{W}_{s,t}^{\text{Strat}} D^T$. One can either define $(B_t - B_s, \mathbb{B}^{\text{Strat}}_{s,t} + \mathsf{s}^\prime(t-s)) := (DW_t - DW_s, D\otimes D(\mathbb{W}_{s,t}^{\text{Strat}}) +\mathsf{s}^\prime(t-s))$ or instead work with $(W_t - W_s, \mathbb{W}_{s,t}^{\text{Strat}} + \mathsf{s}(t-s))$ and reinsert missing factors of $D$ in the coefficients of the homogenised SDE.

We now show that either approach is equivalent. Clearly it holds that $b(x,t) \circ \mathrm{d}B_t = b(x,t) \circ \mathrm{d} DW_t = b(x,t)D\circ \mathrm{d}W_t$ for any matrix $b(x,t)$. If we consider the noise term: \[\left(\sigma(x) \circ \mathrm{d}W_t\right)^{i} = \sigma^{ij}(x) \circ \mathrm{d}W^j_t := \sqrt{2}D^{\alpha j} v^i_\alpha(x) \circ \mathrm{d}W^j_t = \sqrt{2}v_\alpha^i(x) \circ \mathrm{d}D^{\alpha j}W^j_t := \left(\xi(x) \circ \mathrm{d}B_t\right)^i\]

The conclusion from this calculation is one can write the noise matrix $\sigma$ in terms of vectors $\xi_i$ seen in \eqref{abstractEP}. These $\xi_i$'s are, up to scaling, the slow eigenvectors that remain after homogenising their fast eigenvalues, \emph{provided} one is willing to work with a correlated Brownian motion (else we must include additional ``standard deviation" factors $D = \sqrt{M}$ should we require standard Brownian motion $W_t$ that satisfies $\mathbb{E}[W^i_tW^j_t] = \delta^{ij}t$). 

After showing the noise term is consistent we check the Itô-Stratonovich and Wong-Zakai correction terms are in agreement regardless of the representations of the limit of the rough driver $(W^\varepsilon, \mathbb{W}^\varepsilon)$. We make use of the quadratic covariation $\llbracket \cdot, \cdot\rrbracket_t$, and Lévy's characterisation $\llbracket W^i, W^j\rrbracket_t = \delta^{ij}t$ for components of a standard Brownian motion. It follows by bilinearity of $\llbracket \cdot,\cdot \rrbracket_t$, $B_t = DW_t, DD^T = M$ that for correlated Brownian motions, $\llbracket B^i,B^j \rrbracket_t = M^{ij}t$ and so:
\[ \begin{split}\frac{1}{2}\llbracket \sigma^{lk}(Z), W^k\rrbracket_t &= \frac{1}{2}t\frac{\partial \sigma^{lk}}{\partial x^\alpha}\sigma^{\alpha k} = \frac{1}{2}\delta^{mn}t\frac{\partial \sigma^{lm}}{\partial x^\alpha}\sigma^{\alpha n} = \frac{1}{2}\left\llbracket W^m, W^n\right \rrbracket_t\frac{\partial \sigma^{lm}}{\partial x^\alpha}\sigma^{\alpha n} = \frac{1}{2}\delta^{mn}D^{i m}D^{j n}t\frac{\partial \sqrt{2}v^{l}_i}{\partial x^\alpha}\sqrt{2} v^{\alpha}_j\\ &=  \frac{1}{2}D^{ik}D^{jk}t\frac{\partial \sqrt{2}v^l_i}{\partial x^\alpha}\sqrt{2}v^\alpha_j = \frac{1}{2}M^{ij}t\frac{\partial \sqrt{2}v^l_i}{\partial x^\alpha}\sqrt{2}v^\alpha_j =: \frac{1}{2}\left\llbracket B^i , B^j\right\rrbracket_t\frac{\partial \xi^l_i}{\partial x^\alpha}\xi^\alpha_j = \frac{1}{2}\llbracket\xi^l_k(Z), B^k\rrbracket_t .\end{split}\] 
Likewise, for the Wong-Zakai correction term, perturbing the Lévy stochastic area of $\mathbb{B}$ by $\mathsf{s}^\prime(t-s)$ produces a correction term of $\frac{1}{2}(\mathsf{s}^\prime)^{ij}[\xi_i, \xi_j]$, by definition this is the same as perturbing $\mathbb{W}$ by $\mathsf{s}$ with a noise coefficient of $\sigma$, since:

\[\frac{1}{2}(\mathsf{s}^\prime)^{ij}[\xi_i, \xi_j] = D^{i \alpha}\mathsf{s}^{\alpha \beta}D^{j \beta} [v_i, v_j] = \frac{1}{2}\mathsf{s}^{\alpha \beta}[\sqrt{2}D^{i \alpha}v_i^{(\cdot)}, \sqrt{2}D^{j \beta} v^{(\cdot)}_j]= \frac{1}{2}\mathsf{s}^{\alpha \beta}[\sigma^{(\cdot)}_\alpha, \sigma^{(\cdot)}_\beta].\]

Thus either representation is a well defined and consistent definition of a stochastic differential equation with alternative choices of noise and diffusion coefficients.

\section{The Wong-Zakai anomaly for 2D incompressible flows}\label{sec:continuumeuler}

In section \ref{sec:2}, we defined the stochastic Euler-Poincar\'e equations and accompanying variational principle, showing that the Euler equations are a special case. The noise and correction terms were motivated through homogenisation in section \ref{sec:homogenisation}. We now specialise these constructions to the specific case of ideal 2-D fluid dynamics which is the key example studied in this paper. 

The vorticity formulation of the stochastic Euler equation \eqref{EulerwithLevy} can be compactly written as,
\begin{equation}
    \mathrm{d} \omega_t = -\pounds_{\mathrm{d}g_t g^{-1}_t} \omega_t = -\nabla \cdot (\omega_t \mathrm{d}g_t g^{-1}_t) = -\mathrm{d}g_tg^{-1}_t \cdot \nabla \omega_t \,, \quad \mathrm{d}g_tg^{-1}_t = u_t \mathrm{d}t + \frac{1}{2}\mathsf{s}^{ij}[\xi_i, \xi_j] \mathrm{d}t + \xi \circ \mathrm{d} W_t \, ,
\label{vorticity}\end{equation}

where the exterior derivative (curl, for $1$-forms in Euclidean $2$-space) is taken and $\bold{d}u_t^\flat =: \omega_t = \operatorname{curl}(u_t)$ is a $2$-form. 
A key property used is that the exterior derivative commutes with the stochastic differential and the Lie derivative, it follows from equation \eqref{vorticity} that $\omega$ is transported/advected along the stochastic flow and that $\omega_t = (g_t)_* \omega_0$, \eqref{vorticity} is obtained from this relation via an application of the Kunita-Itô-Wentzell formula \cite{de2020implications}. 

The divergence free property in the plane allows a stream function representation of the vector fields $u$ and $\xi_i$. We define the following stream functions $u_t = -\nabla^\perp\psi_t$, $\xi_i = -\nabla^\perp \psi_i$ where $\nabla^\perp = (-\partial_y, \partial_x)$. Under these assumptions, we claim the case of nonzero Wong-Zakai is also given by a potential flow:
\[\frac{1}{2}\mathsf{s}^{ij}[\xi_i, \xi_j] = \frac{1}{2}\mathsf{s}^{ij}[\nabla^\perp \psi_i, \nabla^\perp \psi_j] = -\frac{1}{2}\mathsf{s}^{ij} \nabla^\perp J(\psi_i, \psi_j) \quad \text{where}, \quad J(\psi_i, \psi_j) = -\nabla \psi_i \cdot \nabla^\perp \psi_j.\]

We may prove this in generality. For any symplectic manifold $(M, \theta)$, given $f \in C^\infty(M)$ one defines the Hamiltonian vector field $X_f$ as the unique vector field satisfying $\bold{d}f = \theta(X_f, \cdot)$, with $\bold{d}$ denoting exterior derivative. The commutator of two Hamiltonian vector fields remains Hamiltonian, with: 

$$-X_h = [X_f, X_g], \,\, \text{and}, \,\, h := \{f, g\} := \theta(X_f, X_g).$$

The symplectic form $\theta: TM \times TM \rightarrow \mathbb{R}$, which is a non degenerate, bilinear, induces isomorphisms 

\[\flat_\theta : TM \rightarrow T^*M, \quad X \mapsto X^{\flat_\theta} := \theta(X, \cdot), \quad \sharp_\theta = (\flat_\theta)^{-1} : T^*M \rightarrow TM.\]


This allows a coordinate-free definition of symplectic gradient $\nabla^\perp := \sharp_\theta \bold{d}$, note that no Riemannian metric is required. One may reformulate the operator $f \mapsto X_f$ in terms of $\nabla^\perp$, 

\[\nabla^\perp f = \sharp_\theta \bold{d}f \iff \omega(\nabla^\perp f, \cdot ) = \omega(\sharp_\theta \bold{d}f, \cdot) =: \flat_\theta \sharp_\theta \bold{d} f = \bold{d} f, \quad \text{thus,} \,\, X_f = \nabla^\perp f\]

We may then conclude: 


$$[\nabla^\perp f, \nabla^\perp g] = [X_f, X_g] = -X_h = -\nabla^\perp h = -\nabla^\perp \{f,g\}.$$



The identity above is sufficient to specify the Wong-Zakai anomaly for a symplectic manifold. In the case of $\mathbb{R}^2$, one has the presence of an addditional Riemmanian structure, one may further write $[\nabla^\perp \psi_i, \nabla^\perp \psi_j] = - \nabla^\perp \left( -\nabla \psi_i \cdot \nabla^\perp \psi_j \right)$.


More generally, an even dimensional manifold equipped with both a Riemannian metric $\boldsymbol{g}$ and symplectic form $\theta$ is said to have compatible metric and symplectic structures if there exists a smooth tensor field $\mathbb{J} : TM \rightarrow TM$ satisfying $\mathbb{J}^2 = -\operatorname{id}_{TM}$ such that the following holds: 

\[\boldsymbol{g}(X,Y) = \theta(X, \mathbb{J} Y), \quad \boldsymbol{g}(\mathbb{J}X,\mathbb{J}Y) = \boldsymbol{g}(X,Y), \quad \theta(\mathbb{J}X,\mathbb{J}Y) = \theta(X,Y). \]

The tensor $\mathbb{J}$ is known as an almost complex structure. Compatibility of $\boldsymbol{g}$ and $\theta$ ensures $\sharp_\theta = \mathbb{J} \circ \sharp_{\boldsymbol{g}}$. The skew gradient then can be more recognisably expressed in terms of the perpendicular of the metric gradient, 
\[\nabla^\perp := \sharp_\theta \bold{d} =  \mathbb{J} \circ \sharp_{\boldsymbol{g}} \bold{d} =: \mathbb{J} \circ \nabla.\]
One can then recognise $\mathbb{R}^2$ equipped with Euclidean metric and $\mathbb{J} = \perp$ are compatible to recover the initially shown expression,
\begin{align*}
     \{f,g\} &= \theta(X_f, X_g) = \boldsymbol{g}(\mathbb{J}X_f, X_g), \\
    \text{in } \mathbb{R}^2 \text{ we have,} \quad &= (\nabla^\perp f)^\perp \cdot \nabla^\perp g = -\nabla f \cdot \nabla^\perp g .
\end{align*}



Therefore, the Wong Zakai anomaly may also be represented as the pairwise Jacobian determinants of stochastic stream functions. This is one of many instances where the Wong-Zakai anomaly respects the geometric structure of equations it appears in. It follows that the total stochastic velocity is given by the stochastic stream function

$$\mathrm{d}g_tg^{-1}_t = -\nabla^\perp \Psi_t \circ \mathrm{d}S_t = -\nabla^\perp \left(\psi_t\mathrm{d}t +  \frac{1}{2}\mathsf{s}^{ij}J(\psi_i, \psi_j) \mathrm{d}t + \psi_k \circ \mathrm{d}W^k_t \right)$$

The notation $\circ \mathrm{d}S_t$ denotes compatibility of the semimartingale $\nabla^\perp \Psi_t$ with the semimartingale $S_t = (t, W^1_t, \ldots, W^K_t)$, see Definition 2.3 \cite{crisan2022variational}. In rough paths terms this can also be thought of as constructing a new path $(S, \mathbb{S})$ by defining $S_t = (t, W_t)$ and defining the iterated integrals involving $W_t$ only via $\mathbb{W}_{0,t}$, and remaining cross integrals of $t, W^i_t$ as standard Riemann-Stieltjes integrals. This allows understanding the expression $-\mathrm{d}g_tg^{-1}_t\cdot \nabla \omega_t = -\nabla^\perp \Psi_t \circ \mathrm{d}S_t \cdot \nabla \omega_t =: \{\omega_t, \Psi_t \circ \mathrm{d}S_t \}$ as a Poisson bracket on functionals of the vorticity. By introducing a Hamiltonian such that $\frac{\delta h}{\delta \omega_t} \circ \mathrm{d}S_t = \Psi_t \circ \mathrm{d}S_t$ this bracket, and equation \eqref{vorticity} can be interpreted as Lie-Poisson \cite{holm1998eulerS}.

In order to derive the corresponding Hamiltonian we relate the vorticity and stream functions through $\omega_t = \operatorname{curl}(-\nabla^\perp \psi_t) = -\Delta \psi_t$.  The stream function can then be recovered from the vorticity itself via the Biot-Savart law, using the Green’s function of the Laplacian to compute $(-\Delta)^{-1}\omega_t$ with appropriate boundary
conditions. Consequently the stochastic Hamiltonian for \eqref{vorticity} that produces the correct Lie-Poisson bracket is then given by:

\begin{equation}
    h(\omega_t) \circ \mathrm{d}S_t = \int_M \frac{1}{2}\omega_t (-\Delta)^{-1}\omega_t \mathrm{d}t + \frac{1}{2}\mathsf{s}^{ij} J(\psi_i, \psi_j)\omega_t \mathrm{d}t + \psi_k(x) \omega_t \circ \mathrm{d}W^k_t \mathrm{d}^2x.
\label{hamiltonian}\end{equation}

Note, as in \cite{DHP2023} the Hamiltonian is modified by the Wong-Zakai anomaly. Since $\Delta$ is self adjoint (under the assumption of boundary conditions) it can be shown that $\frac{\delta h}{\delta \omega_t} \circ \mathrm{d}S_t = \psi_0\mathrm{d}t + \frac{1}{2}\mathsf{s}^{ij}J(\psi_i, \psi_j) \mathrm{d}t + \psi_k \circ \mathrm{d}W^k_t = \Psi_t \circ \mathrm{d}S_t$, which recovers the Lie-Poisson equation.

\subsection{Stochastic Point Vortex Dynamics in Two Dimensions}\label{sec:SPV2D}


Stochastic Advection by Lie Transport assumes that the Lagrangian position satisfies the following stochastic particle trajectory mapping \cite{crisan2022variational, holm2015variational}, 
\begin{align}
\b x( \b X,t) = \b x(\b X,0) + \int_{0}^{t} \b u(\b x(\b X,s),s)\mathrm{d}s + \sum_{i=1}^{m}\int_{0}^{t} \b \xi_{i}(\b x(\b X,s)) \circ \mathrm{d} W_s^i, \label{eq:particle traj map}
\end{align}
evolving the initial label $\b X = (X,Y)$ to current configuration $\b x= (x,y)$.
The deterministic stream function $\psi$ is related to vorticity by a differential relationship, solveable using a Green's function as follows 
\begin{align}
\psi(\b x,t)= \int_{\mathbb{R}^2} G(\b x-\b x') \omega( \b x',t) d \b x'.
\end{align}
The negative skew gradient $\b u =-\nabla^{\perp}\psi$ relates the velocity $\b u$ to the vorticity $\omega$ by the kernel $K$
\begin{align}
 \b u(\b x,t)= \int_{\mathbb{R}^2}  K(\b x-\b x') \omega( \b x',t) d\b x', \label{eq:biot savart}
\end{align}
this relationship is known as the Biot-Savart law. Given the Hamiltonian \eqref{hamiltonian} and the Biot-Savart law above one can derive the point vortex model for the Euler equations via the the assumption that all the vorticity is initially concentrated at finitely many points in the domain.

This concentration remains for further times (as a result of Kelvin's theorem, see \eqref{eq: Kelvin circulation theorem} below) at these points following the flow such that, 
\begin{equation}
\omega(\b x,t) = \sum_{\alpha=1}^n \Gamma_{\alpha}\delta (\b x - \b x_{\alpha}(t))
\,.
\label{pointvortex-momap}
\end{equation}

In geometric mechanics terminology, one thinks of this as defining a (singular) momentum map \cite{MARSDEN1983} from phase space points $\lbrace (x_{\alpha} \b, \Gamma_{\alpha} y_{\alpha}) \rbrace_{\alpha = 1}^n$  to the vorticity, which acts as a momentum variable for the Lie-Poisson system. 

The Lagrangian position of an arbitrary position must necessarily obey \eqref{eq:particle traj map}, thus we may consider the evolution of a point vortex $\b x_{\alpha} = (x_{\alpha}(t),y_{\alpha}(t))^T\in \mathbb{R}^2$ with initial position $\b X_{\alpha} = (X_{\alpha}(0),Y_{\alpha}(0))^T\in \mathbb{R}^2$ evolving by the stochastic particle trajectory mapping 
\begin{equation}
\begin{aligned}
    \b x_{\alpha}( \b X_{\alpha},t) = \b x_{\alpha}(\b X_{\alpha},0) &+ \int_{0}^{t} \b u(\b x_{\alpha}(\b X_{\alpha},s),\lbrace \b x_{\beta}(\b X_{\beta},t) \rbrace_{\forall \beta\neq \alpha},s)\mathrm{d}s\\ & + \sum_{i=1}^{m}\int_{0}^{t} \b \xi_{i}(\b x_{\alpha}(\b X_{\alpha},s),\lbrace \b x_{\beta}(\b X_{\beta},t) \rbrace_{\forall \beta}) \circ \mathrm{d} W^i_s \, ,
\end{aligned}
\end{equation}
where we have explicitly denoted the additional functional dependence of $\b \xi_{i}=\b \xi_{i}(\b x_{\alpha}(s),\lbrace \b x_{\beta}\rbrace_{\beta})$ on other points in the domain, similar in analogy to how the deterministic velocity field is implicitly dependent on the other non local points in the domain through the Biot-Savart Law. 
The in-compressible vector fields $\b \xi_{i}= -\nabla^{\perp} \psi_{i}$ are assumed to arise from the skew gradient of a stream function $\psi_{i}$, and are integrated against the components of the $m$ dimensional Brownian motion $W_t$.


Finally, we show the stochastic Kelvin-Noether circulation theorem \cite{holm2015variational} implies the conservation of vortex strength $\mathrm{d} \Gamma_{\alpha}=0$, $\forall \alpha$. Given a closed \emph{stochastically advected} loop $C(t)$ satisfying $\mathrm{d}C = -\pounds_{u \mathrm{d}t + \xi * \mathrm{d}W_t}C$, let $\Omega(t)$ denote the contained region of the loop such that $C(t)= \partial \Omega(t)$. Then the circulation of the loop is the sum of the strengths of the point vortices contained within the loop, and the dynamic definition of the Lie derivative or the pullback by the stochastic flow map, reveals the ``stochastic conservation law"
\begin{align}
\sum_{\alpha\in \Omega(0)}\Gamma_{\alpha} = \int_{\Omega(0)}\omega(\b X, 0) \bold{d}X\wedge \bold{d}Y =\oint_{C(0)}\b u( \b X,0)\cdot \bold{d}\b X = \oint_{C(t)}\b u(\b x, t)\cdot \bold{d}\b x = \int_{\Omega(t)}\omega(\b x, t) \bold{d}x\wedge \bold{d}y  = \sum_{\alpha\in \Omega(t)}\Gamma_{\alpha}. \label{eq: Kelvin circulation theorem}
\end{align}

Where the second equality follows from Stoke's theorem and the fact that $\bold{d}\b u^\flat = \omega$, and the Kelvin-Noether theorem was used in third equality. 

Equation \eqref{eq: Kelvin circulation theorem} holds for all times $t\in \mathbb{R}_{+}$ and all closed, stochastically advected loops $C(t)$, and implies the strength of each point vortex $\Gamma_{\alpha}$ is constant in time (take a loop around each point vortex). 


Putting the ansatz \cref{pointvortex-momap}, the Kelvin theorem \cref{eq: Kelvin circulation theorem}, the Hamiltonian \eqref{hamiltonian}, the Biot-Savart law \cref{eq:biot savart}, and the SALT assumption \cref{eq:particle traj map} together turns the infinite dimensional stochastic fluid dynamics problem into a finite dimensional $\mathbb{R}^{2n}$ system of stochastic ordinary differential equations, described below:
\begin{align}
&\b x_{\alpha}(t) = \b x_{\alpha}(0) + \int_{0}^{t} \b u(\b x_{\alpha}(s),s)\mathrm{d}s + \sum_{i=1}^{m}\int_{0}^{t}\b \xi_{i}(\b x_{\alpha}(s),\lbrace \b x_{\beta}\rbrace_{\forall \beta })\ast \mathrm{d}W^{i}_t, \quad \forall \alpha \in \lbrace 1,...,N \rbrace, \\
&\b u(\b x_{\alpha}(s),s) = \sum_{\beta=1,\beta\neq \alpha}^{n} \Gamma_{\beta} K( \b x_{\alpha}(s) - \b x_{\beta}(s)), \quad \forall \alpha \in \lbrace 1,...,N \rbrace, \\
&\b \xi_{i}(\b x_{\alpha}(s),\lbrace \b x_{\beta}\rbrace_{\beta \neq \alpha}) =-(\nabla^{\perp}\psi_{i})(\b x_{\alpha}(s),\lbrace \b x_{\beta}\rbrace_{\beta\neq \alpha}),\quad \forall \alpha \in \lbrace 1,...,N \rbrace.
\end{align}


\begin{remark}
Here we are differentiating $\psi_{i}(\b x, \lbrace x_{\beta}\rbrace_{\forall \beta})$ in the Eulerian coordinates $\b x$ (ignoring the solution dependent parameter $\lbrace \b x_{\beta}\rbrace_{\forall \beta}$) and substituting in $\b x = \b x_{\alpha}$, to get the velocity at position $\alpha$. In the finite dimensional system it is common to define differentiation with respect to the system variables, for example $\nabla_{\alpha} = (\partial_{x_{\alpha}},\partial_{y_{\alpha}})$. 
It is here we will make the following observation that the gradient of a function evaluated at a position vector, is equivalent to evaluating at a position vector, and taking the gradient with respect to that position vector and reevaluating at this point, $
\nabla f(\b x;\lbrace \b x_{\beta}\rbrace_{\forall \beta})|_{\b x = \b x_{\alpha}} = \nabla_{\alpha}f(\b x_{\alpha})|_{\b x_{\alpha}=\b x_{\alpha}}.$ This allows us to discuss the finite dimensional Hamiltonian structure as follows.
\end{remark}

The canonical Hamilton form, is modified by the stochastic integration in the following way 
\begin{align}
\Gamma_{\alpha} \mathrm{d} \b x_{\alpha} = -\nabla^{\perp}_{\alpha} \left(H \mathrm{d}t + \sum_{\alpha=1}^{n}\sum_{i=1}^{m}\Gamma_{\alpha} \psi_{i}(\b x_{\alpha};\lbrace \b x_{\beta}\rbrace_{\forall \beta })\ast \mathrm{d}W^{i}_t. \right), \quad \forall \alpha\in \lbrace 1,...,n\rbrace, \quad \text{for} \quad \nabla^{\perp}_{\alpha} := (-\partial_{y_{\alpha}},\partial_{x_{\alpha}}),
\label{eq:PV_SDE}\end{align}
 where $H$ is the deterministic Kirchhoff Hamiltonian 
 \begin{align}
     H = \sum_{\alpha,\beta =1,\alpha\neq\beta}^{n}\Gamma_{\alpha}\Gamma_{\beta}G(\b x_{\alpha} - \b x_{\beta}),\label{eq:KirchoffHamiltonian}
 \end{align}
describing the deterministic energy around the point vortices. By considering the time differential $\mathrm{d}$ of some function of state variables $F(t; x_1, ..., x_{\alpha}, ..., x_n,y_1,...,y_{\alpha},...,y_n)$, 
\begin{align}
\mathrm{d}F &=\partial_t F + \sum_{\alpha=1}^{n} \partial_{x_{\alpha}}F \mathrm{d}x_{\alpha} + \partial_{y_{\alpha}}F \mathrm{d}y_{\alpha},\\
&= \partial_t F + \sum_{\alpha=1}^{n}\nabla_{\alpha} F \cdot  - \frac{1}{\Gamma_{\alpha}}\nabla_{\alpha}^{\perp}H \mathrm{d}t + \sum_{\alpha=1}^{n}\nabla_{\alpha} F \cdot - \frac{1}{\Gamma_{\alpha}}\nabla_{\alpha}^{\perp} \sum_{i=1}^{m}\psi_{i} 
\ast \mathrm{d}W^i_t, \\
&=
\partial_t F + \lbrace F,H \rbrace \mathrm{d}t +\sum_{i=1}^{m} \lbrace F, \psi_{i} \rbrace 
\ast \mathrm{d}W^i_t,
\end{align}
we can identify the Poisson bracket as
\begin{align}
\lbrace F,H \rbrace = \sum_{\alpha=1}^{n} \frac{1}{\Gamma_{\alpha}} \bigg(\frac{\partial F}{\partial x_{\alpha}} \frac{\partial H}{\partial y_{\alpha}}-\frac{\partial F}{\partial y_{\alpha}} \frac{\partial H}{\partial x_{\alpha}} \bigg) = \sum_{\alpha=1}^{n} \frac{1}{\Gamma_{\alpha}}\nabla_{\alpha} F \cdot -\nabla_{\alpha}^{\perp} H. \label{eq:poisson bracket}
\end{align}

The Hamiltonian of the point vortex system is modified (or even potentially lost) depending on the choice of stochastic integration, for example in the particular case of Stratonovich noise with a Wong-Zakai anomaly, \begin{align}
H \mapsto \left(H + \frac{1}{2}\sum_{1\leq i,j\leq m}\sum_{\alpha=1}^{n}\Gamma_{\alpha} \mathsf{s}^{ij} J(\psi_i(\b x_{\alpha};\lbrace \b x_{\beta}\rbrace_{\forall \beta }), \psi_j(\b x_{\alpha};\lbrace \b x_{\beta}\rbrace_{\forall \beta }))  \right) \mathrm{d}t + \sum_{\alpha=1}^{n}\sum_{k=1}^{m}\Gamma_{\alpha} \psi_{k}(\b x_{\alpha};\lbrace \b x_{\beta}\rbrace_{\forall \beta })\circ \mathrm{d}W^{k}_t. 
\end{align}
\subsubsection{Deterministic theory and conserved quantities}\label{Example:Deterministic Point Vortex System}

Our analysis of the stochastic system with noises of different types shall be bench-marked against the well studied deterministic dynamics; we remind the reader of its properties and behaviour. Without noise the classical point vortex problem is described by the deterministic equation for particle trajectories
\begin{align}
    \mathrm{d} \b x_{\alpha}( \b X,t) = \b u( \b x_{\alpha}(\b X,t),t)\mathrm{d}t;\quad \b x_{\alpha}(\b X,0) = \b X.
\end{align} 
In the specific case of two dimensional Euler on $\mathbb{R}^2$, the differential relationship relating the vorticity to streamfunction is the translationally invariant ($\psi =(-\Delta)^{-1}\omega$) Poisson's equation with Green's function $G$ described by a radial Newtonian potential,  
\begin{align}
    G(\b x-\b x') = \frac{-1}{2\pi} \log (\|\b x-\b x'\|_2)= \frac{-1}{4\pi} \log (\|\b x-\b x'\|_2^2). \label{eq:greens}
\end{align}

Where $\|\cdot\|_2$ denotes the $\ell^2$ norm. The Biot–Savart law relating velocity to vorticity is defined through the kernel $K$ 
\begin{align}
    K(\b x - \b x') = \frac{-((y-y'),x-x'))^{T}}{2\pi \|\b x-\b x'\|_2^{2}}, \label{eq:Euler BS Kernel}
\end{align}
arising from the skew gradient of the Greens function, this Kernel exhibits blow up when $\b x=\b x'$. 
This Biot-Savart law specifies the velocity and defines the following Hamiltonian system, 
\begin{align}
\frac{d x_{\alpha}}{dt} &= \frac{1}{2\pi}\sum_{\beta=1,\beta \neq \alpha}^n \frac{-\Gamma_{\beta}(y_{\alpha}-y_{\beta})}{(y_{\alpha}-y_{\beta})^2+(x_{\alpha}-x_{\beta})^2}, \\
\frac{d y_{\alpha}}{dt} &= \frac{1}{2\pi}\sum_{\beta=1,\beta \neq \alpha}^n \frac{\Gamma_{\beta}(x_{\alpha}-x_{\beta})}{(y_{\alpha}-y_{\beta})^2+(x_{\alpha}-x_{\beta})^2}. \label{eq:det velocity}
\end{align}
 The conserved quantities associated with in-variance of the Hamiltonian \cref{eq:KirchoffHamiltonian} under translation and rotation are
\begin{align}
T_x = \sum_{\alpha=1}^{n}\Gamma_{\alpha} x_{\alpha}, \quad T_y = \sum_{\alpha=1}^{n}\Gamma_{\alpha} y_{\alpha}, \quad R = \frac{1}{2}\sum_{\alpha=1}^{n}\Gamma_{\alpha}(x_{\alpha}^2+y_{\alpha}^2).
\end{align}
these are the conservation of linear and angular impulse. Directly from the definition of the Poisson bracket \cref{eq:poisson bracket}
one can verify the conservation laws $\partial_t (T_x,T_y,R,H)^T=0$ using the properties
$\lbrace T_x,H \rbrace$, $\lbrace T_y,H \rbrace$, $\lbrace R,H \rbrace$, $\lbrace H,H \rbrace = 0$ after directly computing
$\lbrace x_{\beta},y_{\gamma} \rbrace =  \frac{1}{\Gamma_{\beta}}\delta_{\beta,\gamma}
$, $\lbrace x_{\beta},x_{\gamma} \rbrace = 0$, $\lbrace y_{\beta},y_{\gamma} \rbrace = 0$. The quantities $T_x^2+T_y^2, R,H$ are independent integrals in involution $
\lbrace T_x^2+T_y^2, R \rbrace , \lbrace T_x^2+T_y^2, H \rbrace,  \lbrace R, H \rbrace = 0$ giving integrability for the $3$ point vortex problem by Liouville’s theorem. These properties have all been well established in the literature \cite{aref2007point,poincare1893theorie,grobli1877specielle}.

Specific to the case of $n = 3$ unit strength $\Gamma_i=1$ point vortices initially located on the roots of the unit circle $z^3 = -1$, $z\in \mathbb{C}$, it has been proven that the point vortices remain in a rotating equilateral triangle with the area $A_{\alpha \beta \gamma}$, angles $\alpha_{\alpha \beta \gamma}$, center of vorticity $\b x_c = (\sum_{\alpha=1}^{n}\Gamma_{\alpha})^{-1} (T_x , T_y)^T$ and intervorticial distance $l_{\alpha \beta}$ conserved. The center of vorticity is equivalent to the center of mass, $(x_{c},y_{c}) := (\frac{1}{n}\sum_{\alpha=1}^{n} x_{\alpha}, \frac{1}{n} \sum_{\alpha=1}^{n}y_{\alpha}) = (0,0)$, and the canonical coordinates $(p_{\alpha},q_{\alpha}):=(x_{\alpha},\Gamma_{\alpha} y_{\alpha})$ are the Cartesian coordinates of the point vortices. 

It is in this classically understood setting we set up our case study into the effect of the Wong Zakai anomaly, Numerical Lévy area, and other choices of stochastic integration.

\subsection{Case study : Stochastically stable and unstable configurations of three unit strength point vortices} \label{sec: case study}


We consider the evolution of $n=3$ point vorticies of unit vortex strength $\Gamma_{\alpha}=1$ initially located on the roots of the unit circle $z^3 = -1$, namely $\b X_1,\b X_2,\b X_3  = (-1,0)^T, (1/2,\sqrt{3}/2)^T,(1/2,-\sqrt{3}/2)^T$. Whose evolution is described by the  SDE
\begin{align}
\mathrm{d} \b x_{\alpha}(\b X,t) = \b u(\b x_{\alpha}(\b X,t))\mathrm{d}t + \sum_{i=1}^{2} \b \xi_{i}(\b x_{\alpha}(\b X,t),\b x_{c}(\lbrace \b x_{\beta}(\b X,t)\rbrace_{\forall \beta }))  \ast \mathrm{d}W_t^i, \quad  \b x_{\alpha}(\b X,0) = \b X_\alpha, \label{eq:general case}
\end{align}

where $\b x_{\alpha} = (x_{\alpha},y_{\alpha})\in \mathbb{R}^2$, denotes the position of the particle $\alpha$ for $\alpha \in \lbrace 1,2,3\rbrace$ after some time $t$ so that we have 6 degrees of freedom. As before, the deterministic contribution of velocity is described by the Biot-Savat law whilst the vector fields $\b \xi_{1},\b \xi_{2}$ (with associated stream functions $\psi_1, \psi_2$), are integrated against the two dimensional driving signal $ (W_t^{1}, W_t^2)$ enhanced as a Itó, Stratonovich or Statonovich with Wong-Zakai anomaly rough path. 



We define a rotational stream function $\psi_1$ about the center of vorticity $(x_c,x_c)$, such that $\forall (x,y)\in \mathbb{R}^2$ 
\begin{align}
\psi_{1}(\b x, \b x_{c}(\lbrace \b x_{\beta}\rbrace_{\forall \beta}); A, r) = A\exp \left( -\frac{r}{2} \|\bx - \b x_c(\{\b x_{\beta} \}_{\forall \beta})\|^2_2 \right) ,\quad A,r > 0.
\end{align}

The above definition of $\psi_1$ is defined for an arbitrary value $\b x$, but is typically evaluated at each point vortex $\b x = \b x_{\alpha}=(x_{\alpha},y_{\alpha})^T$ in the numerical method eg:\eqref{eq:PV_SDE}. The value of $\psi_1$ has a dependence on the other particles by the center of vorticity $\b x_c = (\sum_{\alpha=1}^{n}\Gamma_{\alpha})^{-1} (T_x , T_y)^T$. We denote the additional dependence of the function $\psi_1$ on the other particles as $\lbrace \b x_{\beta}\rbrace_{\forall \beta}$. 
The rotation is defined such that in the frame of reference of the center of vorticity (denoted by a $\widehat{\cdot}$ symbol) the resulting vector field for \emph{arbitrary Eulerian coordinates} rotate the point vortices by $\widehat{\b \xi}_1 :=(y,-x)^{T}  A r\exp \left(- \frac{r}{2} \|\b x\|^2_2 \right)$, and leaves the center of vorticity invariant. 




We define a translation stream function $\psi_{2}$
\begin{align}
\psi_{2}(\b x; \b a) &= \b x \cdot \b a , \quad \b a = (a,b)^T \in \mathbb{R}^2
\end{align}
giving a constant velocity field ${\b \xi}_2 = (- b,a)^{T}$ which moves points and center of vorticity in a uniform direction. All positions relative to the center of vorticity remain unchanged under $\b \xi_2$. For a diagram of both vector fields, see \cref{fig:velocity fields} below. 

\begin{figure}[H]
  \centering
  \begin{minipage}[b]{0.45\textwidth}
    \includegraphics[width=\textwidth]{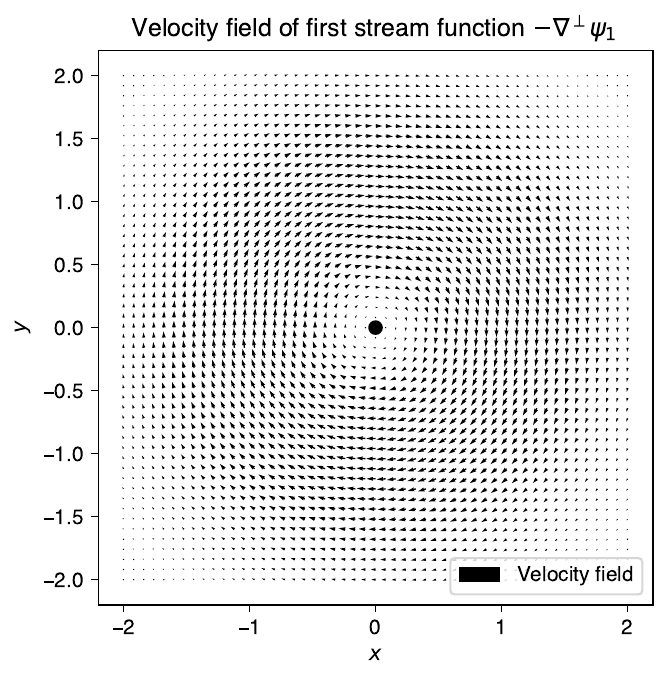}
    \caption{The radial vector field associated to noise that induces rotation around $\b x_c = 0$ times a multiple of $\mathrm{d}W^1_t$ to the velocity field. $r = 1, A = \frac{1}{2}$, giving $\widehat{\b \xi}_1 = (y, -x) \frac{1}{2}\exp(-\frac{1}{2}\|\b x\|)$.  }
  \end{minipage}
  \hfill
  \begin{minipage}[b]{0.45\textwidth}
    \includegraphics[width=\textwidth]{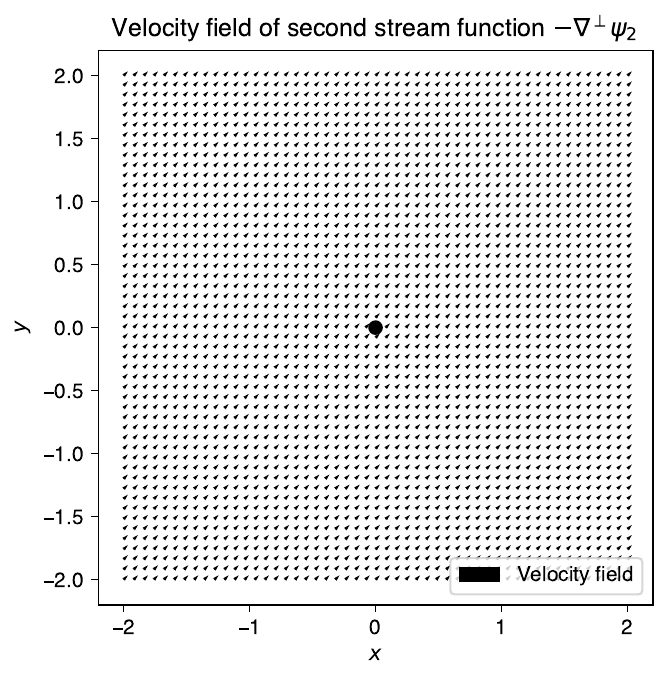}
    \caption{The constant vector field associated to noise that translates in a multiple of $\mathrm{d}W^2_t$ by the direction $\b \xi_2 = (-b,a) = (a,b)^\perp$ shown for $a=1, b=-1$.}\label{fig:velocity fields}
  \end{minipage}
\end{figure}

It is worth noting that any set of translational vector fields have vanishing Lie brackets (as constant fields commute), 
\begin{align}
\left[\nabla^{\perp}\psi_{2}(\bx ; \b a_1),\nabla^{\perp}\psi_{2}(\b x; \b a_2)\right] = 0, \quad \forall \b a_1, \b a_2 \in \mathbb{R}^2,
\end{align}
and any set of rotational vector fields about the origin will also have vanishing Lie brackets (as both are functions of $\|\b x - \b x_c\|$). 
\begin{align}
\left[\nabla^{\perp}\psi_{1}(\b x; A_1, r_1),\nabla^{\perp}\psi_{1}(\b x; A_2, r_2)\right] = 0, \quad \forall r_1, r_2, A_1, A_2 > 0. 
\end{align}

Such constraints prevent being able to examine the Wong-Zakai anomaly with $\psi_1$ or $\psi_2$ alone. In the case where there is a rotation and a translation the commutator
\begin{align}
   [-\nabla^{\perp} \psi_1,-\nabla^{\perp} \psi_2] = A r e^{-  \frac{r}{2} \left(x^{2} + y^{2}\right)} \left(- r b x y + a \left(r y^{2} - 1\right) ,  - r a x y + b \left(r x^{2} - 1\right)\right)^T \label{WZcommutator}
\end{align}
does not vanish, see \cref{fig:Wong_Zakai drift Lie Braket}. 
As per \cref{sec:continuumeuler}, the commutator of incompressible vector-fields has a stream function relationship $[\b \xi_1,\b \xi_2] = \nabla^{\perp} J(\psi_1, \psi_2)$, and for our particular vector fields, we have the corresponding streamfunction relationship
\begin{align}
\psi_{WZ} = J(\psi_1, \psi_2) = A r \left(- a (y-y_c) + b (x-x_c)\right) e^{-  r \left((x-x_c)^{2} + (y-y_c)^{2}\right)/2}
\end{align}

One notes that $\psi_{WZ}$ is not invariant under rotation and translation.



\subsubsection{Stable Stratonovich configurations }\label{Example:Stratonovich Point Vortex System}

 
Our key comparison of the rough point vortex system \eqref{eq:general case} will be in reference/comparison to being driven by Stratonovich Brownian motion, $\mathrm{d}\boldsymbol{Z} = \circ \mathrm{d}W$. This choice of noise is expected to arise from any limit (homogenisation being one such example) where the noise vector fields $\b \xi_i$ commute \cite{ikeda2014stochastic}, or in the special case where the approximation of the Weiner process (or homogenisation, respectively) produces no perturbation to the Lévy area, $\mathsf{s} = 0$. The evolution of the point vortices is given by
\begin{align}
\mathrm{d} \b x_{\alpha}(\b X,t) = \b u(\b x_{\alpha}(\b X,t),t)\mathrm{d}t + \b \sum_{i=1}^2\b \xi_{i}(\b x_{\alpha}(\b X,t),\lbrace \b x_{\beta}(\b X,t)\rbrace_{\forall \beta })\circ \mathrm{d}W^{i}_t;\quad \b x_{\alpha}(\b X,0) = \b X.\label{stratparticle}
\end{align}

First we note that in the case of no translation, that is, $\psi_2 = 0$ we preserve the conserved quantities $T_x, T_y, R, H$. The conservation of $T_x, T_y, R$ (in the sense they obey the SDE $\mathrm{d}\{T_x, T_y, R\} = 0$) follows from the fact that the entire stochastic Hamiltonian $\Psi \circ \mathrm{d}S_t= \psi_0 \mathrm{d}t + \psi_1 \circ \mathrm{d}W$ is $SE(2)$ invariant and applying Noether's theorem for stochastic systems \cite{takao2020stochastic}[Prop. 2.18]. To verify the conservation of the deterministic energy $H$ we directly compute using the stochastic analogue of \eqref{eq:poisson bracket}:


$$\mathrm{d}H = \{H, \Psi \circ \mathrm{d}S_t\} = \{ H, H \mathrm{d}t\} + \{ H, \psi_1 \circ \mathrm{d} W^1_t\} = \sum_{\alpha,\beta =1, \alpha\neq \beta}^n  Ar \frac{\Gamma_{\alpha} \Gamma_{\beta}}{2 \pi } \frac{y_{\alpha} x_{\beta} - x_{\alpha} y_{\beta}}{\|\bx_{\alpha} - \bx_{\beta}\|^2} e^{-\frac{r}{2}\|\b x_\alpha\|^2} \circ \mathrm{d}W^1_t$$ 

Note, that in the absence of translation we have $\bx_c = 0$, when $n=3$ one also has each intervorticial distance $l_{\alpha \beta}:=\|\b x_{\alpha}-\b x_{\beta}\|_2$ and distance from the origin $\|\bx_\alpha-\b 0\|$ conserved and index independent with $\Gamma_\alpha \equiv 1$. It follows that \[ \frac{Are^{-\frac{x_{\alpha}^2 + y_\alpha^2}{2}}}{2\pi \|\bx_{\alpha} - \bx_{\beta}\|^2} = C\] for some $C$ independent of $\alpha, \beta$ and since $ \sum_{\alpha,\beta = 1,\alpha\neq \beta}^n y_\alpha x_\beta - x_\alpha y_\beta = 0$, we conclude $\mathrm{d}H = 0$. A variant of this argument also holds in the case that $\psi_2 \neq 0$ and $\b x_c$ is not the origin, see \cref{sec:supplementary calcs}.

For the case of $\psi_2 \neq 0$ we make use of an argument attributed to Takao and Holm (unpublished masters thesis 2016), proving that for a constant vector field $\b \xi_{i}(\b x) = (a_{i},b_{i})$,  driving the flow of the point vortex problem, the system is in a moving frame with the same deterministic conservation laws. The conservation laws and invariance must be viewed in a stochastically moving coordinate frame. Let $F \in \{H, T_x, T_y, R\}$ denote a quantity that was conserved in the deterministic system, the evolution in the case of \eqref{stratparticle} is given as:

$$\mathrm{d} F = \{F , \Psi \circ \mathrm{d}S_t\} = \{F, H\mathrm{d}t\} + \{F, \psi_1 \circ \mathrm{d}W^1_t\} + \{F, \psi_2 \circ \mathrm{d}W^2_t\}$$

As expected $\{F, H\mathrm{d}t\} = 0$ continues to hold from the deterministic theory. Although $\psi_2 \neq 0$ means $\b x_c$ is no longer identically zero, we still have $\{F, \psi_1 \circ \mathrm{d}W^1_t\} = 0$ as $SE(2)$ symmetry of $\psi_1$ still holds, since $\|R \b x - \b x^\prime - (R \b x_c - \b x^\prime) \| = \|\b x - \b x_c \|$, noting $\b x_c$ moves with any transformation of the point vortices. See \cref{sec:supplementary calcs} for a direct calculation.

It is important to note that $\psi_2$ has no $SE(2)$ symmetry, and we do not expect zero Poisson brackets from its contribution. Indeed it can be immediately calculated from the definitions that $\{T_x, \psi_2 \circ \mathrm{d}W^2_t \} = bn \circ \mathrm{d}W^2_t$ and $\{T_y, \psi_2 \circ \mathrm{d}W^2_t \} = - na \circ \mathrm{d}W^2_t$ to which we conclude the evolution of the quantities $T_x, T_y$:




$$\mathrm{d} T_x -b n \circ \mathrm{d}W^2_t = 0, \quad \mathrm{d} T_y + a n \circ \mathrm{d}W^2_t = 0.$$ 

By changing into the moving co-ordinate frame $(\widehat{x}_\alpha,\widehat{y}_\beta) := (x_{\alpha}(0)-b W^2_t, y_{\alpha}(0) + a W^2_t)$, reveals a conservation law $\mathrm{d} (T_{\widehat{x}},T_{\widehat{y}}) = (0,0)$ in the moving frame. Furthermore, we have: 
$$\{R, \psi_2 \circ \mathrm{d}W^2_t \} = -\sum_{\alpha=1}^n \nabla_\alpha \left(\frac{1}{2}\sum_{\beta=1}^n x_\beta^2 + y_\beta^2 \right) \cdot \nabla_\alpha^\perp \psi_2 \circ \mathrm{d}W^2_t = \sum_{\beta=1}^n (b x_\beta - a y_\beta) \circ \mathrm{d}W^2_t = (bnx_c -any_c  ) \circ \mathrm{d}W^2_t$$

Resulting in a stochastic evolution equation of $\mathrm{d} R = (bn x_c - an y_c ) \circ \mathrm{d}W^2_t$. By switching to the stochastic coordinate frame we enforce $(\widehat{x}_c, \widehat{y}_c) = 0$ and correspondingly $\widehat{R} := R(\cdot - \widehat{\b x}_c)$ is constant in the moving frame. Similarly for the deterministic energy with translation we have:
$$\mathrm{d}H = \{H, \Psi \circ \mathrm{d}S_t\} =\{ H, \psi_2 \circ \mathrm{d} W^2_t\} = -\sum_{\alpha = 1}^n\nabla_\alpha H \cdot \nabla^\perp \psi_2 \circ \mathrm{d}W^2_t = \sum_{\alpha=1}^n\left(b \frac{\partial H}{\partial x_\alpha}-a \frac{\partial H}{\partial y_\alpha}\right) \circ \mathrm{d} W^2_t = -(b \mathrm{d}T_y + a \mathrm{d}T_x)$$ 

One substitutes the previously calculated equations for $T_x, T_y$ and finds $H$ conserved.

$$ -b \mathrm{d}T_y - a \mathrm{d}T_x = (abn - abn) \circ \mathrm{d}W^2_t = 0$$





To summarise, we verified the following evolution equations for $T_x,T_y,R,H$ 

\begin{align}
    \mathrm{d}T_x &=  bn \circ \mathrm{d}W^2_t\\
    \mathrm{d}T_y &=  -an \circ \mathrm{d}W^2_t \\
    \mathrm{d}R &= \left(any_c  - bnx_c \right) \circ \mathrm{d}W^2_t \\
    \mathrm{d}H &=  b \mathrm{d}T_y + a \mathrm{d}T_x = 0 
\end{align}

By defining the frame $(\widehat{x}_\alpha,\widehat{y}_\alpha) := (x_{\alpha}(0)-b W^2_t, y_{\alpha}(0) + a W^2_t)$ we correspondingly define stochastically translated analogues. Note that $\Gamma_\alpha \equiv const, \forall \alpha$ is necessary to define these as conserved (see \cref{sec:supplementary calcs}) in the sense of a stochastic Hamiltonian system. 

\begin{align}
T_{\widehat{x}} = \sum_{\alpha=1}^{n} \widehat{x}_{\alpha}, \quad T_{\widehat{y}} = \sum_{\alpha=1}^{n} \widehat{y}_{\alpha},\quad \widehat{\b x}_c = \frac{1}{n}(T_{\widehat{x}} , T_{\widehat{y}})^T, \quad \widehat{R} = \frac{1}{2}\sum_{\alpha=1}^{n}(\widehat{x}_{\alpha} - \widehat{x}_{c})^2+(\widehat{y}_{\alpha} - \widehat{y}_{c})^2.
\end{align}

We recover an interesting example in which the Stratonovich system of point vortices once again conserves quantities associated with invariance under translations and $\b x_c $ centered rotations. Furthermore, by a geometric argument one gets preservation of area, angle and intervorticial distance for the point vortex initially located on $z^3 = 1$ in the $(\widehat{x},\widehat{y})$ frame.

\begin{remark}
    Similar arguments are expected to hold for this configuration when the driving signal is a strongly geometric rough path in light of \cite{crisan2022variational}. A particular case of interest is FBM, which possesses a canonical geometric lift to an $\alpha\in (1/3,1/2)$ Hölder rough path. This also has a motivation via homogenisation theory, (see \ref{sec:WIPs}) and the case of Hurst index $0.4$ is tested in \cref{Method:FBM}. 
\end{remark}

\subsubsection{Pseudostable Itô configurations} \label{Example:Itô Point Vortex System}
In this section, we examine the properties of Itô noise in comparison to the Stratonovich system. Note that a variational approach such as SALT requires Stratonovich noise to derive equations to perform integration by parts and satisfy the Kunita-Itô-Wentzell formula \cite{de2020implications}, however one can ameliorate this difficulty to derive Itô equations by writing in terms of Stratonovich plus a correction term, whos inclusion we shall investigate the dynamical consequences.

We also remark that Itô noise is one possible result of a homogenisation procedure (see \cite{Pavliotis2008}[Sec 11.7.6]), albeit not the one considered in this paper in \ref{sec:homogenisation}.
\begin{align}
\mathrm{d}\b x_{\alpha}( \b X,t) = \b u( \b x_{\alpha}(\b X,t),t)\mathrm{d}t + \sum_{i=1}^2 \b \xi_{k}(\b x_{\alpha}(\b X,t),\lbrace \b x_{\beta}(\b X,t)\rbrace_{\forall \beta }) \mathrm{d}W^k_t; \quad \b x_{\alpha}(\b X,0) = \b X.
\label{ItoPV}\end{align}
Interpreting this in terms of Stratonovich noise has the effect of removing a drift, via the Itô-Stratonovich correction 
$\frac{1}{2}\sum_{k=1}^{2} (\b \xi_k \cdot \nabla )\b \xi_k = \frac{1}{2}\b \xi_1 \cdot \nabla \b \xi_1$, where the second term in the sum vanishes because $\b {\xi}_2$ is constant. A physical interpretion of the Itô-Stratonovich correction in the context of fluid dynamics is examined in \cite{Holm2020a}. Constant Itô noise has been considered as early as 1973 \cite{chorin1973numerical}, and has the interpretation of adding viscosity by allowing the point vortices to have a random walk, simulating  diffusion, and giving the corresponding Fokker-Plank equation a viscosity term \cite{majda2002vorticity}. 

In the particular construction we have chosen, the drift associated with Itô Stratonovich correction acts in the frame of reference of the center of vorticity and causes point vortices to drift further apart in a radial manner see, \cref{fig:Itô-stratonovich}, such that when $\b x_c = 0$, the Itô Stratonovich correction for the vector field's takes the following form
\begin{align}
\frac{1}{2}\sum_{k=1}^{2} \b \xi_k \cdot \nabla \b \xi_k = 
\frac{1}{2} A^{2} r^{2} e^{- r \left(x^{2} + y^{2}\right)}(-x,-y)^T,\label{eq:ito-strat-vectorfield}
\end{align}
at the level of the continuum. One remarks that the Itô-Stratonovich correction of the point vortex system is the regular Itô-Stratonovich correction of vector fields evaluated and summed over all point vortices (see \cref{Sec:Lévy Area contribution in higher order Stratonovich Integrators}). Consequently, one could conjecture the loss of this particular Itô-Statonovich correction (of the form \cref{eq:ito-strat-vectorfield}), causes points to drift further away from the origin (center of vorticity), preserving the angle between the three point vortices, but not the area of the inscribed equilateral triangle. Furthermore \cref{eq:ito-strat-vectorfield} has a non zero divergence, resulting in a Itô-Stratonovich correction without a stream-function or (Stratonovich) Hamiltonian representation. The divergence is explicitly calculated to be,




$$\nabla \cdot \left(\frac{1}{2}\sum_{i=1}^{2} (\b \xi_k \cdot \nabla )\b \xi_k \right) = A^{2} r^{2} \left(r x^{2} + r y^{2} - 1\right) e^{- r \left(x^{2} + y^{2}\right)} \neq 0.$$
\subsubsection{Unstable Wong-Zakai anomaly configuration}\label{Example:Type I }

In this section we examine the point vortex system in the presence of a Wong-Zakai anomaly. Motivated from the generic case that a smooth approximation of Brownian motion (one such case being homogenisation) converges to a geometric lift with $\mathsf{s} \neq 0$, in terms of SDE this corresponds to a correction term appearing distinct from the Itô-Stratonovich correction. 

We are interested in investigating the following system 
\begin{align}
d \b x_{\alpha} =\b u(\b x_{\alpha},t)\mathrm{d}t + \sum_{k=1}^{2} \b{\xi}_k (\b x_{\alpha},\lbrace \b x_{\beta} \rbrace_{\forall \beta})\circ \mathrm{d}W^k_t + \frac{1}{2} \sum_{1\leq i, j \leq 2}\mathsf{s}^{ij}[\b{\xi}_i,\b{\xi}_j] (\b x_{\alpha},\lbrace \b x_{\beta} \rbrace_{\forall \beta})\mathrm{d}t.
\end{align}
The bracket computed in \cref{WZcommutator} admits a stream function representation, but this stream function fails to exhibit the same symmetries of $\psi_1$ or $\psi_2$. As a result, we no longer expect angle, inscribed area, and other quantities to be conserved.
Regardless of whether the modified system possess novel conserved quantities, or none at all, the ``energy surface" formed by intersecting level sets of conserved quantities is entirely different to that of the Stratonovich or deterministic systems, and one conjectures the dynamical behaviour can be drastically different.

\begin{figure*}
\centering
\begin{subfigure}[t]{0.45\textwidth}
        \centering
        \includegraphics[width=.90\textwidth]{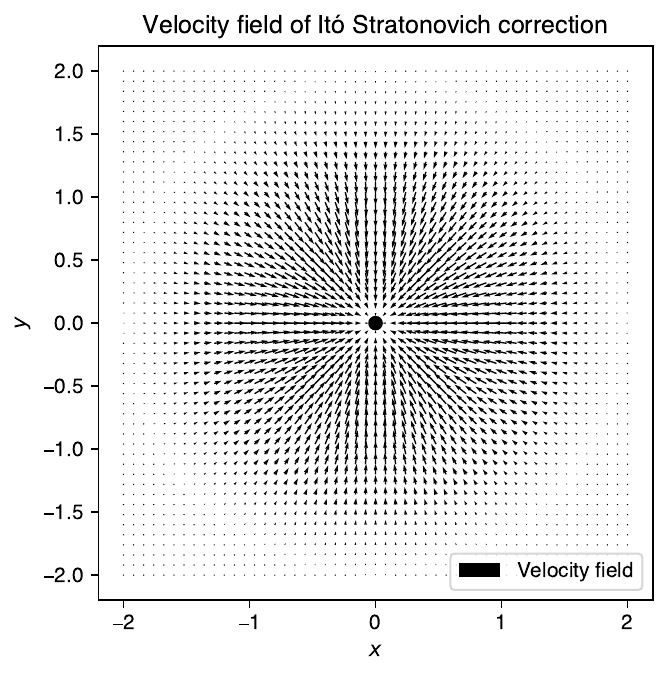}
        \caption{The vector field associated with converting Itô to Stratonovich, around the center of vorticity initially at $(0,0)$. }
\label{fig:Itô-stratonovich}
\end{subfigure}
\hfill
\begin{subfigure}[t]{0.44\textwidth}
        \centering
       \includegraphics[width=.95\textwidth]{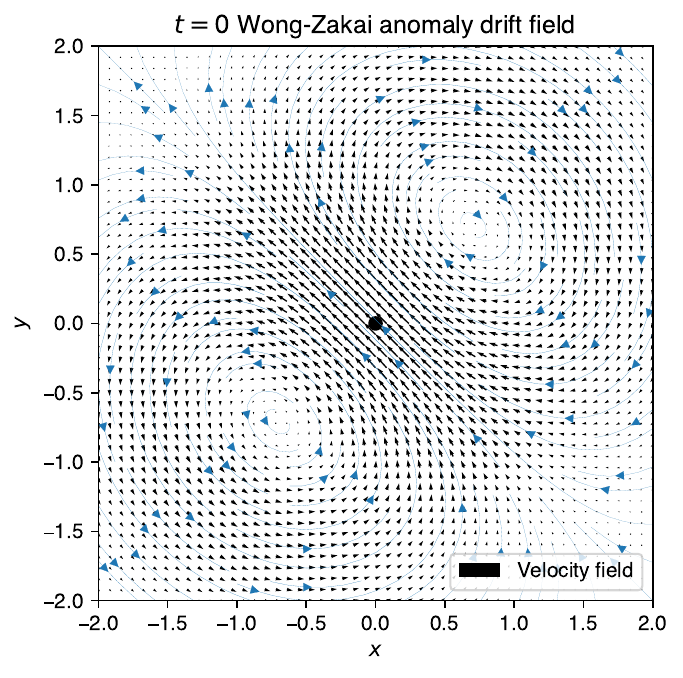}
\caption{Shown is the velocity field and streamlines of the Lie Bracket of vector fields at time $0$ where the center of vorticity is initially located at $(0,0)$. }
\label{fig:Wong_Zakai drift Lie Braket}
\end{subfigure}
\caption{ The subtraction of the Itô-Stratonovich correction drift term \cref{fig:Itô-stratonovich} (relative to the Stratonovich noise \cref{stratparticle}) caused the Itô point vortices to remain circular but spread out into a larger circle. The appearance of the Wong-Zakai drift field \cref{fig:Wong_Zakai drift Lie Braket} and the location of the point vortices within it, caused additional dynamical behaviour. }
\end{figure*}

\Cref{fig:Wong_Zakai drift Lie Braket} visualises the resulting Wong-Zakai drift field in the frame of the center of vorticity.
Here, for the numerical example we will assume 
\begin{align}
\mathsf{s} = 
\begin{bmatrix}
    0,&1\\
    -1,&0
\end{bmatrix},
\end{align}
given. In practice the magnitude of this peturbation is not known, the true size of the perturbation could be assimilated from data/observations or a priori knowledge about the system. We pick a normalised value of 1, for which the effects are visible. 
\paragraph{Wong-Zakai anomaly from a pure area process}
\label{Example: Type II}

It is indeed possible to consider the Wong-Zakai anomaly drift in the absence of any noise, this can be interpreted as taking the driving path $\boldsymbol{Z}_{t,s} = (0, \mathsf{s}(t-s))$ to be a pure area process \cite{Hairer2020, friz2010multidimensional}. Such examples can be constructed as limits of oscillatory smooth paths of the form  $Z^1_n(t) = 1/n \cos(2 n^2 t)$, $Z^2_n(t) = 1/n \sin(2 n^2 t)$ converging uniformly to 0 pathwise, but in the rough paths topology to 
\begin{align}
(Z^1_t,Z^2_t) =\left(\begin{bmatrix}
0\\
0\\
\end{bmatrix}, 
\begin{bmatrix}
0,t\\
-t,0\\
\end{bmatrix}\right)\in \mathbb{R}^2\otimes \mathfrak{s}\mathfrak
{o}(2).
\end{align}
See \cite{friz2010multidimensional}[Ex. 0.1-0.2-0.4]. 
The associated point vortex dynamics driven by time and $\boldsymbol{Z}$ can be identified with the ODE solution of 
\begin{align}
\mathrm{d}\b x = (\b u(\b x,t) + [\b{\xi}_1, \b{\xi}_2])\mathrm{d}t,
\end{align}\label{type2}
we consider this the footprint of the previously defined Stratonovich noise with Wong-Zakai anomaly drift, as to examine the drift in the absence of fluctuations. This is a entirely different deterministic Hamiltonian system $\Gamma_{\alpha} \partial_{t} \b x_{\alpha} = -\nabla^{\perp}H_{WZD}$ with Hamiltonian given by \begin{align}
    H_{WZD} = H + \sum_{\alpha=1}^{n} \Gamma_{\alpha}\psi_{WZ}(\b x_{\alpha};\lbrace \b x_{\beta} \rbrace_{\forall \beta}),
\end{align} in which we expect different dynamical behaviour. This Hamiltonian system does not posess the same translational and rotational in-variances as the deterministic system and we do not have reason to expect the circular motion or triangular configuration to be preserved. Nevertheless, the conservation of energy $\partial_t H_{WZD}= \lbrace H_{WZD},H_{WZD}\rbrace = 0$ follows, and it is plausible that this new Hamiltonian system may have other conserved quantities.


\section{Numerical investigation}\label{Numerical investigation}

In this section we numerically investigate the changes in dynamical behaviour caused by the different interpretations of noise, in reference to the expected dynamics from the calculations in sections \ref{Example:Stratonovich Point Vortex System}, \ref{Example:Itô Point Vortex System} and \ref{Example: Type II}. We first outline the algorithms and discretisation procedure in \cref{sec: numerics}. The results of the simulation are presented in \cref{sec:results} 

\subsection{Numerical Methods}\label{sec: numerics}


The point vortex system is converted into a state valued SDE system 
\begin{align}
\mathrm{d} \b q_t = \b f (\b q_t)\mathrm{d}t + G(\b q_t)\circ \mathrm{d} W_t,
\end{align}
where the state variable $\b q = (x_1,...,x_{\alpha},...,x_n ,y_1,...,y_{\alpha},...,y_n)^T\in \mathbb{R}^{2n}=\mathbb{R}^{d}$ consist of the stacked components of the location of point vortices. Where $\b f(\b q):\mathbb{R}^{2n}\rightarrow \mathbb{R}^{2n}$, $G(\b q):\mathbb{R}^{2n} \rightarrow \mathbb{R}^{2n\times m}$ and $ W_t$ is a m-dimensional Brownian motion, are constructed as in \cref{Sec:Lévy Area contribution in higher order Stratonovich Integrators}, with the exception that the velocity uses the regularised velocity
$$\b u(\b x_{\alpha}) = \sum_{\beta =1, \beta \neq \alpha}^{n} \Gamma_{\beta}K_{\delta_{1}}(\b x_{\alpha}-\b x_{\beta}).$$
Whose construction relies on vortex blob regularisation of the Greens function similar to that of Krasney 1986 \cite{krasny1986study}. However the Euler Kernel $K$ (\cref{eq:Euler BS Kernel}) is adaptively mollified as follows,
\begin{align}
    K_{\delta_{1}}(\b x - \b x') = \frac{(-(y-y'),x-x'))^{T}}{2\pi (\|\b x-\b x'\|_2^{2} + \delta_{1}^2)} \mathbb{I} ( \|\b x-\b x'\|_2^{2} \leq \delta_{1}^2) +  K(\b x - \b x')\mathbb{I} ( \|\b x-\b x'\|_2^{2} > \delta_{1}^2),
\end{align}
where $\mathbb{I}$ is an indicator function, and $\delta_1 \in \mathbb{R}^{+}$ is a fixed (small) number. This prevents division by zero (blow up) by changing to the regularised kernel when points get within a fixed distance of one another, otherwise the unregularised kernel is used. 

The timestepping method follows Chorin's \cite{chorin1973numerical} use of a high order explicit Runge-Kutta integration scheme. However we require convergence to several interpretations of noise, to do so we will use an explicit additive Runge-Kutta (ARK) scheme where the drift and diffusion are attributed their own Butcher tableau, see \cite{ruemelin1982numerical} for order conditions on the coefficients. We define the intermediary $s$ substages $ \b Q^k \in \mathbb{R}^{d}$ by
\begin{align}
\b Q^{k} =  \b q^{n} + \sum_{j=1}^{s} a_{k,j} \b f (\b Q^{k})\Delta t  + \sum_{j=1}^{s} \tilde{a}_{k,j}\sum_{i=1}^{m} \b g_i (\b Q^{k})  \Delta W^i, \quad \forall k\in \lbrace 1,...,s \rbrace, \label{eq: Semi-Discrete}
\end{align}
and the final update step by
\begin{align}
\b q^{n+1} = \b q^n + \sum_{k=1}^s b_k \b f(\b Q^k)\Delta t + \sum_{k=1}^s \tilde{b}_k \sum_{i=1}^{m} \b g_{i} (\b Q^k)  \Delta W^i.
\end{align}
The coefficients $\b c,\tilde{\b c} \in\mathbb{R}^s$, $A,\tilde{A} 
\in\mathbb{R}^{s\times s}$, $\b b,\tilde{\b b} \in\mathbb{R}^s$ will displayed with the following double Butcher tableau format
\begin{tabular}{c|c}
$\b c$ & $A$ \\
\hline & $\b b^T$
\end{tabular},
\begin{tabular}{c|c}
$\tilde{\b{c}}$ & $\tilde{A}$ \\
\hline & $\tilde{\b b}^T$
\end{tabular}
$=$
\begin{tabular}{c|cccc}
$c_1$ & $a_{11}$ & $a_{12}$ & $\ldots$ & $a_{1 s}$ \\
$c_2$ & $a_{21}$ & $a_{22}$ & $\ldots$ & $a_{2 s}$ \\
$\vdots$ & $\vdots$ & $\vdots$ & $\ddots$ & $\vdots$ \\
$c_s$ & $a_{s 1}$ & $a_{s 2}$ & $\ldots$ & $a_{s s}$ \\
\hline & $b_1$ & $b_2$ & $\ldots$ & $b_s$
\end{tabular},
\begin{tabular}{c|cccc}
$\tilde{c}_1$ & $\tilde{a}_{11}$ & $\tilde{a}_{12}$ & $\ldots$ & $\tilde{a}_{1 s}$ \\
$\tilde{c}_2$ & $\tilde{a}_{21}$ & $\tilde{a}_{22}$ & $\ldots$ & $\tilde{a}_{2 s}$ \\
$\vdots$ & $\vdots$ & $\vdots$ & $\ddots$ & $\vdots$ \\
$\tilde{c}_s$ & $\tilde{a}_{s 1}$ & $\tilde{a}_{s 2}$ & $\ldots$ & $\tilde{a}_{s s}$ \\
\hline & $\tilde{b}_1$ & $\tilde{b}_2$ & $\ldots$ & $\tilde{b}_s$
\end{tabular}
and uniquely determine the ARK method.
\begin{method}[Deterministic]\label{Method:Deterministic} We discretise the deterministic point vortex system, using RK4 method, 
\begin{align}
\begin{tabular}{c|c c c c}
0&0&0&0&0 \\
1/2&1/2&0 &0&0\\
1/2&0&1/2&0&0\\
1&0&0&1&0\\
\hline
&1/6&1/3&1/3&1/6\\
\end{tabular},
\begin{tabular}{c|c c c c}
0&0&0&0&0 \\
0&0&0&0&0 \\
0&0&0&0&0 \\
0&0&0&0&0 \\
\hline
&0&0&0&0 \\
\end{tabular}.
\end{align}
\end{method}

\begin{method}[Stratonovich method 1]\label{Method:Stratonovich integration scheme} We discretise the Stratonovich point vortex system introduced in \cref{Example:Stratonovich Point Vortex System}, using a RK4 scheme where the right hand side accepts both the drift and diffusion, described by the RK4-RK4, double Butcher tableau,
\begin{align}
\begin{tabular}{c|c c c c}
0&0&0&0&0 \\
1/2&1/2&0 &0&0\\
1/2&0&1/2&0&0\\
1&0&0&1&0\\
\hline
&1/6&1/3&1/3&1/6\\
\end{tabular},
\begin{tabular}{c|c c c c}
0&0&0&0&0 \\
1/2&1/2&0 &0&0\\
1/2&0&1/2&0&0\\
1&0&0&1&0\\
\hline
&1/6&1/3&1/3&1/6\\
\end{tabular}.
\end{align}    
\end{method}

\begin{method}[Itô system]\label{Method:type ito}
    For the discretisation of the Itô point particle system introduced in \cref{Example:Itô Point Vortex System}, we use an additive RK method, where RK4 does the deterministic drift and an Euler scheme treats the stochastic diffusion term. This RK4-EM approach is described by the double Butcher Tableau,
\begin{align}
\begin{tabular}{c|c c c c}
0&0&0&0&0 \\
1/2&1/2&0 &0&0\\
1/2&0&1/2&0&0\\
1&0&0&1&0\\
\hline
&1/6&1/3&1/3&1/6\\
\end{tabular},
\begin{tabular}{c|c c c c}
0&0&0&0&0 \\
1&1&0 &0&0\\
0&0&0&0&0 \\
0&0&0&0&0 \\
\hline
&1&0&0&0\\
\end{tabular}.
\end{align}rrR
\end{method}

\begin{method}[Type I]\label{Method:type I}
For the simulation of the Stratonovich and Wong-Zakai anomaly drift point particle system introduced in \cref{Example:Type I }, we absorb the Wong Zakai anomaly drift  into the drift term as follows $$
\b f(\b q) = \big( (\b u,\b v)^{T} + 1/2([\b \xi_1,\b \xi_2]^x,[\b \xi_1,\b \xi_2]^y)^{T} \big) (\b q),$$ and use the Stratonovich integrator \cref{Method:Stratonovich integration scheme}. Where the superscript $x,y$ denotes the $x$ and $y$ components of the vector field respectively, and $[\cdot,\cdot]$ denotes the regular commutator of vector fields. 
\end{method}

\begin{method}[Type II]\label{Method:type II}
For the simulation of the 
deterministic ``Area process" driven point particle system introduced in \cref{Example:Type I }, we absorb the Wong-Zakai anomaly terms into the drift term and use the deterministic \cref{Method:Deterministic}. 
\end{method}

\begin{method}[Stratonovich method 2 (Lévy Area)]\label{Method:NLA}
Denoted by the RK4-RK4-p-NLA, scheme, we use the existing Stratonovich RK4-RK4 \cref{Method:Stratonovich integration scheme}, and add on an approximation of the Lévy area. To do so, we use Kloeden Platen and Wright's truncated Fourier series expansion of a Brownian Bridge process up to $K$ terms (see \cite{kloeden1992approximation}, and \cref{sec:Numerical simulation of the Lévy area}), and use the direct evaluation of the commutator of vector fields at point vortex positions. 
\end{method}

\begin{method}[Canonical geometric lift of FBM] \cite{friz2010multidimensional} \label{Method:FBM}
We use the existing \cref{Method:Stratonovich integration scheme}, and replace the normally distributed increments by increments of realisations of fractional Brownian motion, generated using the Davies-Harte algorithm \cite{davies1987tests}. 
\end{method}

 Due to the Wagner-Platen / Stratonovich-Taylor expansion, it is understood \cref{Method:Stratonovich integration scheme} will capture the Itô-Stratonovich correction and converge to the Stratonovich system with strong order $1/2$, and \cref{Method:type ito} will converge to the Itô system with strong order $1/2$, this is explicitly shown as a subcase of the work in \cite{ruemelin1982numerical}. \Cref{Method:FBM} when driven with $\alpha$-Hölder signals such as FBM with Hurst parameter between $1/3$ and $1/2$, is consistent, as the order of the entire scheme is greater than or equal to three \cite{redmann2022runge}. (see \ref{Li2022} for one possible motivation for this generalisation, and \cite{redmann2022runge} for the strong (local and global) convergence rates). The Wong-Zakai anomaly is captured as a drift and numerically treated with 4'th order accuracy. When approximating Lévy's Stochastic area for a Weiner process in the context of the numerical method, we use a truncated Fourier series of a Brownian bridge process with enough terms for strong order 1 convergence \cite{kloeden1992approximation}, additional details \cref{sec:Numerical simulation of the Lévy area}. \Cref{Method:NLA} is technically higher order than the other methods in this paper, approximating the Lévy area such that the convergence of the scheme is strong order 1. 
 \Cref{Method:Stratonovich integration scheme} captures other higher order (symmetric) terms in the Stratonovich Taylor expansion such that if the commutator vanishes the scheme becomes strong order 2 for Brownian motion. If the commutator does not vanish then the scheme is strong order $1/2$ when approximating Brownian motion due to the missing Lévy area correction terms, see \cref{Sec:Lévy Area contribution in higher order Stratonovich Integrators}.



\begin{figure}[H]
\centering
\begin{subfigure}[t]{0.325\textwidth}
\centering
\includegraphics[width=.95\textwidth]{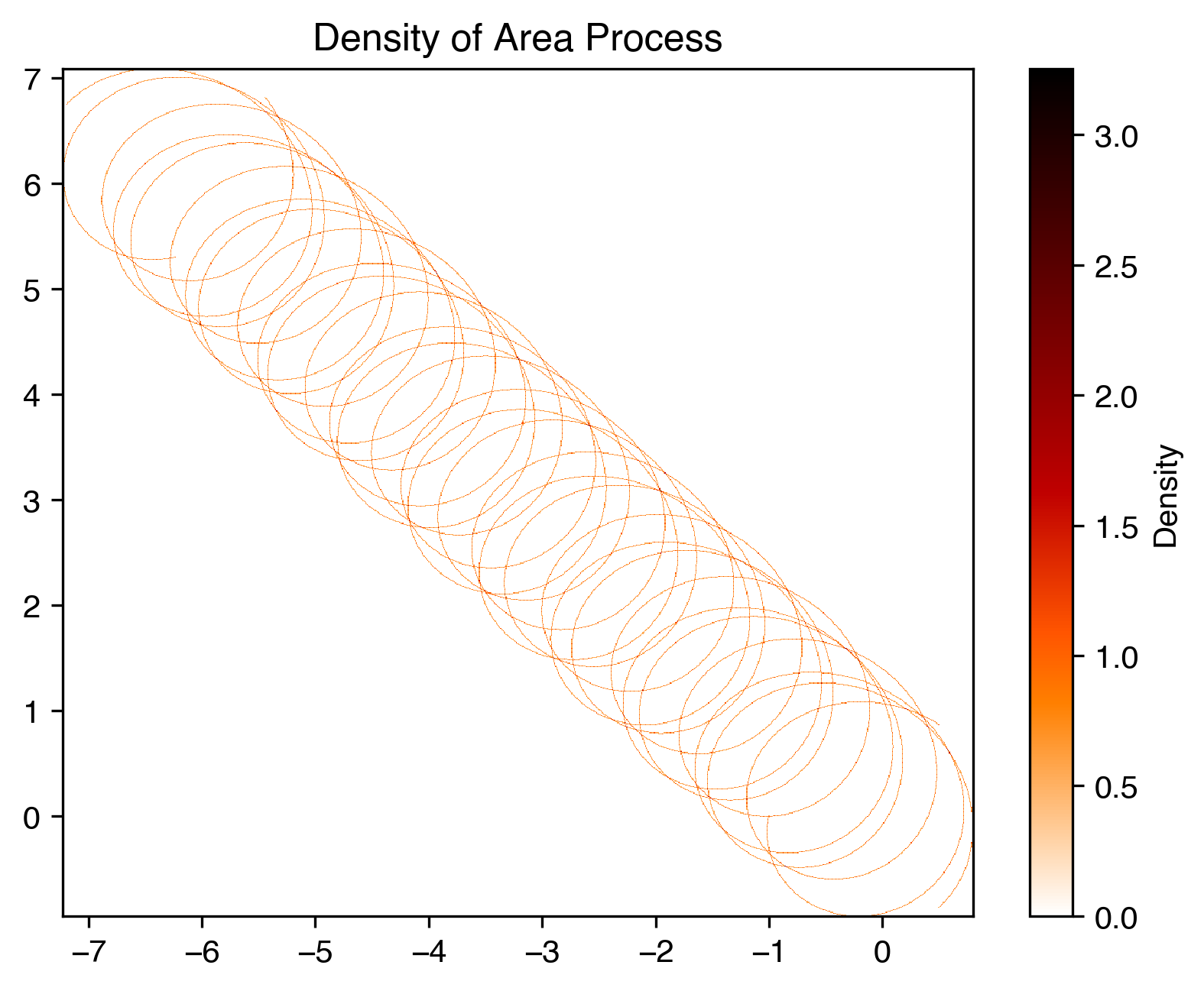}
\caption{ }
\label{fig:PDFArea Process}
\end{subfigure}
\begin{subfigure}[t]{0.325\textwidth}
\centering
\includegraphics[width=.95\textwidth]{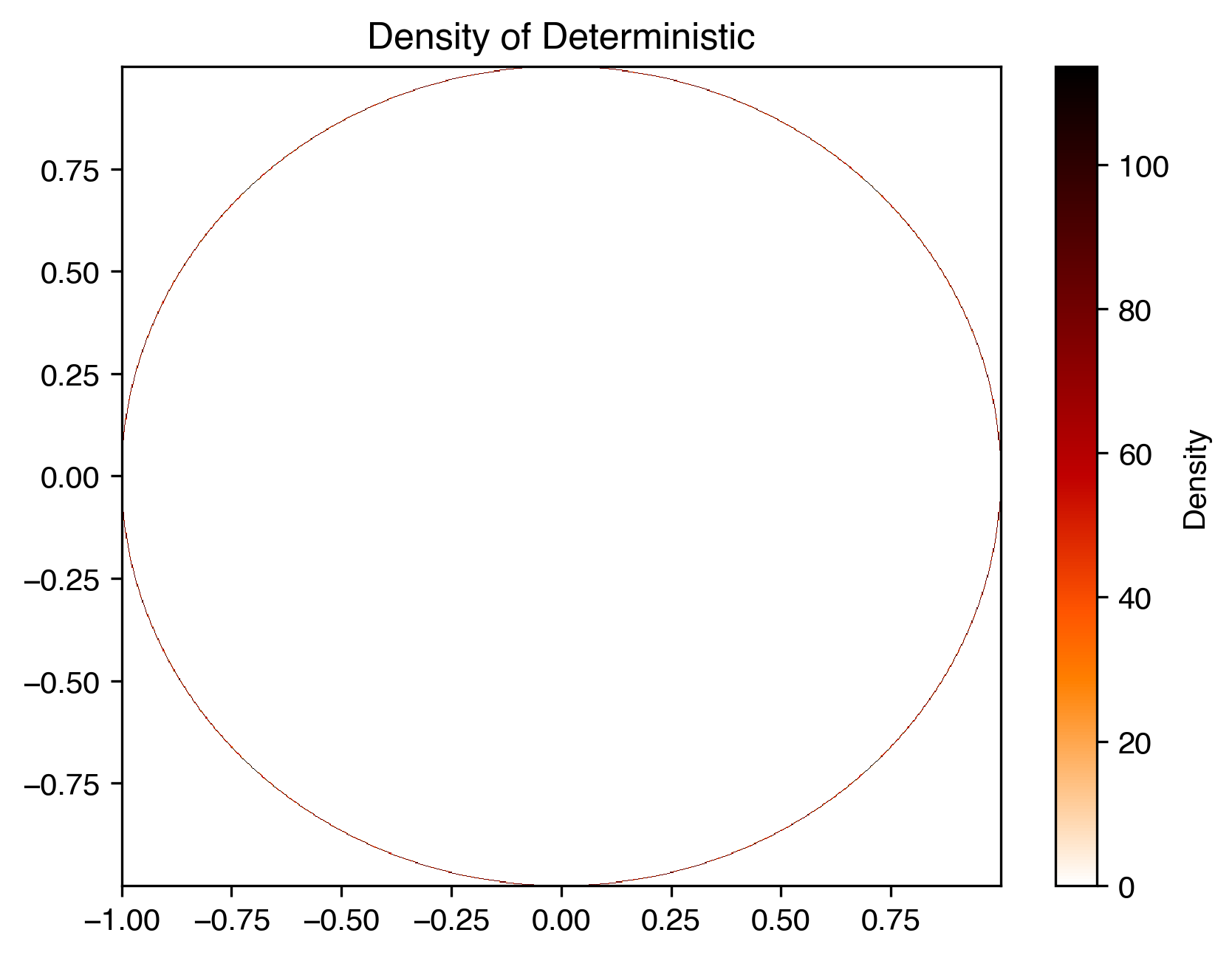}
\caption{}
\label{fig:PDFDeterministic}
\end{subfigure}
\begin{subfigure}[t]{0.325\textwidth}
\centering
\includegraphics[width=.95\textwidth]{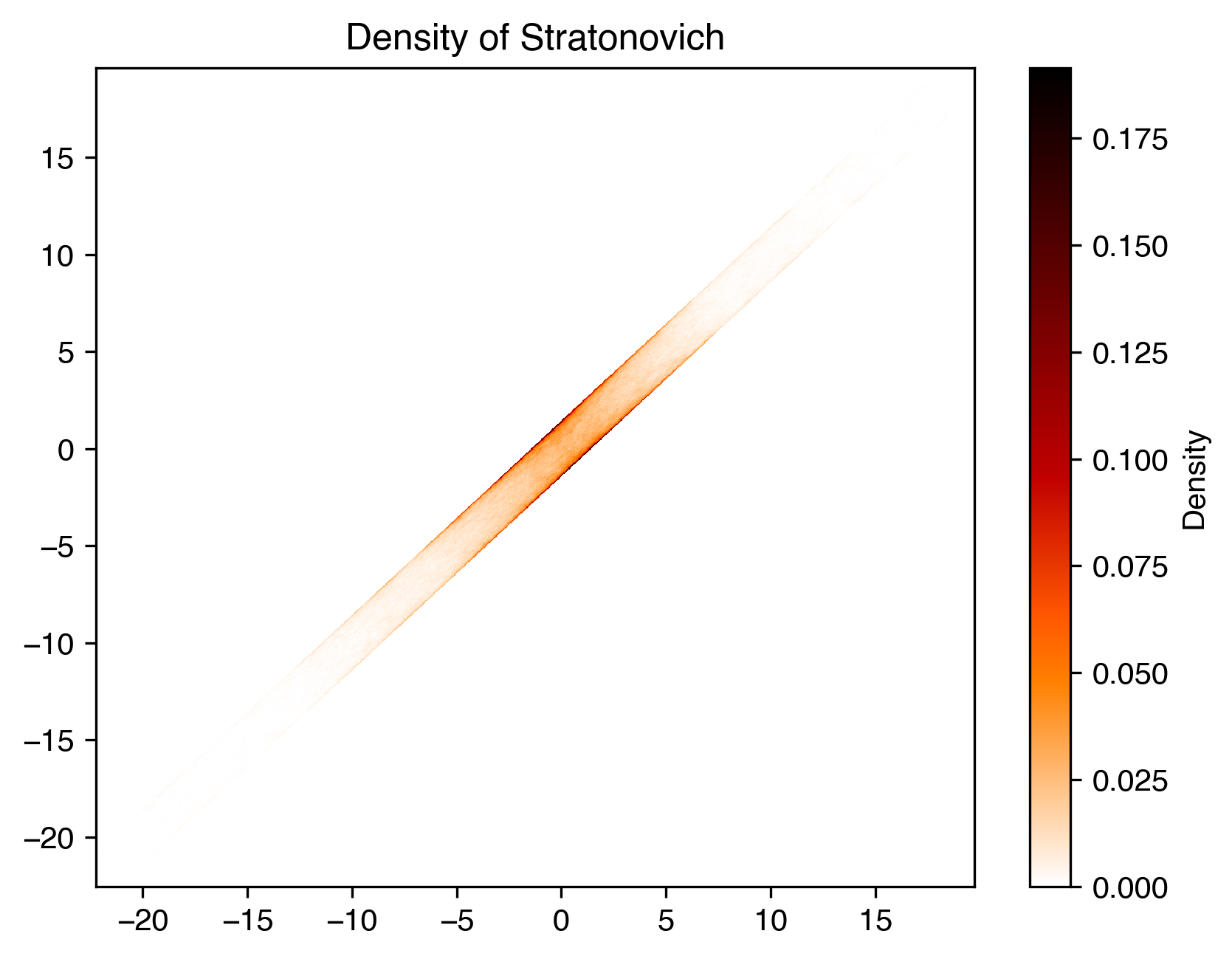}
\caption{}
\label{fig:PDFStratonovich}
\end{subfigure}\\
\begin{subfigure}[t]{0.325\textwidth}
\centering
\includegraphics[width=.95\textwidth]{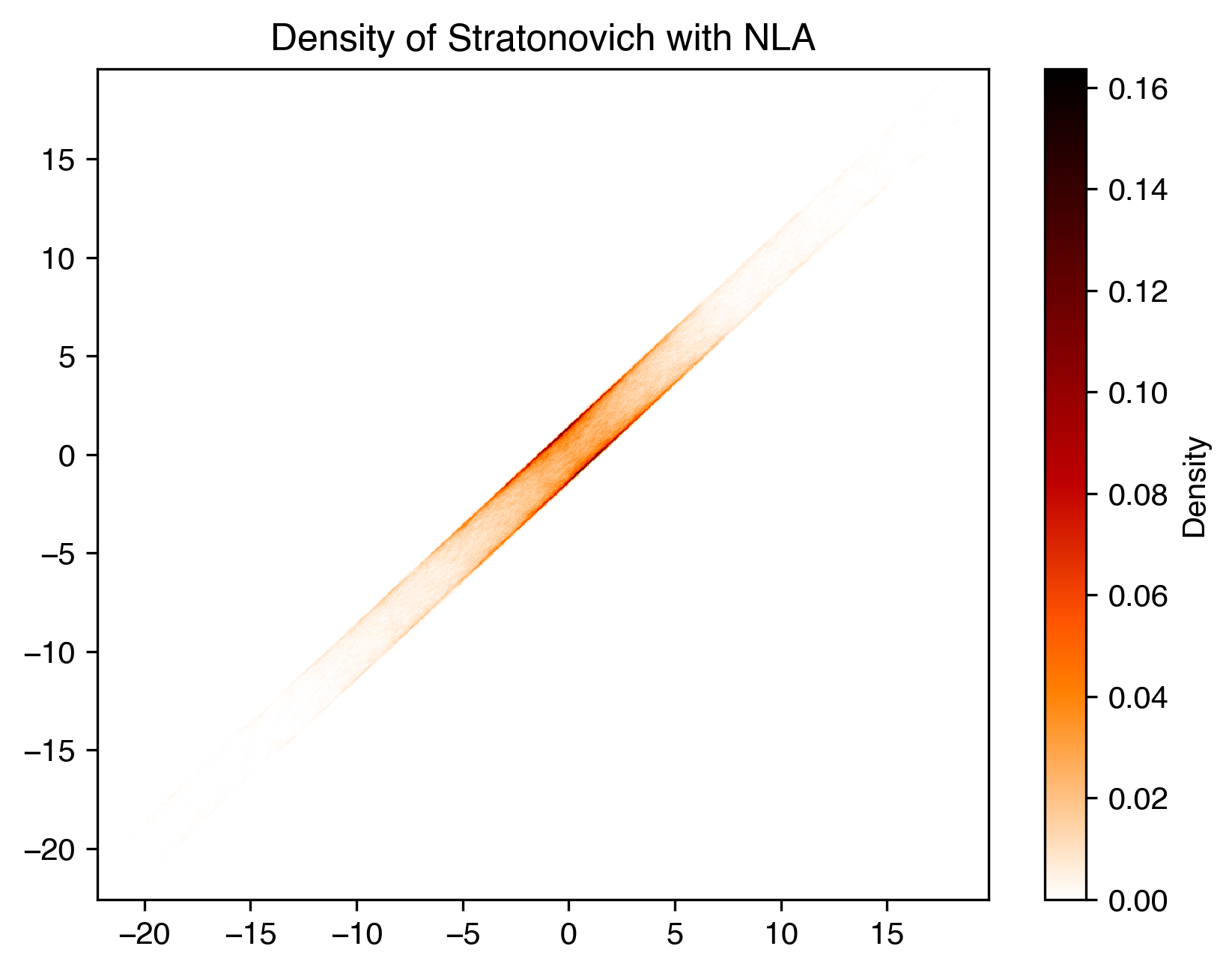}
\caption{ }
\label{fig:PDFStratonovich NLA}
\end{subfigure}
\begin{subfigure}[t]{0.325\textwidth}
\centering
\includegraphics[width=.95\textwidth]{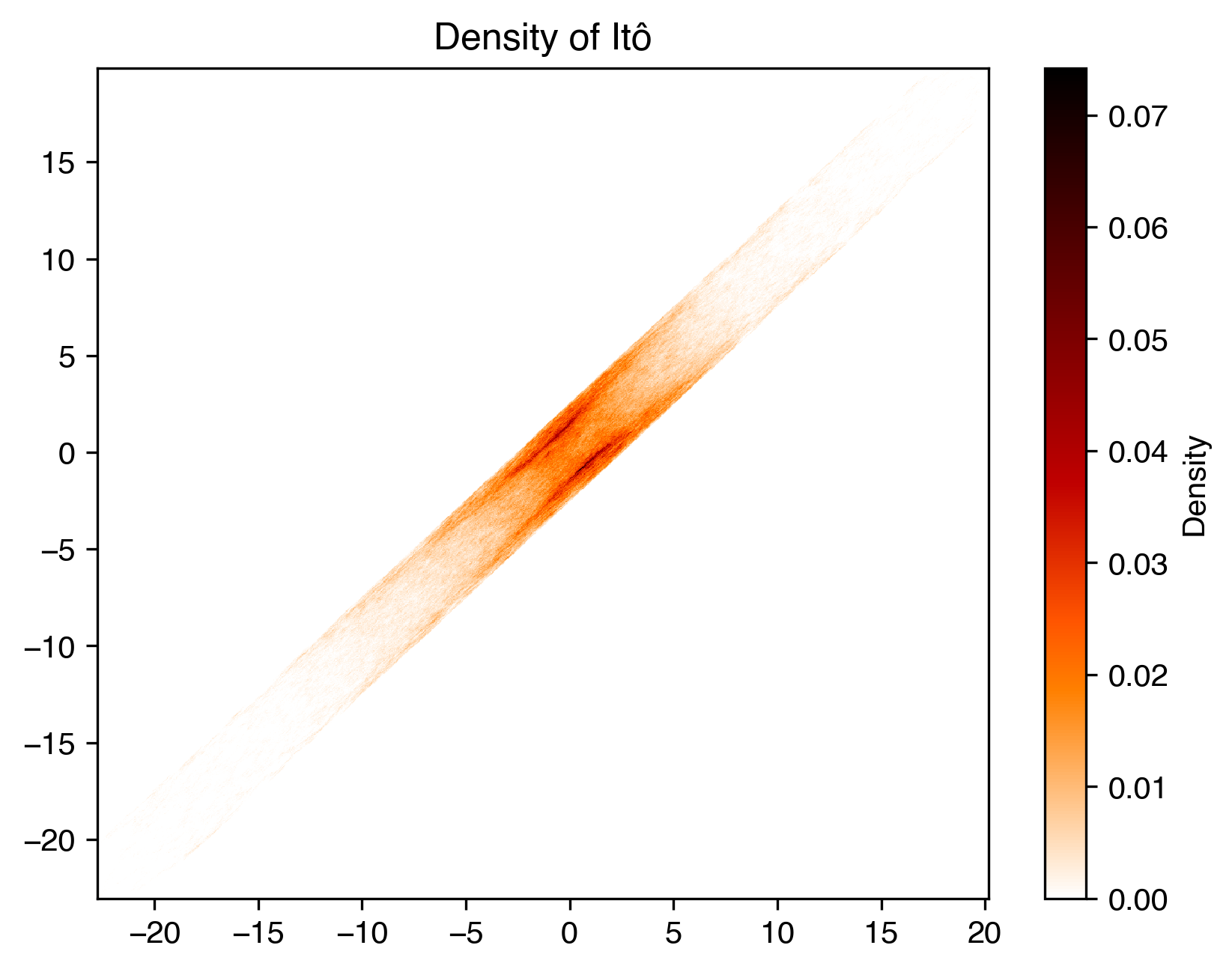}
\caption{ }
\label{fig:PDFItô}
\end{subfigure}
\hspace*{\fill}   
\begin{subfigure}[t]{0.325\textwidth}
\centering
\includegraphics[width=.95\textwidth]{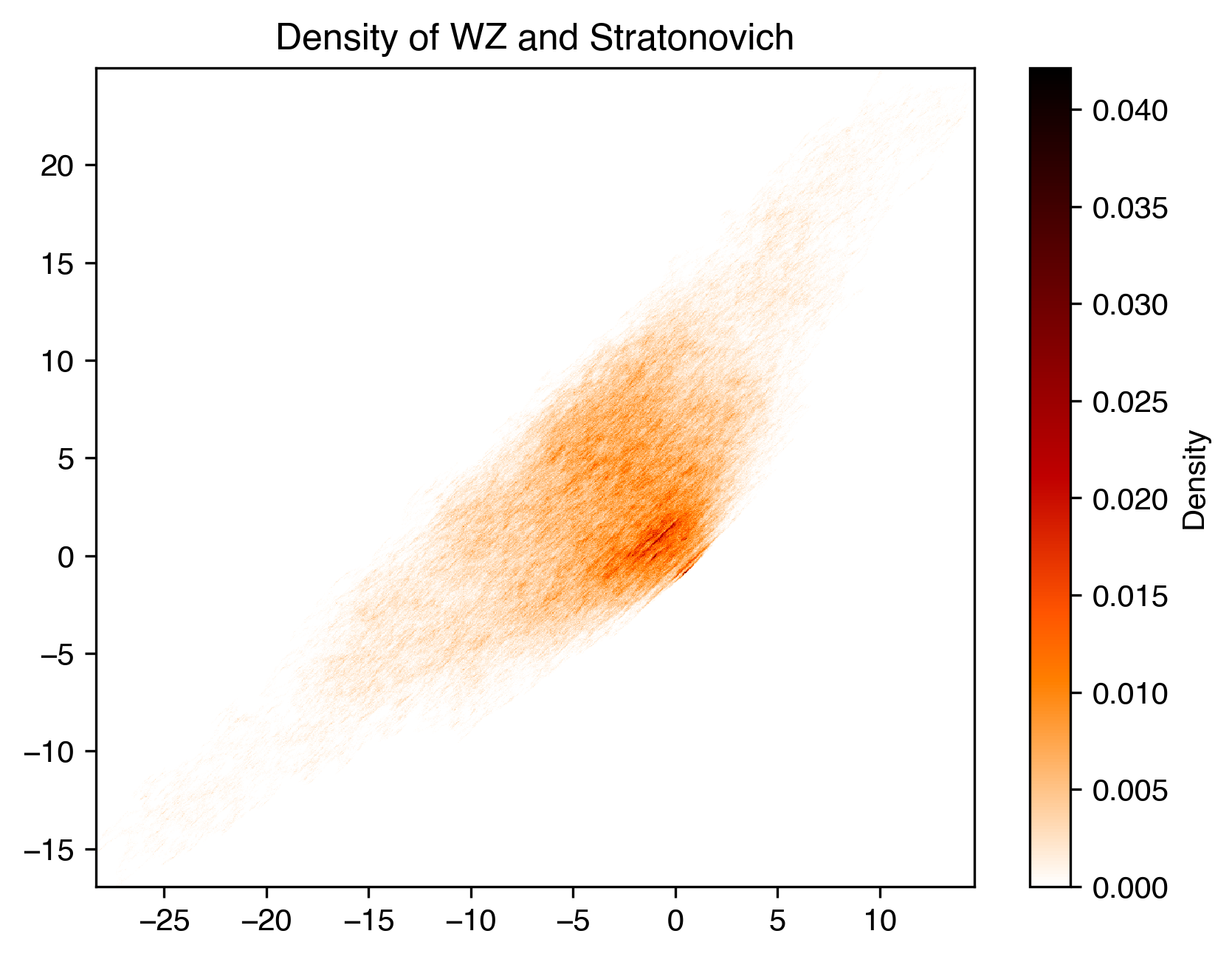}
\caption{ }
\label{fig:PDFWZ and Stratonovich}
\end{subfigure}\\
\begin{subfigure}[t]{0.325\textwidth}
\centering
\includegraphics[width=.95\textwidth]{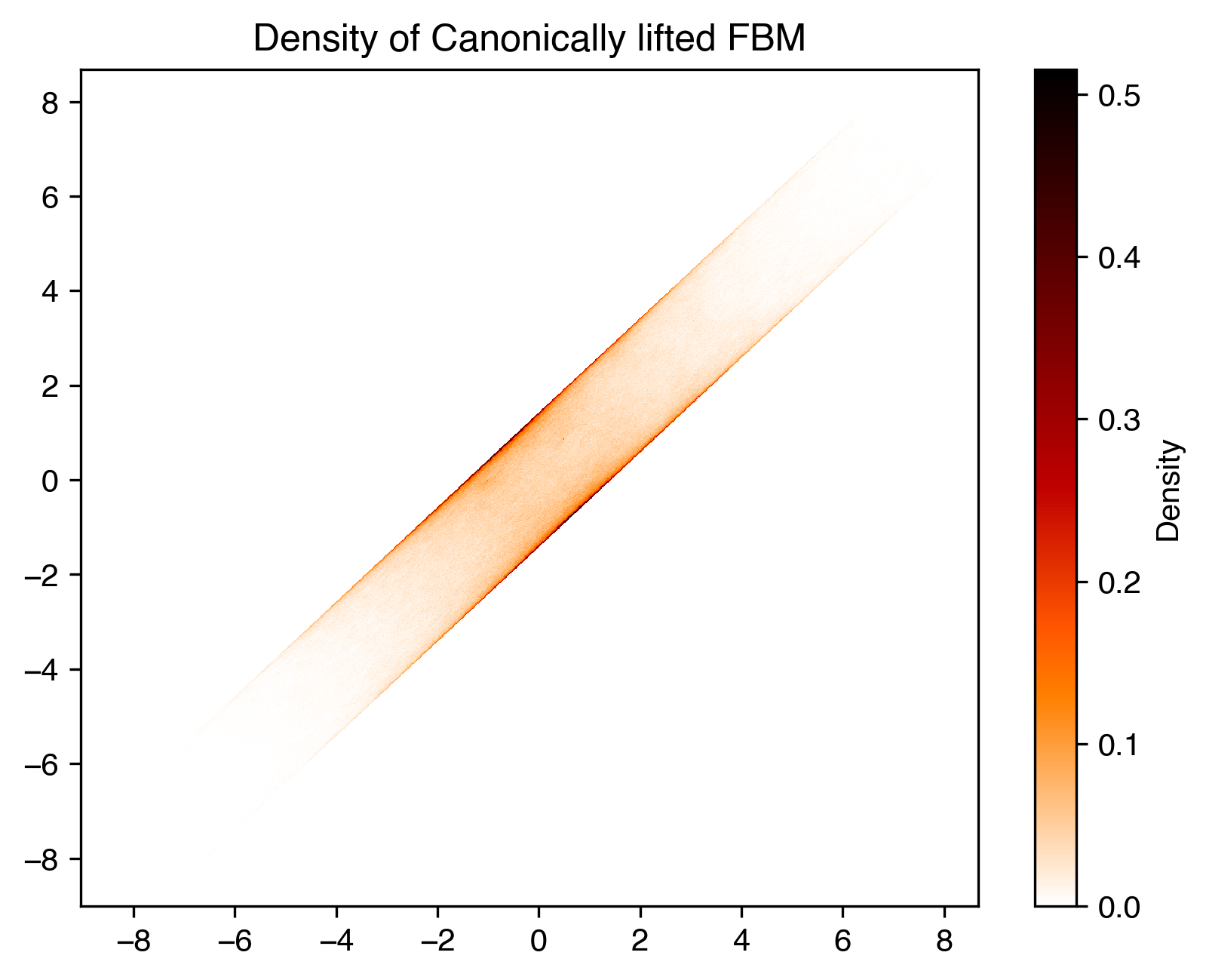}
\caption{ }
\label{fig:PDFfbm}
\end{subfigure}
\caption{ The area weighted histograms of particle trajectories of a 100 member ensemble over the time interval $[0,40]$. Values are binned into an array of size $1024 \times 1024$ and the number of occurrences are used to compute the area weighted density. }
\label{figs:PDF's}
\end{figure}

\begin{figure}[H]
\centering
\begin{subfigure}[t]{0.495\textwidth}
\centering
\includegraphics[width=.95\textwidth]{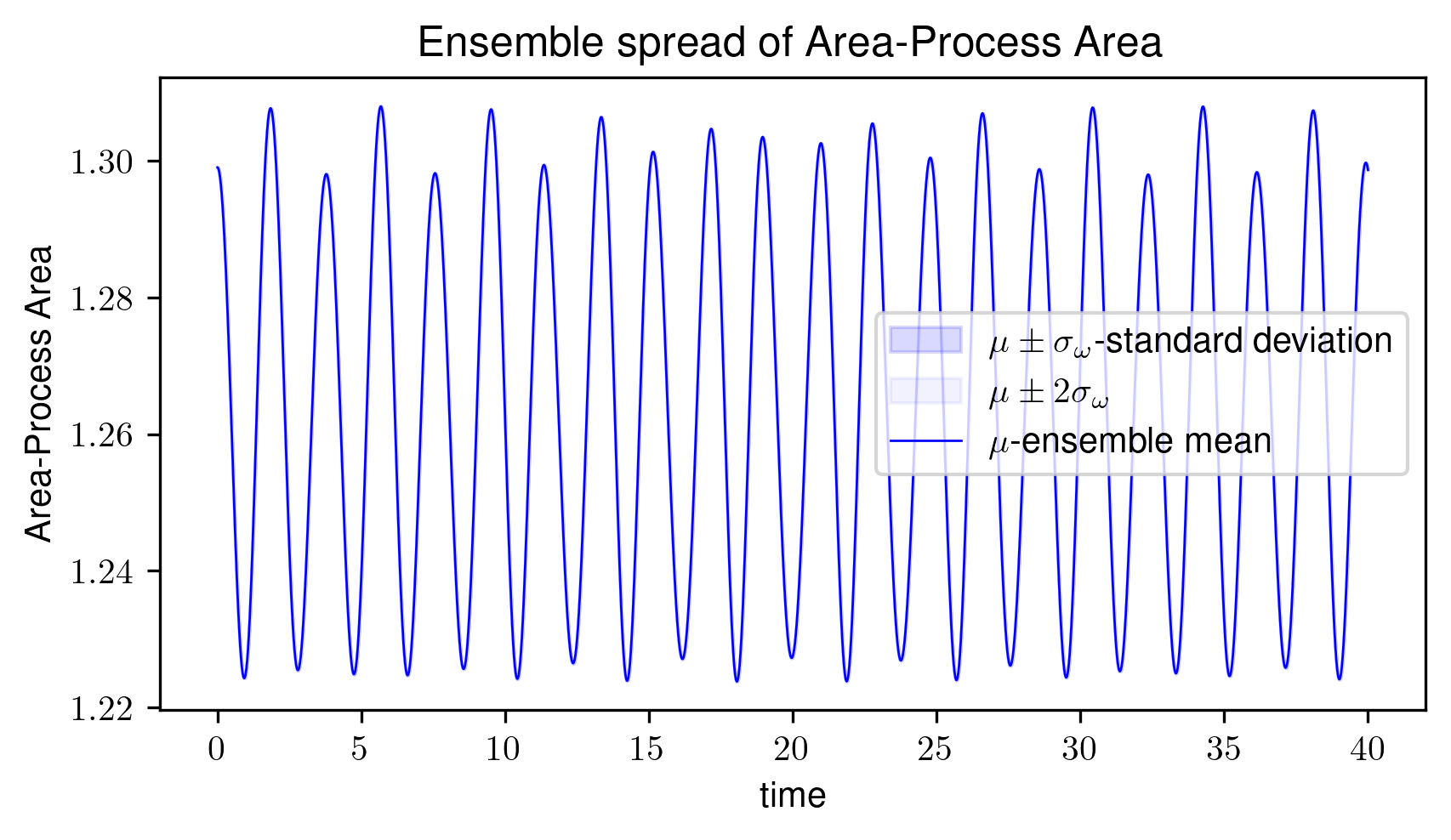}
\caption{ }
\label{fig:_spaghetti_Area-Process Area}
\end{subfigure}
\begin{subfigure}[t]{0.495\textwidth}
\centering
\includegraphics[width=.95\textwidth]{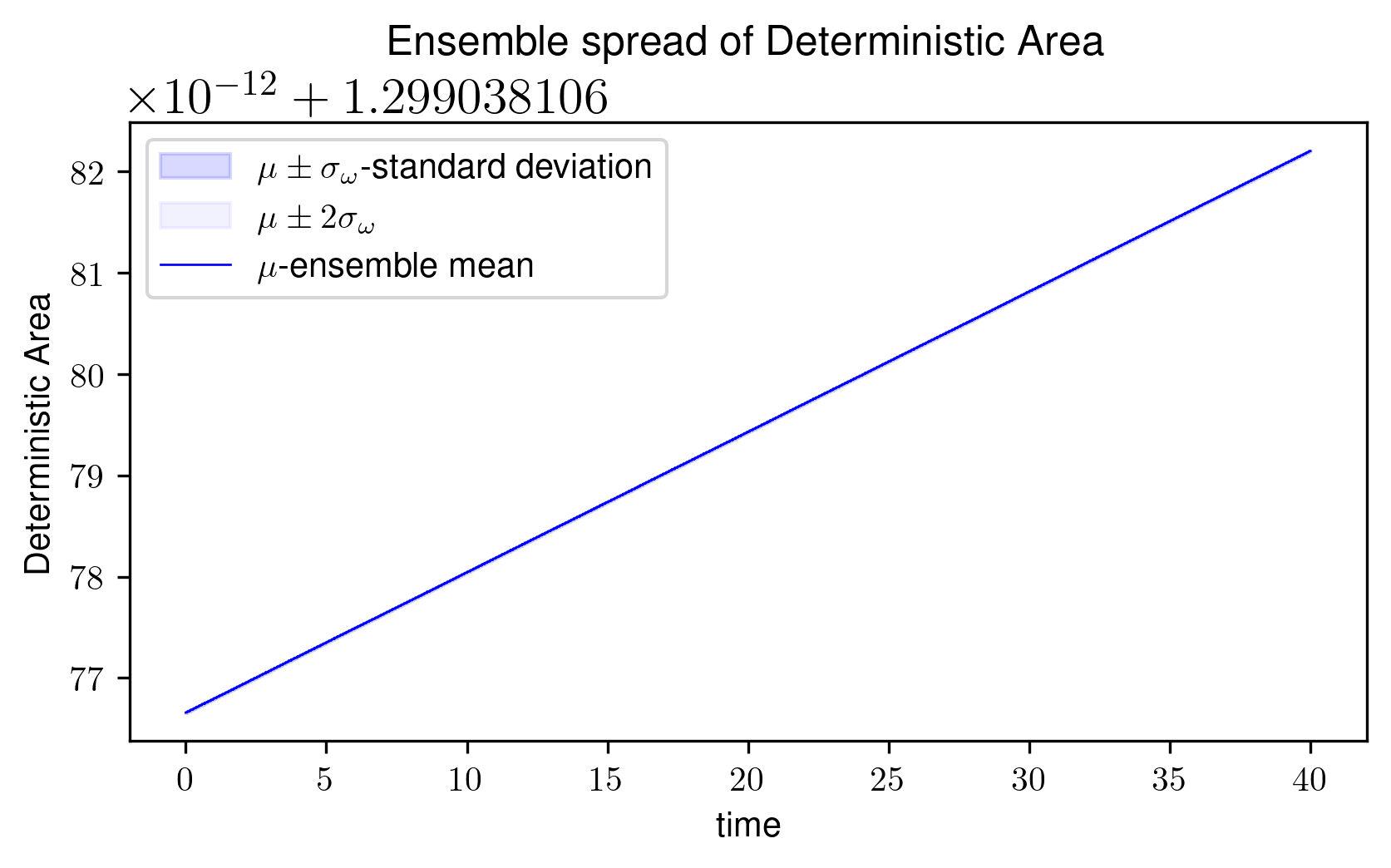}
\caption{}
\label{fig:_spaghetti_Deterministic Area}
\end{subfigure}\\
\begin{subfigure}[t]{0.495\textwidth}
\centering
\includegraphics[width=.95\textwidth]{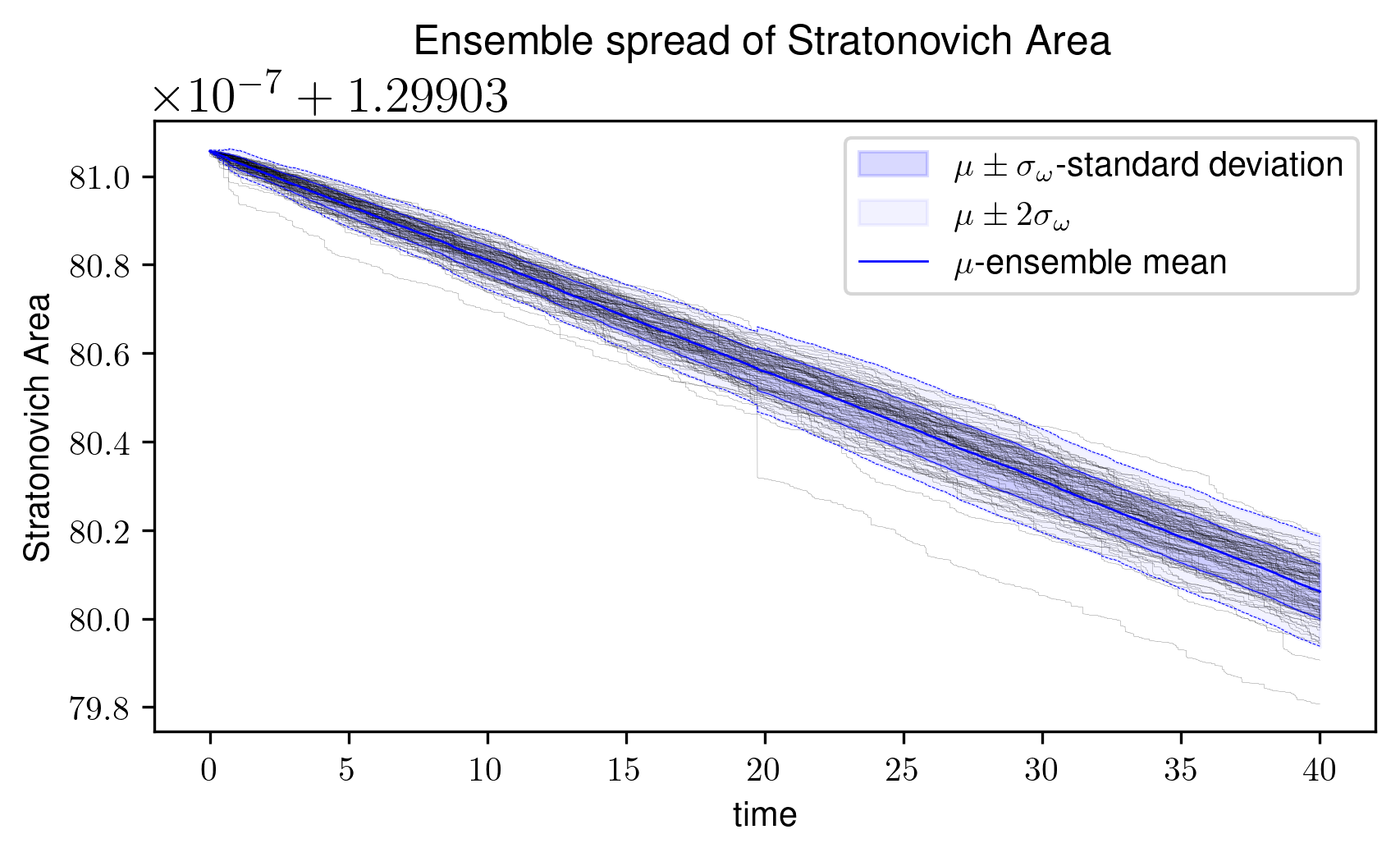}
\caption{}
\label{fig:_spaghetti_Stratonovich Area}
\end{subfigure}
\begin{subfigure}[t]{0.495\textwidth}
\centering
\includegraphics[width=.95\textwidth]{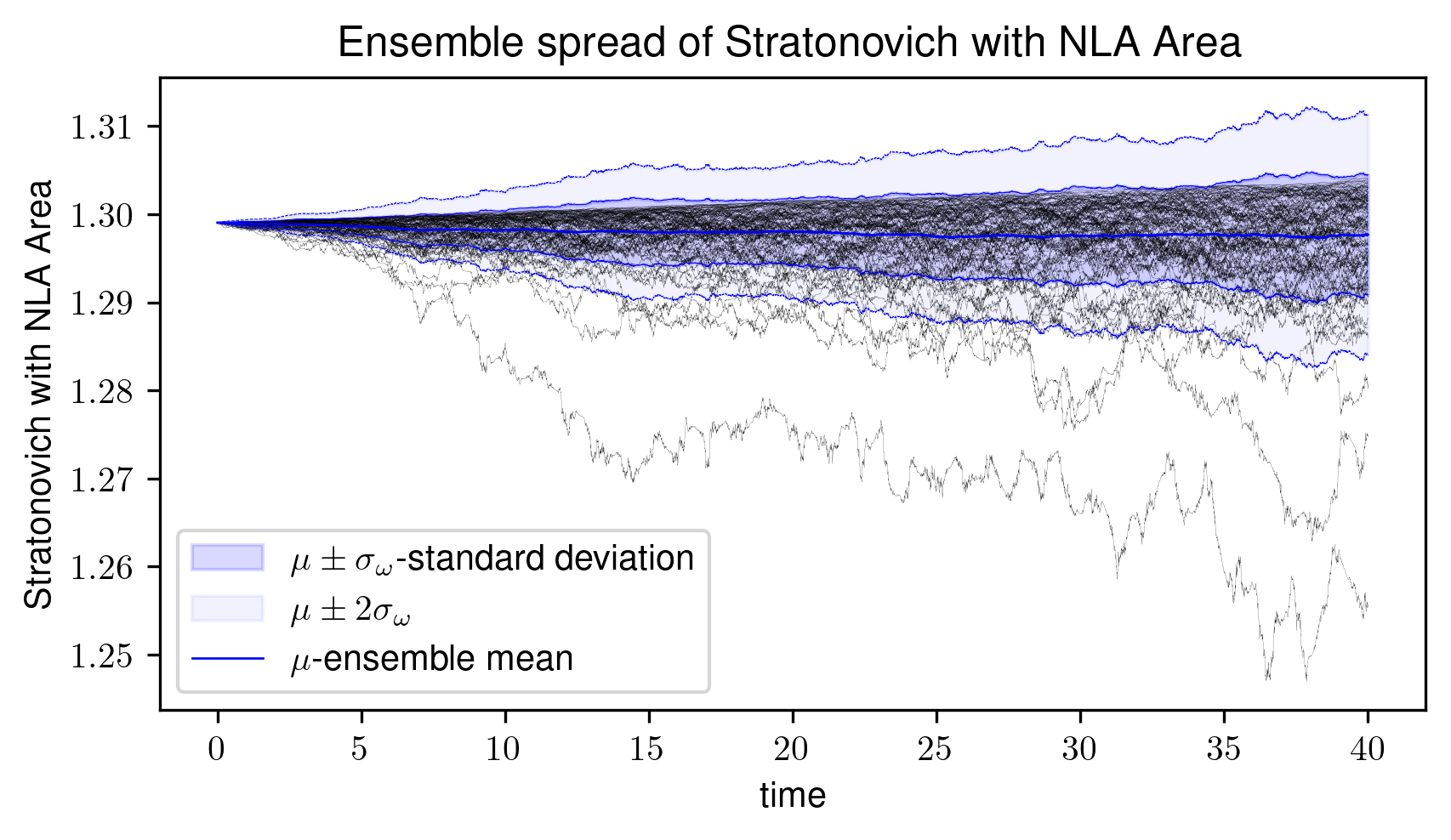}
\caption{ }
\label{fig:_spaghetti_Stratonovich with NLA Area}
\end{subfigure}\\
\begin{subfigure}[t]{0.495\textwidth}
\centering
\includegraphics[width=.95\textwidth]{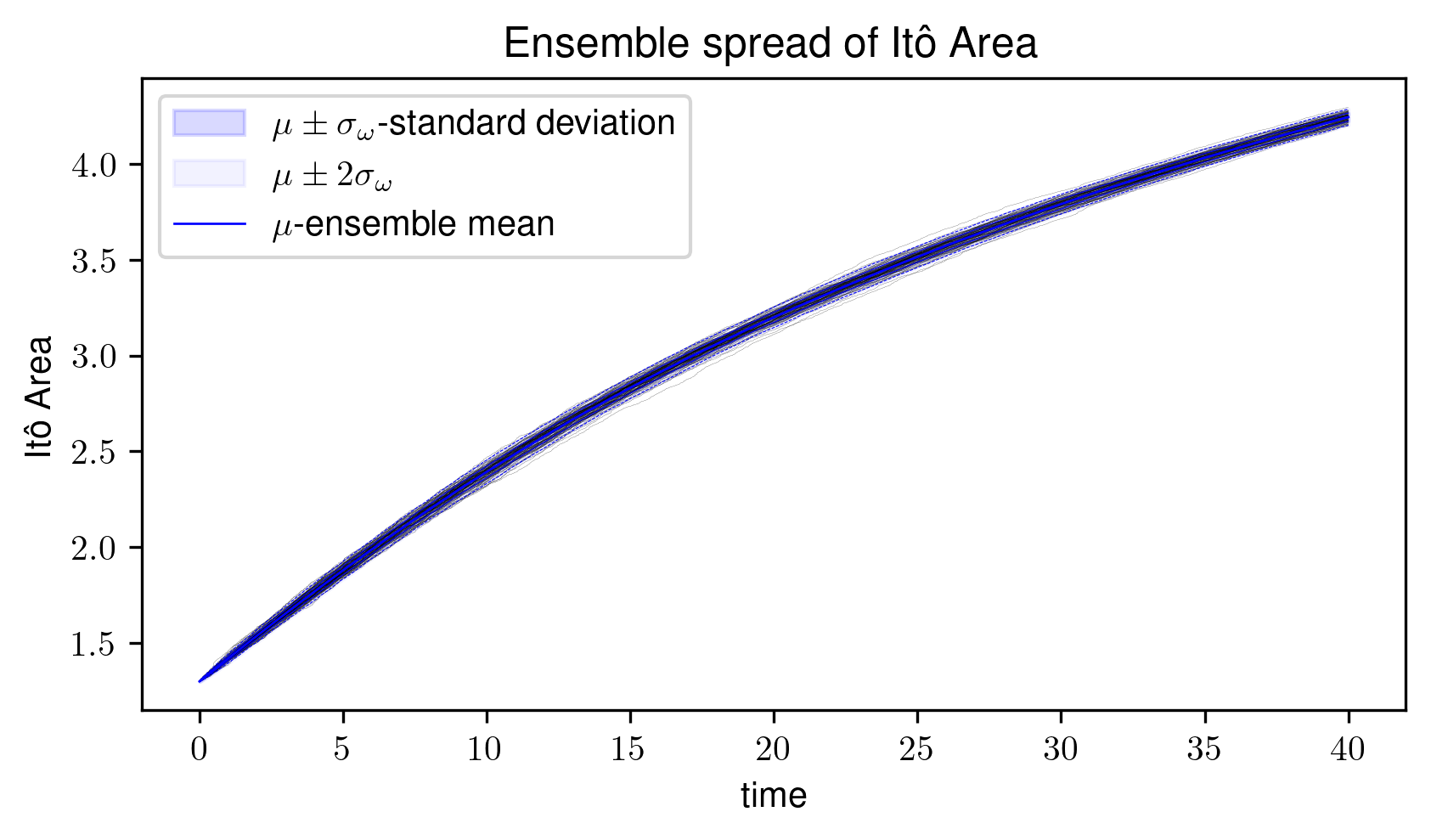}
\caption{ }
\label{fig:_spaghetti_Itô Area}
\end{subfigure}
\begin{subfigure}[t]{0.495\textwidth}
\centering
\includegraphics[width=.95\textwidth]{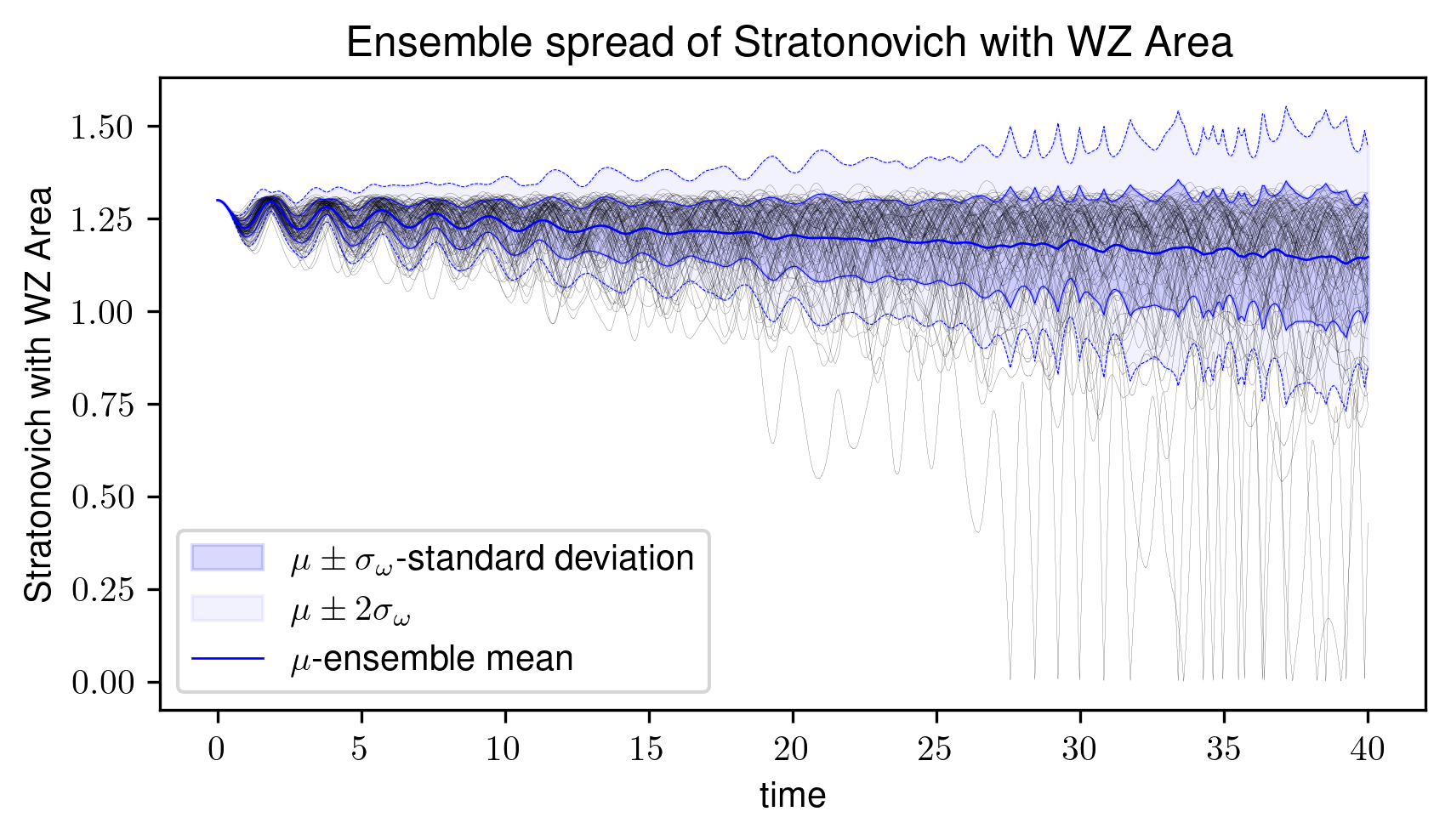}
\caption{ }
\label{fig:_spaghetti_Stratonovich with WZ Area}
\end{subfigure}\\
\begin{subfigure}[t]{0.495\textwidth}
\centering
\includegraphics[width=.95\textwidth]{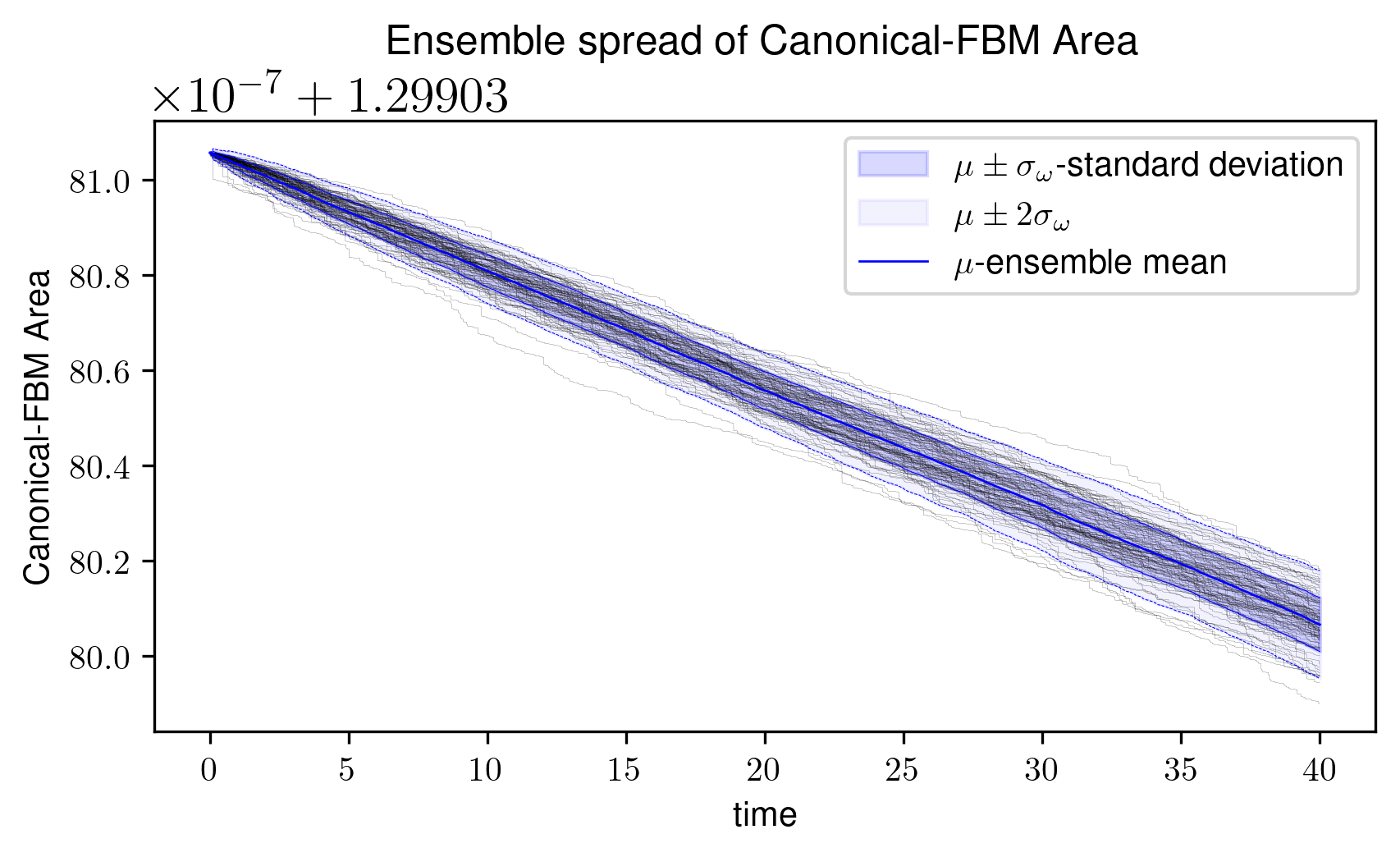}
\caption{ }
\label{fig:_spaghetti_Canonical-FBM Area}
\end{subfigure}
\caption{ Area between the three point vortices as time evolves, for a 100 member ensemble, for \cref{Method:type II},  \cref{Method:Deterministic}, \cref{Method:Stratonovich integration scheme}, \cref{Method:NLA}, \cref{Method:type ito},  \cref{Method:type I}, and \cref{Method:FBM} respectively. We plot the traced paths of ensemble member area in black, also plotted is the mean, the standard deviation, and twice the standard deviation. Data is not necessarily normal and the variance doesn't necessarily provide accurate representation of the ensemble statistics.  }
\label{figs:Area}
\end{figure}


\begin{figure}[H]
\centering
\begin{subfigure}[t]{0.495\textwidth}
\centering
\includegraphics[width=.95\textwidth]{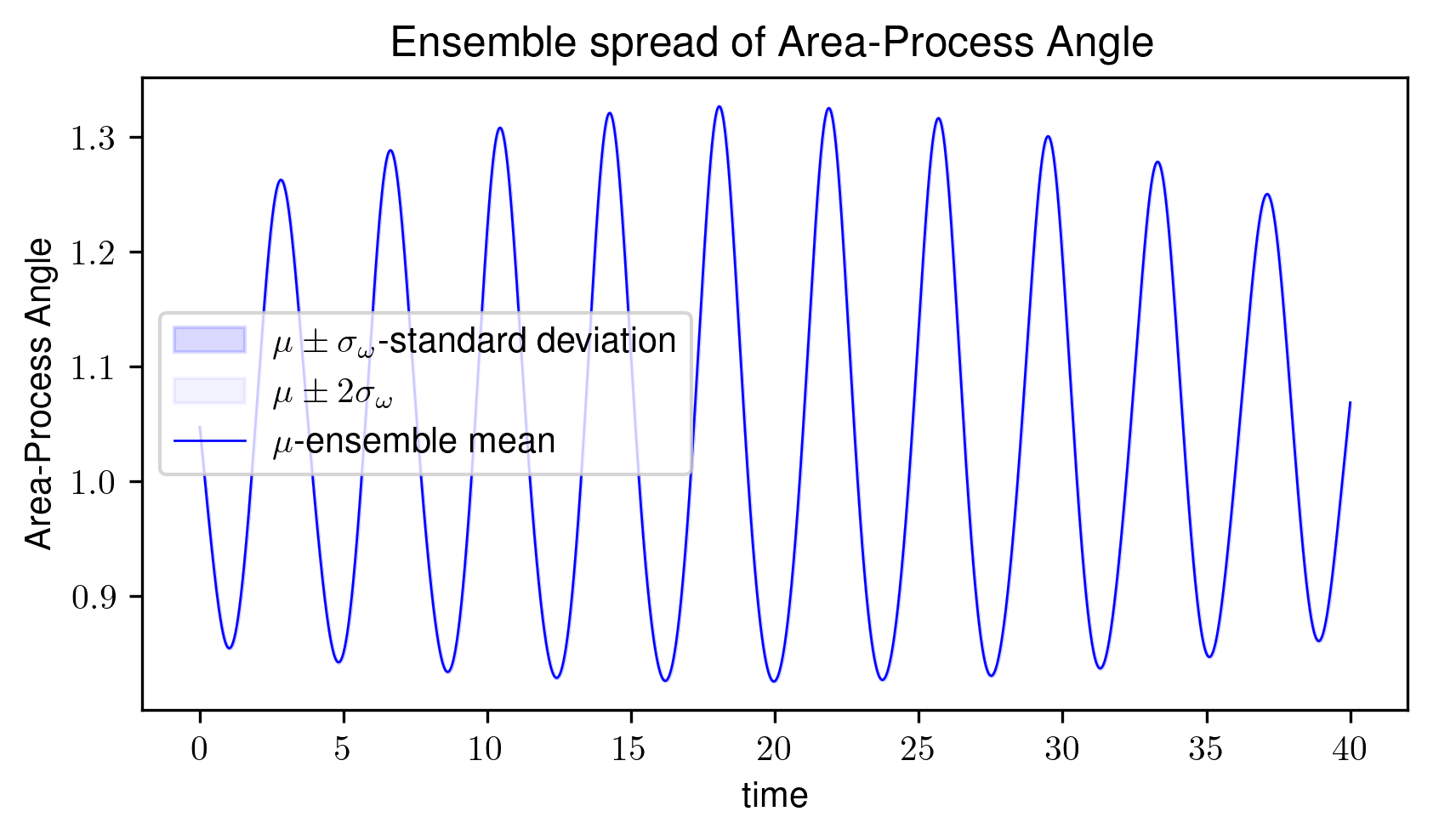}
\caption{ }
\label{fig:_spaghetti_Area-Process Angle}
\end{subfigure}
\begin{subfigure}[t]{0.495\textwidth}
\centering
\includegraphics[width=.95\textwidth]{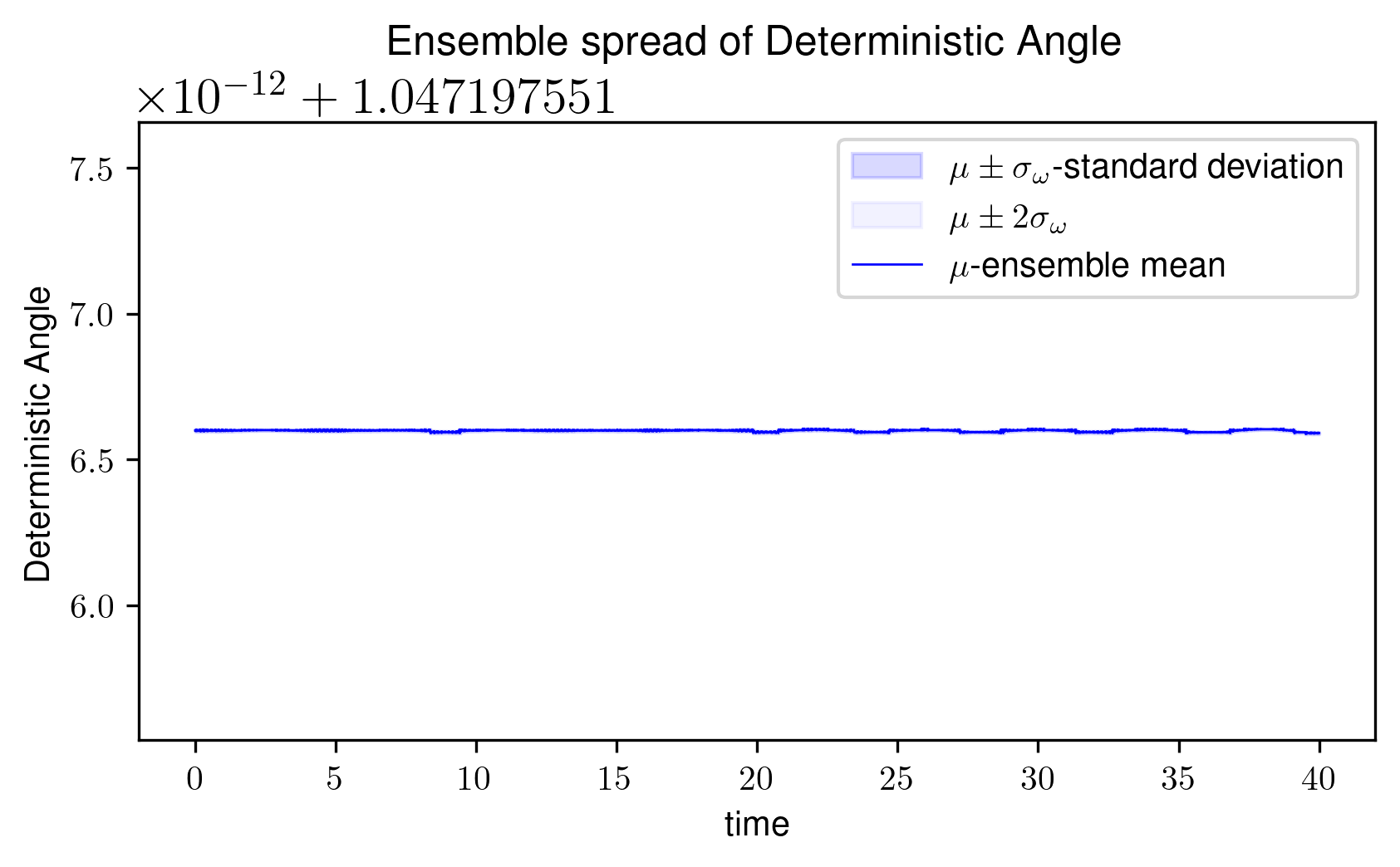}
\caption{}
\label{fig:_spaghetti_Deterministic Angle}
\end{subfigure}\\
\begin{subfigure}[t]{0.495\textwidth}
\centering
\includegraphics[width=.95\textwidth]{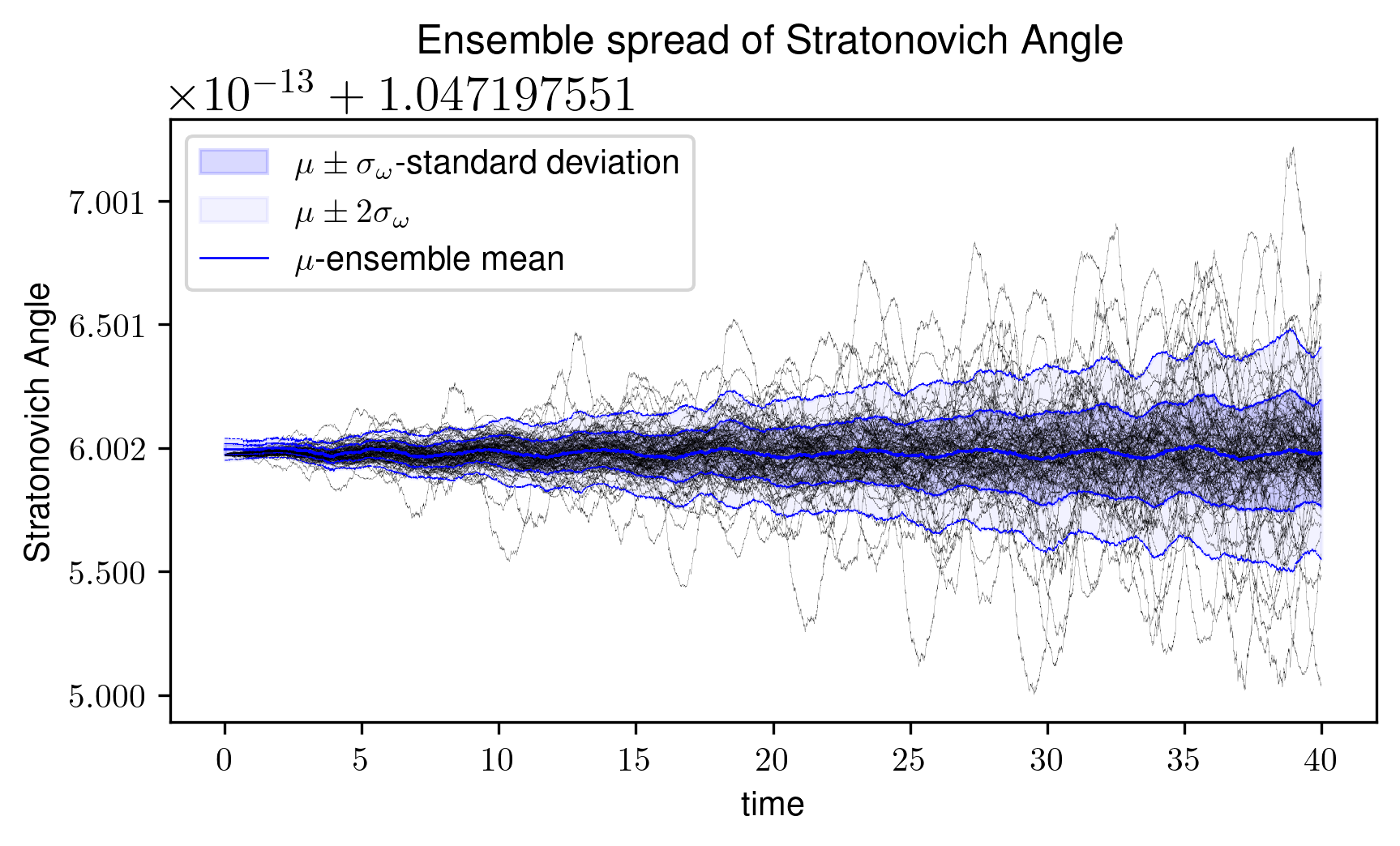}
\caption{}
\label{fig:_spaghetti_Stratonovich Angle}
\end{subfigure}
\begin{subfigure}[t]{0.495\textwidth}
\centering
\includegraphics[width=.95\textwidth]{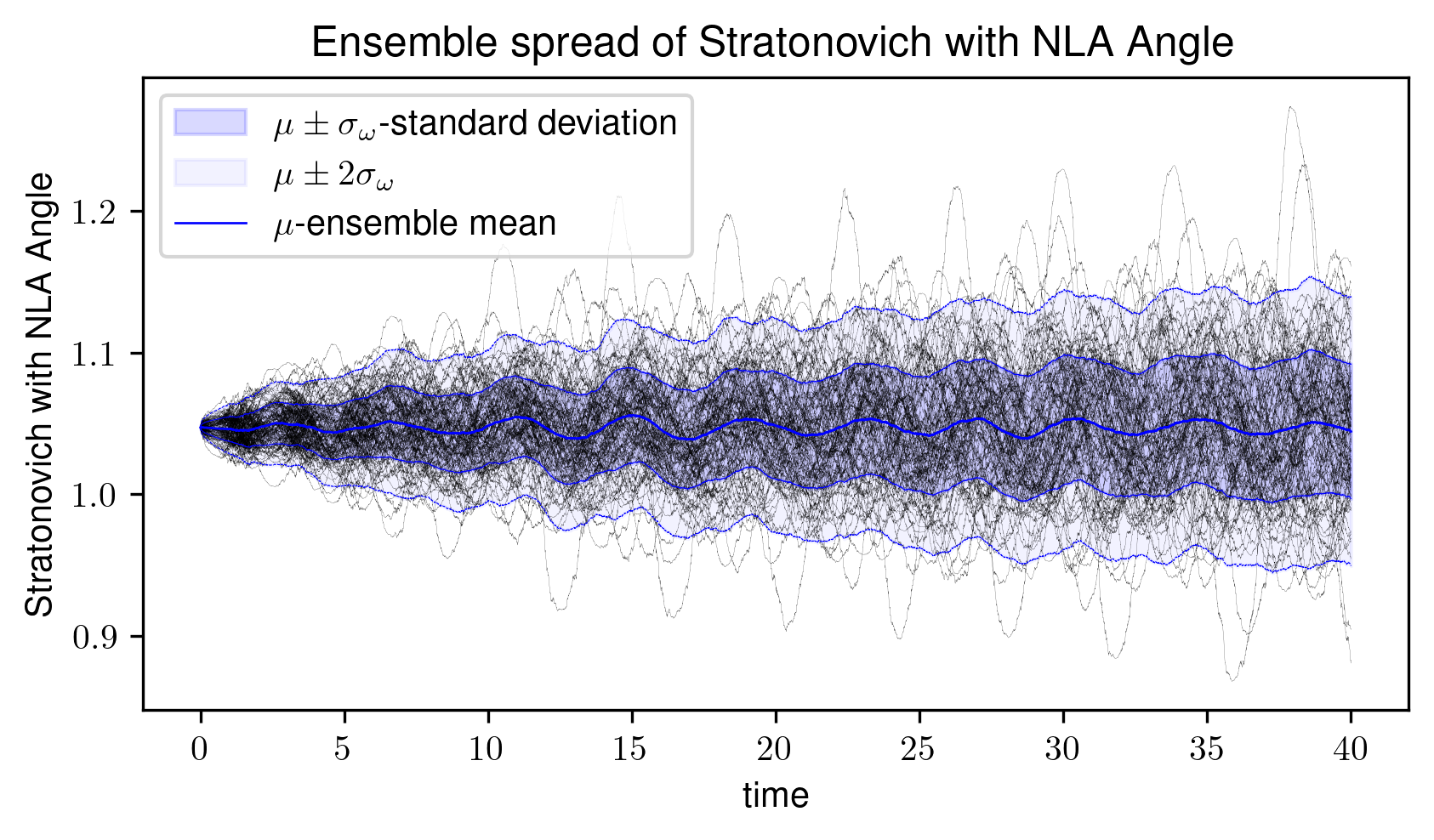}
\caption{ }
\label{fig:_spaghetti_Stratonovich with NLA Angle}
\end{subfigure}\\
\begin{subfigure}[t]{0.495\textwidth}
\centering
\includegraphics[width=.95\textwidth]{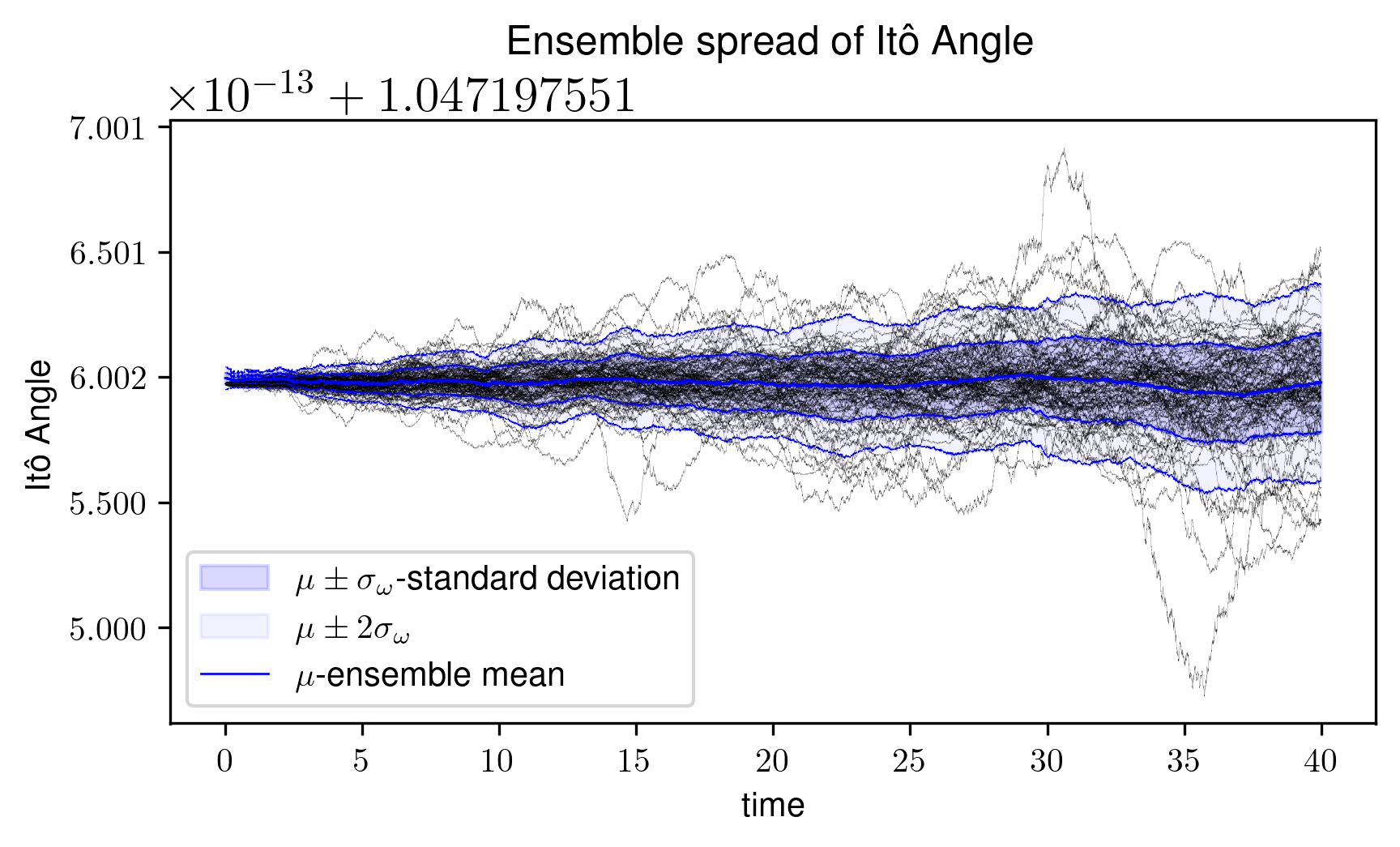}
\caption{ }
\label{fig:_spaghetti_Itô Angle}
\end{subfigure}
\begin{subfigure}[t]{0.495\textwidth}
\centering
\includegraphics[width=.95\textwidth]{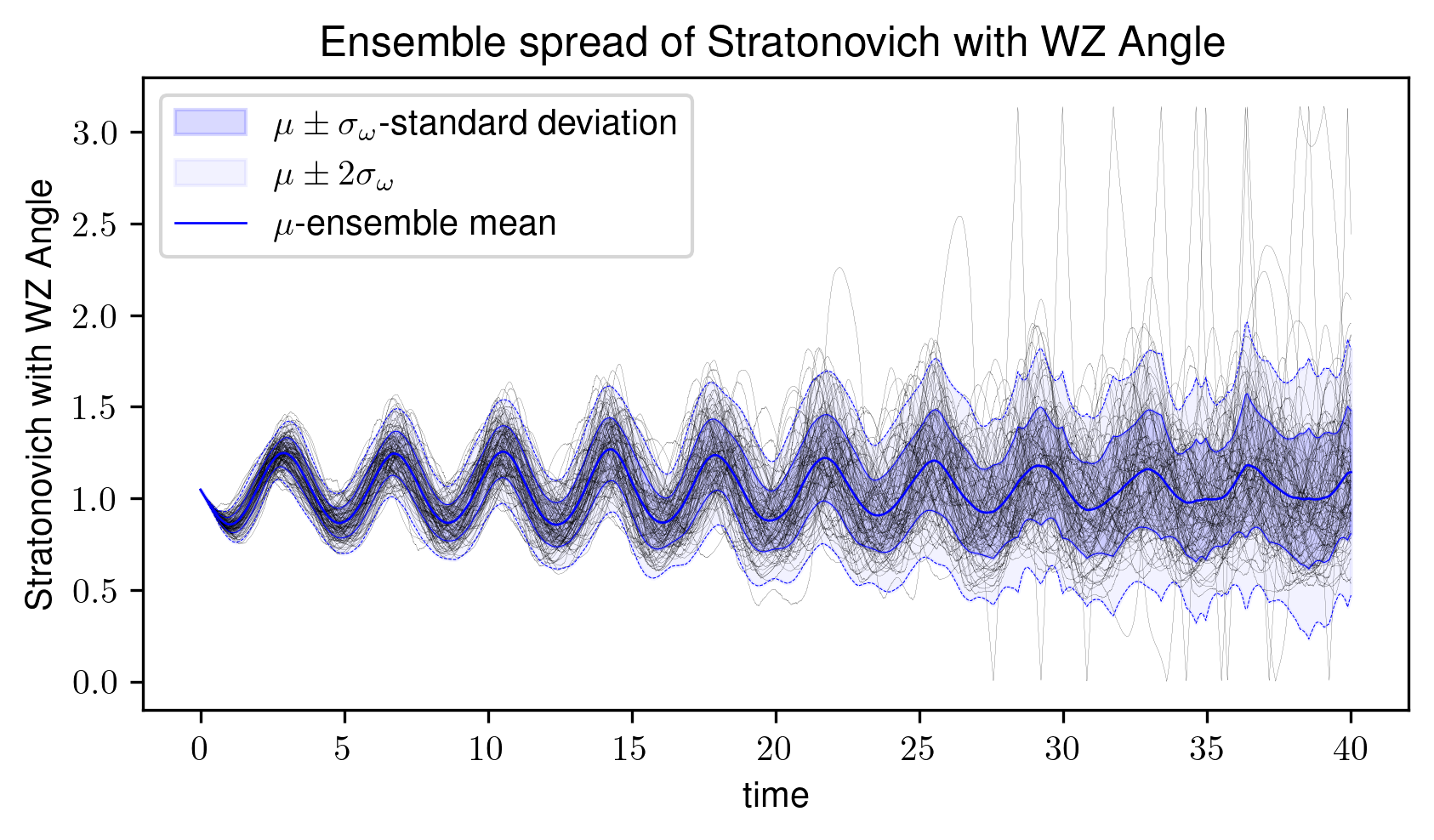}
\caption{ }
\label{fig:_spaghetti_Stratonovich with WZ Angle}
\end{subfigure}\\
\begin{subfigure}[t]{0.495\textwidth}
\centering
\includegraphics[width=.95\textwidth]{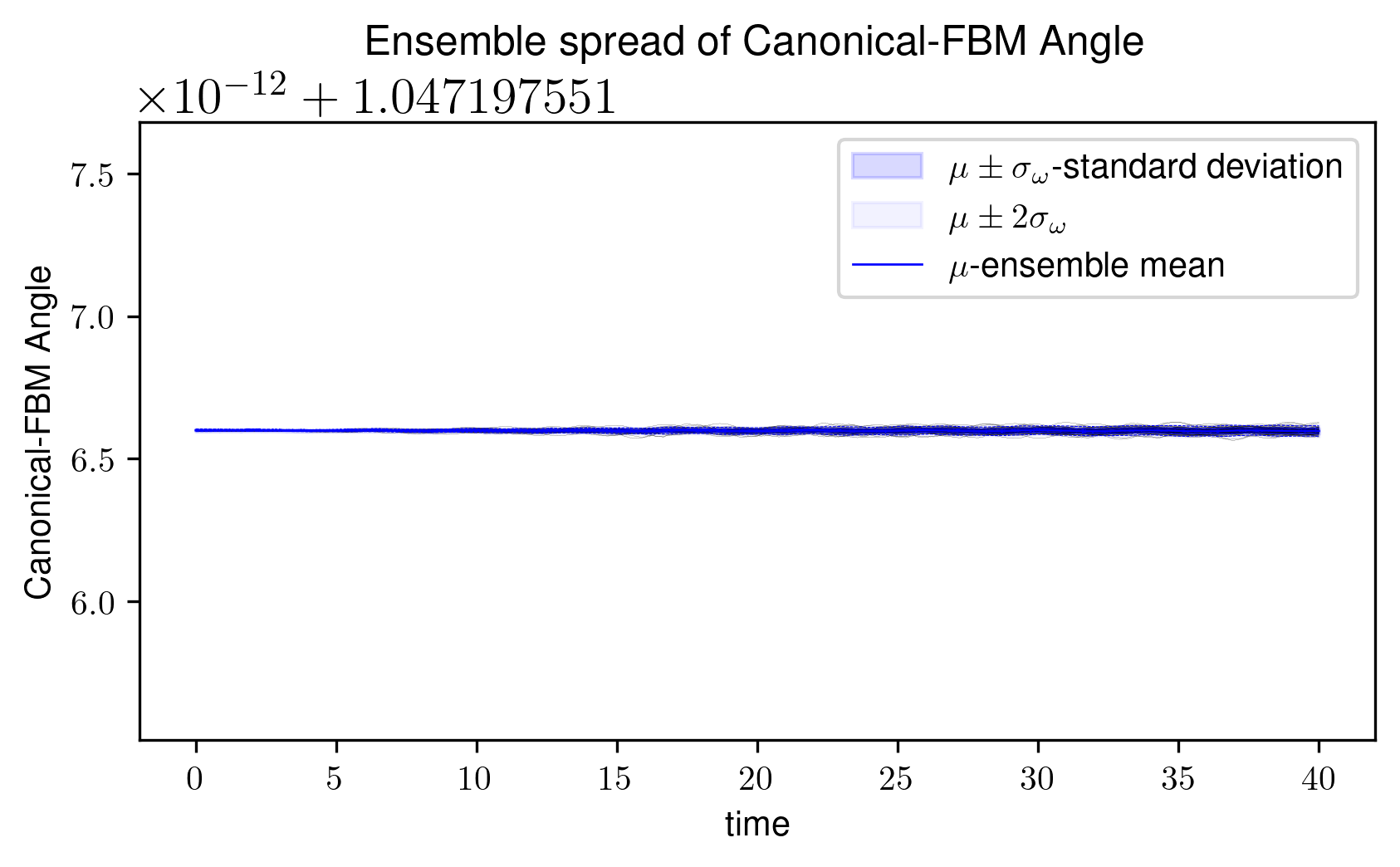}
\caption{ }
\label{fig:_spaghetti_Canonical-FBM Angle}
\end{subfigure}
\caption{  Angle between the three point vortices as time evolves, for a 100 member ensemble, for \cref{Method:type II},  \cref{Method:Deterministic}, \cref{Method:Stratonovich integration scheme}, \cref{Method:NLA}, \cref{Method:type ito},  \cref{Method:type I}, and \cref{Method:FBM} respectively. We plot the traced paths of ensemble member angle in black, also plotted is the mean angle, the standard deviation, and twice the standard deviation of angle. Data is not necessarily normal and the variance doesn't necessarily provide accurate representation of the ensemble statistics.  }
\label{figs:Angle}
\end{figure}

\subsection{Results}\label{sec:results}


For the experiment described in \cref{sec: case study} the stream functions $\psi^1,\psi^2$ have parameters, $r=1$, $A = 1/2$, $a_1=1$, $a_2=-1$. The time interval $t\in [0,40]$ is split into $10000$ timesteps, with step-size $\Delta t = (250)^{-1}$. The $10000$ increments of Brownian motion $\Delta W^1, \Delta W^2$ are generated as realisations of a normal distribution and scaled by $\sqrt{\Delta t}$.  The canonical FBM-driven point vortex system simulation also uses $10000$ increments generated using the Davies-Harte algorithm, with Hurst parameter $0.4$. The adaptive regularisation of the Greens function uses parameter $\delta_0 = \sqrt{10^{-10}}$. In \cref{Method:NLA} we take $K=10000>(500)$ terms in the truncated Fourier series of the Brownian bridge process, consistent with the bound required for strong order one convergence \cite{kloeden1992approximation}.

\Cref{figs:PDF's} contains the area-weighted histograms of all particle trajectories of a 100-member ensemble over the time interval $[0,40]$. Values are binned into an array of size $1024 \times 1024$ and the number of occurrences are counted and used to compute an area-weighted density, serving as a proxy for probability of occurrence in a region of phase space over the given time interval. \Cref{figs:Area} contains the computed area between the three points, as time evolves for all 100 members of an ensemble, in particular the mean, standard deviation, and twice the standard deviation, as well as the trajectories are plotted. \Cref{figs:Angle} contains the computed angle between three points, as time evolves for all 100 members of an ensemble, plotted is the mean, standard deviation, twice the standard deviation, and trajectories. The 
angle is computed as $\alpha_{123} = \cos^{-1}(\frac{\b x_1 \cdot \b x_2}{\|\b x_1 \|_2 \| \b x_2 \|_2 })$, the area is computed by Herons formula, $A_{123} = (s (s-\|\b x_1-\b x_2\|_2) (s-\|\b x_2-\b x_3\|_2) (s-\|\b x_3-\b x_1\|_2))^{\frac{1}{2}}$ where $s = \frac{1}{2} (\|\b x_1-\b x_2\|_2 + \|\b x_2-\b x_3\|_2 + \|\b x_3-\b x_1\|_2)$ is the semiperimiter.

The Deterministic RK4 scheme \cref{Method:Deterministic} observes the expected behavior, the point vortices trace a circle (as seen in \cref{fig:PDFDeterministic}) and preserve the area $3\sqrt{3}/4 \approx 1.047$ and angle $\pi/3\approx 1.299$ in the inscribed equilateral triangle to $10^{-12}$ over the interval $[0,40]$ as can be seen in \cref{fig:_spaghetti_Deterministic Area,fig:_spaghetti_Deterministic Angle}. The explicit integrator does exhibit (small) numerical drift in the area associated with this type of scheme as seen in \cref{fig:_spaghetti_Deterministic Area}. 

The Stratonovich RK4-RK4 \cref{Method:Stratonovich integration scheme} has a traced density observed in \cref{fig:PDFStratonovich} aligning with the expected behaviour predicted in \cref{Example:Stratonovich Point Vortex System}, (the points vortices remain in a triangle and rotate around one another and experience North-East South-West stochastic translational motion and a stochastic rotation about the center of vorticity). The explicit integrator preserves the area $3\sqrt{3}/4 \approx 1.047$, contained in the inscribed equilateral triangle to $10^{-7}$, and exhibits a numerical drift in the other direction as compared with the deterministic scheme as seen in \cref{fig:_spaghetti_Stratonovich Area}. Also observed in \cref{fig:_spaghetti_Stratonovich Area} the drift appears realisation dependent. The angle between points $\pi/3 \approx 1.299$ is preserved to $10^{-13}$ and is also Brownian motion dependent. The angle does not exhibit a noticeable directional numerical drift in the ensemble behaviour as seen in \cref{fig:_spaghetti_Area-Process Angle}. This experimental result aligns (aside from small drift) with the expected behaviour of angle and area preservation predicted in \cref{Example:Stratonovich Point Vortex System}.

The Stratonovich system with the Wong Zakai anomaly (introduced in \cref{Example:Type I }) when approximated with \cref{Method:type I}, exhibits a North-West drift in addition to the North-East to South-West stochastic translation, as seen in \cref{fig:PDFWZ and Stratonovich}. \Cref{fig:_spaghetti_Stratonovich with WZ Area} shows the area taking values in $[0,1.5]$ indicating the loss of the initial equilateral triangular configuration. \Cref{fig:_spaghetti_Stratonovich with WZ Angle} shows the area taking values from $[0,\pi]$ indicating the loss of the initial triangular configuration and co-linearity occurring. The interaction between noise and the Wong-Zakai drift caused additional destabilising behaviour as compared with the Wong-Zakai drift on its own \cref{fig:_spaghetti_Stratonovich with WZ Area} to \cref{fig:_spaghetti_Area-Process Area}. Furthermore, the area and angle are not well represented by the mean and variance alone, as seen in \cref{fig:_spaghetti_Stratonovich with WZ Area}, \cref{fig:_spaghetti_Stratonovich with WZ Angle}. The Stratonovich system with the Wong Zakai anomaly drift, experiences periodic fluctuation to the angle as observed in \cref{fig:_spaghetti_Stratonovich with WZ Angle}, this is initially imparted in phase with the area process driven system \cref{fig:_spaghetti_Area-Process Angle}, but grows significantly in magnitude.

The Itô system described in \cref{Example:Itô Point Vortex System} when approximated by \cref{Method:type ito}, has a density in phase space plotted in \cref{fig:PDFItô} similar to that of the Stratonovich system's \cref{fig:PDFStratonovich}, however one can see that the points have become further apart and have a less concentrated density. The dynamic behaviour is more clearly understood by the area \cref{fig:_spaghetti_Itô Area} and angle \cref{fig:_spaghetti_Itô Angle} ensemble plots, showing the equilateral configuration still intact (angle preserved to $10^{-13}$) but the area increasing beyond twice the initial area. 

The Stratonovich system described in \cref{Example:Stratonovich Point Vortex System} when \cref{Method:NLA} is used to approximate Lévy area terms in the equation, has a similar traced density to that of the Stratonovich system observed in \cref{fig:PDFStratonovich NLA}. In \cref{fig:_spaghetti_Stratonovich with NLA Area} the area varies by $5\times 10^{-1}$. In \cref{fig:_spaghetti_Stratonovich with NLA Angle} the angle is varies by $5\times 10^{-1}$. Indicating that the higher order (Lévy area approximating) integrator lost some known dynamical behaviour, despite exploring a similar region in phase space over a short time interval. 

The Area process driven system described in \cref{Example: Type II}, when integrated by the deterministic RK4 scheme \cref{Method:type II}, has a North West motion observed in \cref{fig:PDFArea Process}. The area between the points plotted in \cref{fig:_spaghetti_Area-Process Area} oscillates between 1.22 and 1.3, and the angle oscillates between 0.8 and 1.3 radians see \cref{fig:_spaghetti_Area-Process Angle}. This is specific to the point vortex initial configurations in the commutator of vector fields \cref{fig:Wong_Zakai drift Lie Braket}.

The Canonically lifted FBM driven point vortex system when solved with \cref{Method:FBM}, also appears to exhibit similar structure in phase space as the normally driven point vortex system  as seen in \cref{fig:PDFfbm}. The area and angle preservation by \cref{Method:FBM} occurred at $10^{-7}$ with the same path dependent drifting behavior as the Stratonovich system. 

\subsection{Pathwise comparisons}

In this section we drive the point vortex systems (with exception of the fractional Brownian motion driven system) by the same realisation of Brownian motion. This allows path-wise comparison in solution behaviour.   

\begin{figure}[H]
\centering
\includegraphics[width=220pt]{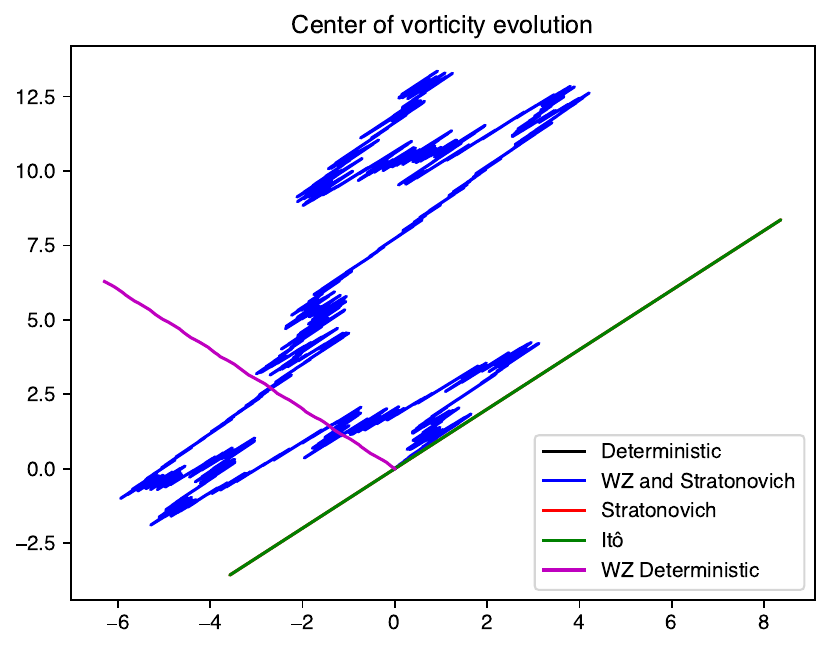}
\includegraphics[width=220pt]{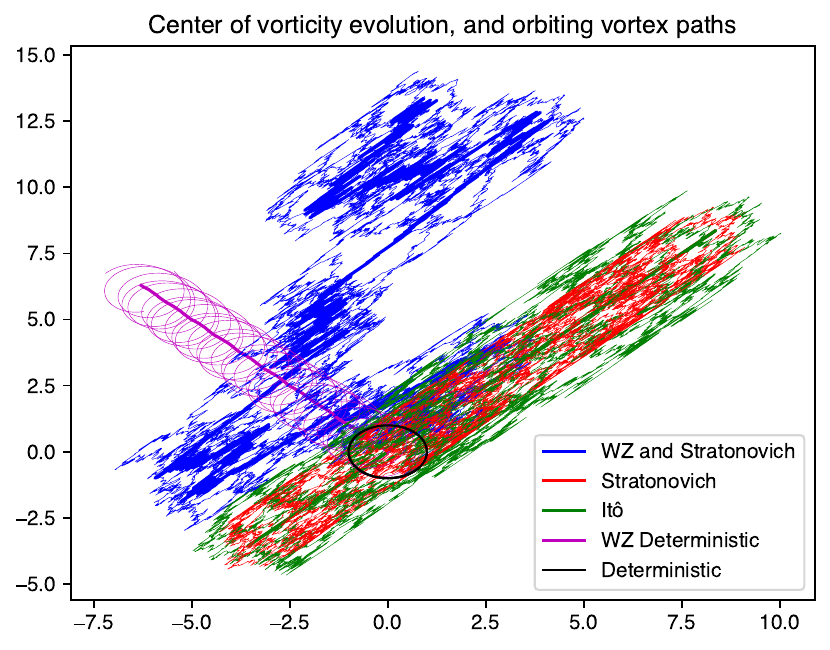}
\caption{Evolution of the cell mean and the trace of the surrounding three
vortices over the time interval respectively. The effect of moving north west is a result of the particular point vortex location in the Wong-Zakai anomaly field, the effect of moving north east is a result of the translational vector field $\xi_2$, and the particular instance of noise. The Stratonovich and Itô have the same center of vorticity and remains on a straight line. The deterministic dynamics remain rotating on the unit circle. }
    \label{fig:Evolution of cell means}
\end{figure}

\begin{figure}[H]
    \centering
\includegraphics[width=80mm]{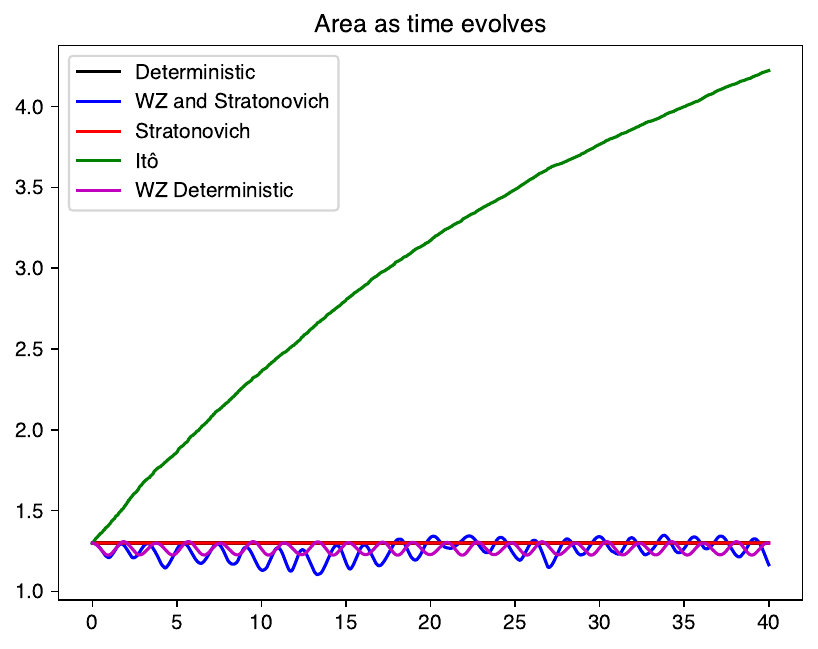}
\includegraphics[width=80mm]{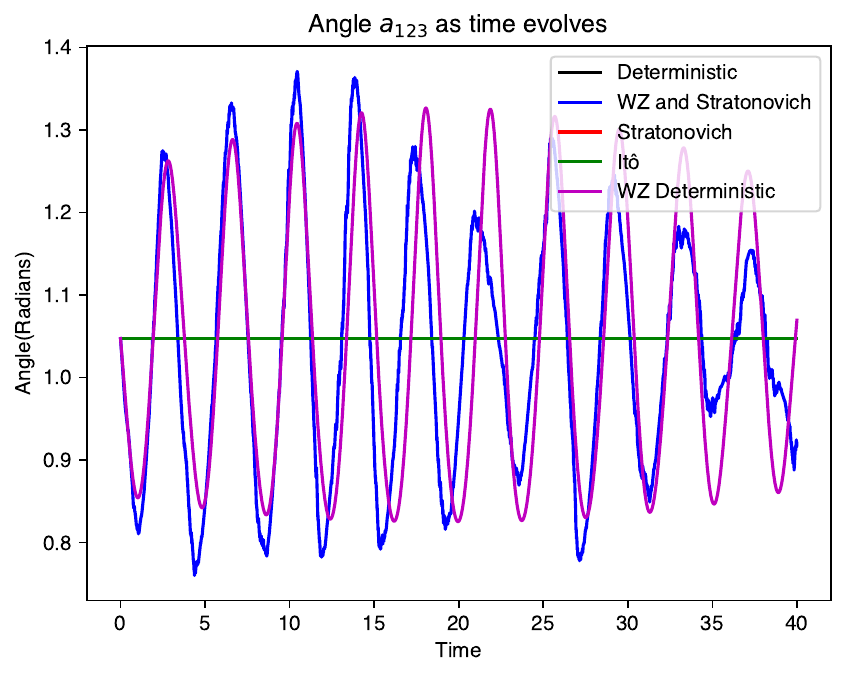}

\caption{Area $A_{123}$ and angle $\alpha_{123}$ between the three point vortices as time evolves. The deterministic scheme preserves angle and area (RK4 \cref{Method:Deterministic} preserved the area in phase space to $10^{-10}$ over the time interval of interest), Stratonovich preserves angle and area (the Stratonovich integrator \cref{Method:Stratonovich integration scheme} preserved the symplectic form (Area to $10^{-8}$), Itô preserves angle but increases the area, the Wong Zakai systems (deterministic and stochastic) both periodically vary in both angle and area, initially in phase.}
\label{fig:Area and Angle}
\end{figure}

\Cref{fig:Area and Angle} contains two figures, plotting the evolution of the angle $\alpha_{123}$ and area $A_{123}$ between the three point vortices $\alpha = 1,2,3$ for \cref{Method:Deterministic}, \cref{Method:Stratonovich integration scheme}, \cref{Method:type ito}, \cref{Method:type I}, \cref{Method:type II}, for the systems introduced in \cref{Example:Deterministic Point Vortex System}, \cref{Example:Stratonovich Point Vortex System}, \cref{Example:Itô Point Vortex System}, \cref{Example:Type I }, \cref{Example: Type II} respectively using the same samples of Brownian motion. \Cref{fig:Evolution of cell means} contains two plots, consisting of the center of vorticity evolution and the trace of the surrounding three point-vortices over the entire time interval using the same samples of Brownian motion.

The deterministic RK4 \cref{Method:Deterministic} preserves the area and angle of the deterministic system introduced in \cref{Example:Deterministic Point Vortex System} (to $10^{-12}$) and the three points remain in a circle about the origin \cref{fig:Area and Angle,fig:Evolution of cell means}. 

The Stratonovich RK4-RK4 \cref{Method:Stratonovich integration scheme} preserves the area and angle between the point vortices (to $10^{-8}$), for the Stratonovich system introduced in \cref{Example:Stratonovich Point Vortex System} in \cref{fig:Area and Angle}. The points and remain in a equilateral triangle following the translating center of vorticity with an additional stochastic rotation on the surrounding point vortices \cref{fig:Evolution of cell means}.

The RK4-EM scheme \cref{Method:type ito} for the Itô system \cref{Example:Itô Point Vortex System}, preserves the angle between the point vortices $\pi/3 \approx 1.299$ (to $10^{-8}$) and the points remain in a triangle in a stochastic moving frame \cref{fig:Area and Angle,fig:Evolution of cell means}. The center of vorticity is in the same location as the Stratonovich system, however the area is not preserved. The Stratonovich-Itô correction caused the point vortex system to move away from the center of vorticity until a larger equilateral triangle configuration is reached \cref{fig:Area and Angle,fig:Evolution of cell means}.

 \Cref{Method:type I} for the Stratonovich with Wong-Zakai drift system introduced in \cref{Example:Type I }, does not conserve angle or area (\cref{fig:Area and Angle}), and there is a drift North-West \cref{fig:Evolution of cell means} due to the position of the vortices in the Wong Zakai anomaly drift field visualised in \cref{fig:Wong_Zakai drift Lie Braket}. 

The \cref{Method:type II} scheme for the deterministic ODE system \cref{Example: Type II} attained in the limit of a ``Area Process'' driven signal. Does not preserve the angle or area \cref{fig:Area and Angle} and also experiences the North-West drift \cref{fig:Evolution of cell means} due to the location of the vortices in the Wong Zakai anomaly drift field. 

\begin{figure}[H]
    \centering
\includegraphics[width=220pt]{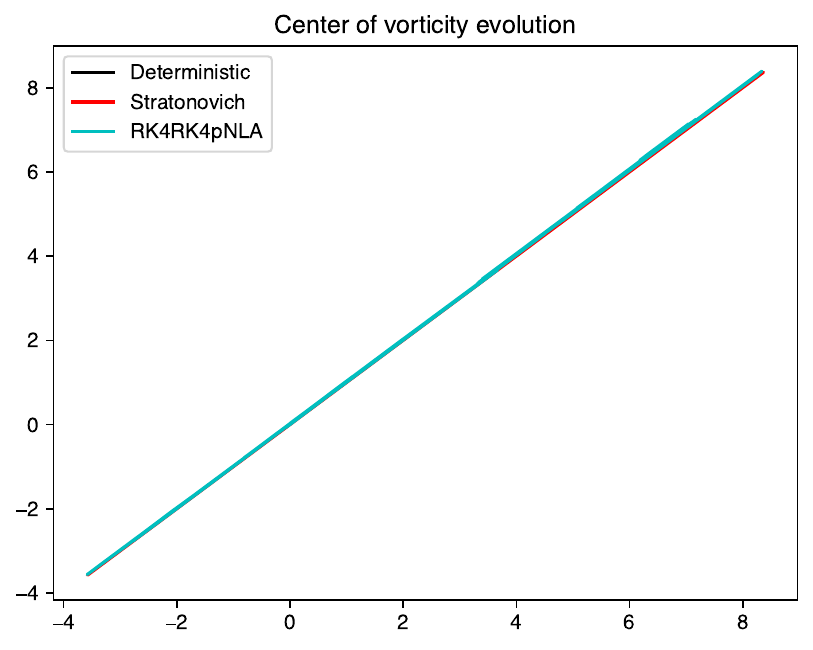}
\includegraphics[width=220pt]{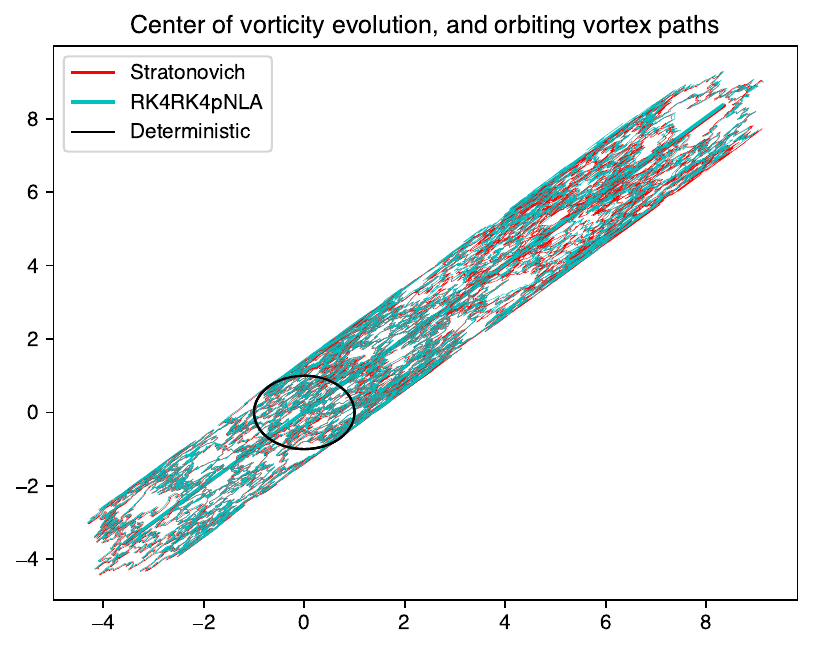}
\caption{Numerical Lévy Areas effect on cell means and orbits. Similar to having a time dependent Wong Zakai anomaly $s_{ij} = J_{i,j}/(\Delta t) $ proportional to nested integration of Brownian motion and inversely proportional to the step size. We see that the addition of the NLA imprints small North-West and South-East motion due to the commutator of vector fields.}
\label{fig:Numerical_Lévy_areacell_means_and_orbits}
\end{figure}

\begin{figure}[H]
    \centering
\includegraphics[width=220pt]{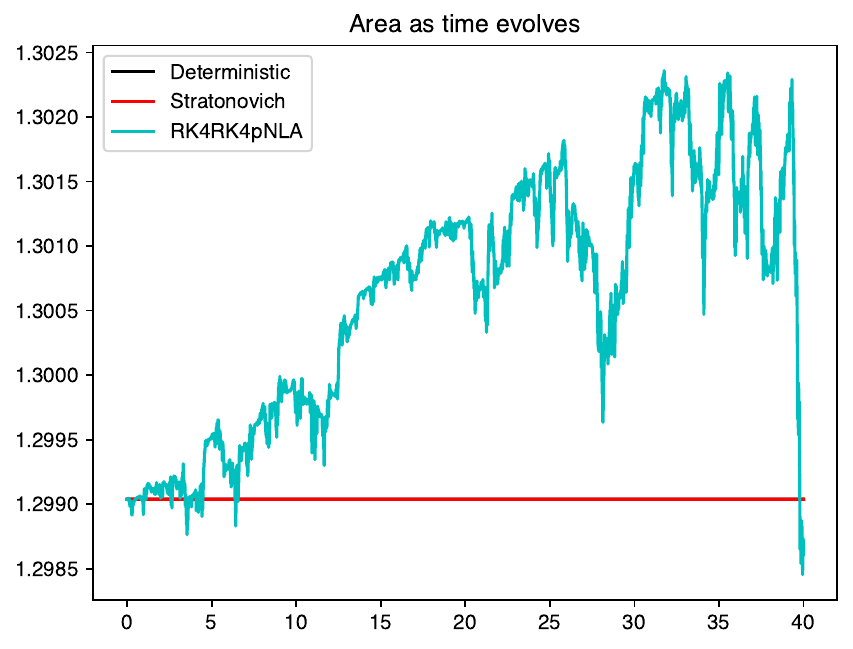}
\includegraphics[width=220pt]{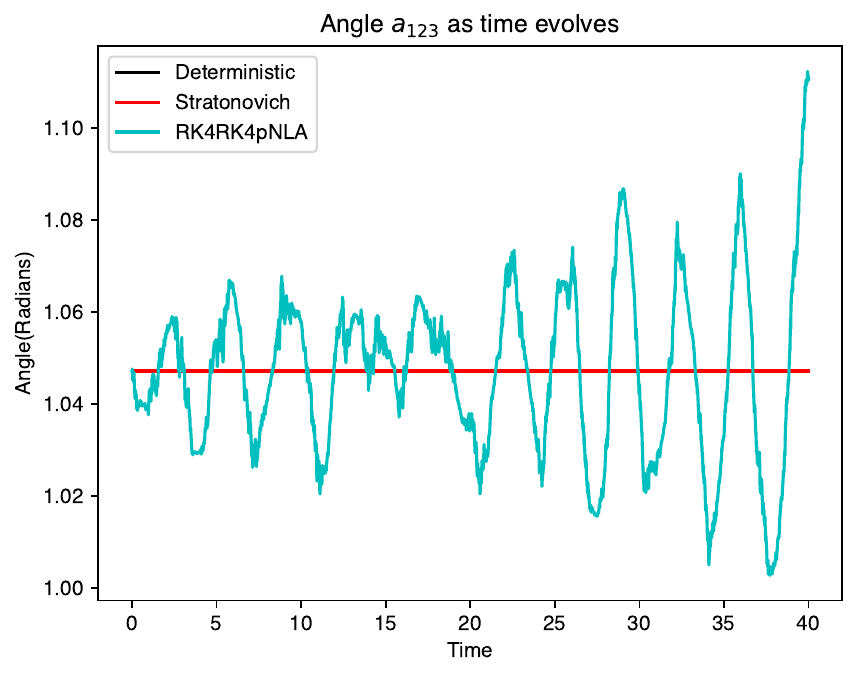}
\caption{Area, and angle as time evolves for the RK4RK4 scheme and the RK4RK4pNLA scheme with the addition of the Numerical Lévy area. The introduction of the higher order terms in the Stratonovich Taylor expansion cause the area and angle between three point vortices to change non-physically.}
    \label{fig:NLA area and angle as time evolves.}
\end{figure}

\Cref{fig:Numerical_Lévy_areacell_means_and_orbits} shows the center of vorticity and orbiting point vortices, for the Deterministic RK4 system, Stratonovich RK4RK4 system and the Stratonovich system where the Lévy area is approximated (RK4RK4 plus Numerical Lévy area Approximation is denoted RK4RK4pNLA in the figure label). Observed is similar pathwise behavior between the Stratonovich and the (RK4RK4pNLA) scheme, with the Lévy area imparting a small North-West time dependent motion to the system. In \cref{fig:NLA area and angle as time evolves.} we plot the area and angle of these same three integrators in the time window $[0,40]$. We observe that the addition of the Lévy in the RK4RK4pNLA scheme caused the loss of the equilateral triangular configuration despite pathwise similar behaviour to the Stratonovich RK4-RK4 integrator observed in \cref{fig:Numerical_Lévy_areacell_means_and_orbits}.

\section{Conclusions} \label{sec: conclusion}
Similar to \cite{cotter2017}, in section \ref{sec:homogenisation} we motivated the modified SALT velocity decomposition in the variational principle \eqref{abstractEP} \cite{DHP2023} via homogenisation theory. 
We showed the Wong-Zakai anomaly and Itô-Stratonovich correction are instances of the additional drift first identified in \cite{cotter2017}. 
From the homogenisation of a deterministic fluid relation, we verified the noise could arise as a limit of time rescaled eigenvalues and the vector fields $\xi$ to be the spacial velocity-velocity correlation eigenvectors as first remarked in \cite{cotter2017, SALTalgo}. Furthermore, we identified the Wong-Zakai anomaly to be the antisymmetric part of the variance covariance matrix of these eigenvalues. The arguments in \cref{sec:homogenisation} provide additional motivation from a modelling viewpoint for a Wong-Zakai anomaly considered in \cite{DHP2023}. Thermal effects are one such example where correlations are introduced in fluid systems, where the power spectrum of observed data is ``red" \cite{TheBasicEffectsofAtmosphereOceanThermalCouplingonMidlatitudeVariability}, such effects may have implications for global ocean circulation \cite{Ditlevsen2023}.

The numerical experiments described in \cref{sec: numerics} support the theoretical predictions of \cref{sec:SPV2D}. In addition, they examine the extent to which integrators agnostic of Lévy area are exhibit desirable solution behaviour. We considered one classic problem in fluid dynamics and showed a wide variation in outcomes regarding solution symmetry, area conservation that was conditional to the Lévy area's inclusion in the numerics.

We interpret our results in two parts. The first is that including higher order terms in a Stratonovich-Taylor expansion in the pursuit of increasing the formal order of accuracy can result in additional non physical behaviour, not associated with the continuum model. The higher order stochastic integrator adopted in this paper, adds an approximation of the Lévy area to an existing Additive Runge Kutta method. This approach captures higher order (symmetric) terms than a Milstein type discretisation whilst still including the antisymetric approximation of the Lévy area. As shown in our experiments, the inclusion of these ``typically small" nested integral terms had little effect in the short time pathwise behavior of the system, yet we observed a loss in some aspect of structure preservation, namely the conservation of area and angle. 


The second, parallel, conclusion is that the Wong-Zakai anomaly, which is of higher order, produces substantial dynamical effects which must not be discarded. The presence of Wong-Zakai anomalies, which are present generically in systems from homogenisation theory or as smooth approximation of SDE's, have a dramatic effects on solution behaviour. This is particularly relevant to geophysical simulations with SALT arising from ``real" data, we do not expect the perfect scale seperation corresponding to the $\varepsilon = 0$ homogenised limit, rather one observes deterministic velocities closely approximated by a realisation of a SDE:

$$u_t(x^\varepsilon_t) = \dot{g}_t{g}^{-1}_t(x^\varepsilon_t) = \dot{x}^\varepsilon_t =  \frac{1}{\varepsilon}\lambda^{k,\varepsilon}_t v_k(x) + \operatorname{Ad}_{\widetilde{g}_{t/\varepsilon}}\overline{u}_t (x^\varepsilon_t) \approx U(x) \mathrm{d}t + \xi(x) \circ \mathrm{d} \widetilde{\boldsymbol{B}},$$

where the Brownian rough path $\widetilde{\boldsymbol{B}}_{t.s} = (B_t - B_s, \mathbb{B}^{\text{Strat}}_{t,s} + \mathsf{s}(t-s))$ has perturbed Lévy area. One can choose to represent a RDE driven by $\widetilde{\boldsymbol{B}}_{t.s}$ as an SDE with a Wong-Zakai anomaly in its drift, however proceeding as in \cite{SALTalgo} one may take $u_t(x^\varepsilon_t)$ as a ``fine scale truth" and only have access to the near noise term $\lambda^{k,\varepsilon}_t v_k(x) \approx \xi(x) \circ \mathrm{d} \widetilde{\boldsymbol{B}}_t$ with no knowledge of how to compute commutators $[\xi_i, \xi_j]$ or procedure to extract $\mathsf{s}$ from iterated integrals of $\lambda^{k,\varepsilon}_t$. Therefore, if one adopts the approach of SALT using Brownian paths with modified Lévy area, it is necessary to compute the iterated integrals to accurately capture the effects of $\mathsf{s}$ (cf. Equation \eqref{EQ: Stratonovich Taylor Expansion}). Indeed, failing to do so ignores physical effects due to changing in the Hamiltonian of the system \cite{DHP2023}, these appear at a lower order than that of the numerical Lévy area.

The example tested in this paper of point vortex dynamics in \ref{sec:SPV2D} is a particularly tractable for analytic computations, where the destabilising effect of the Wong-Zakai anomaly is easily noted. We hypothesise that for fluid models of much higher complexity than the Euler point vortex system, the Wong-Zakai anomaly produces continuum analogues of the chaotic effects seen here. 


\begin{figure*}
\centering
    \begin{subfigure}[t]{0.295\textwidth}
\centering
\includegraphics[width=.95\textwidth]{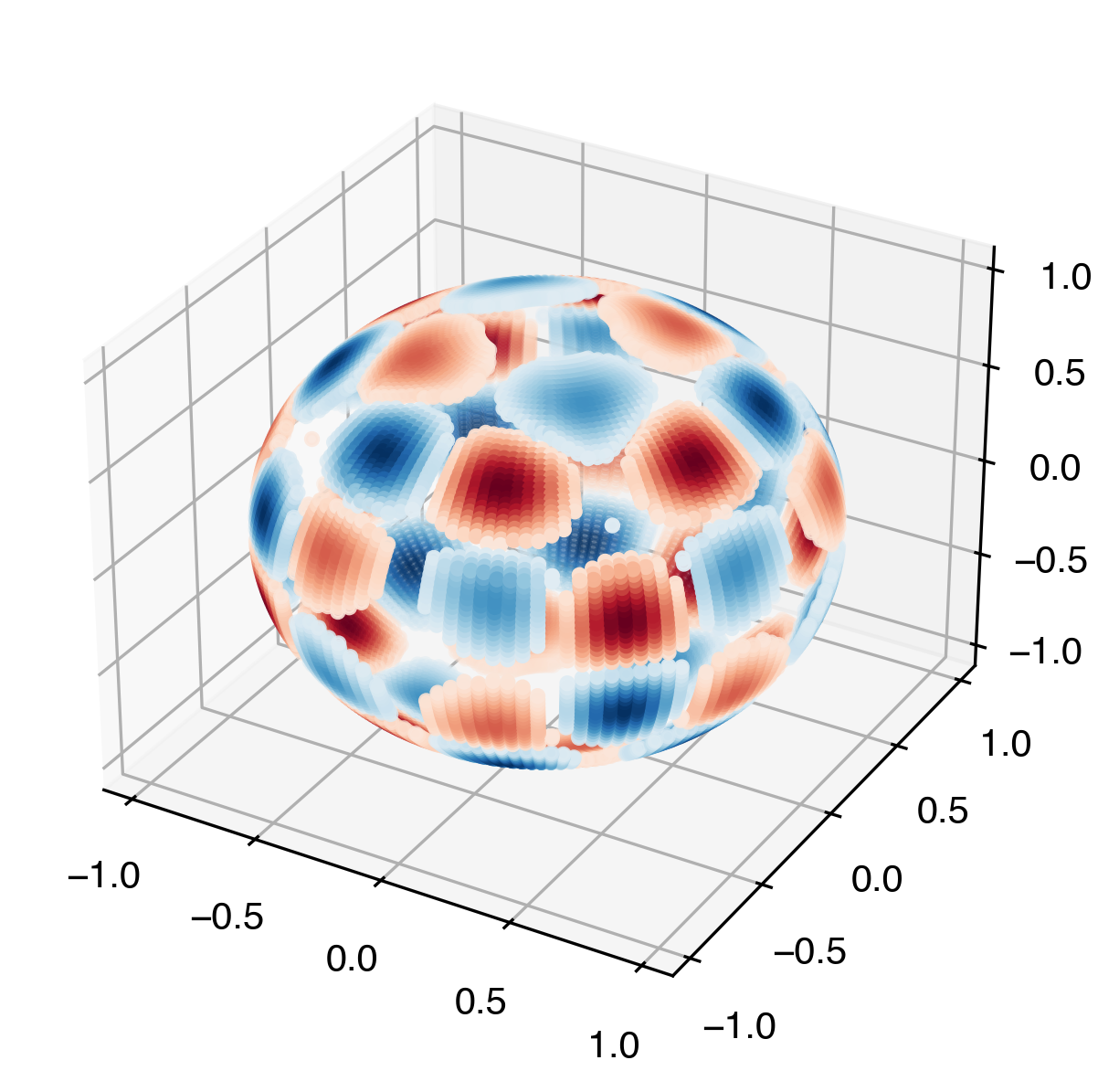}
\label{fig:time zero many points}
\end{subfigure}
\begin{subfigure}[t]{0.295\textwidth}
\centering
\includegraphics[width=.95\textwidth]{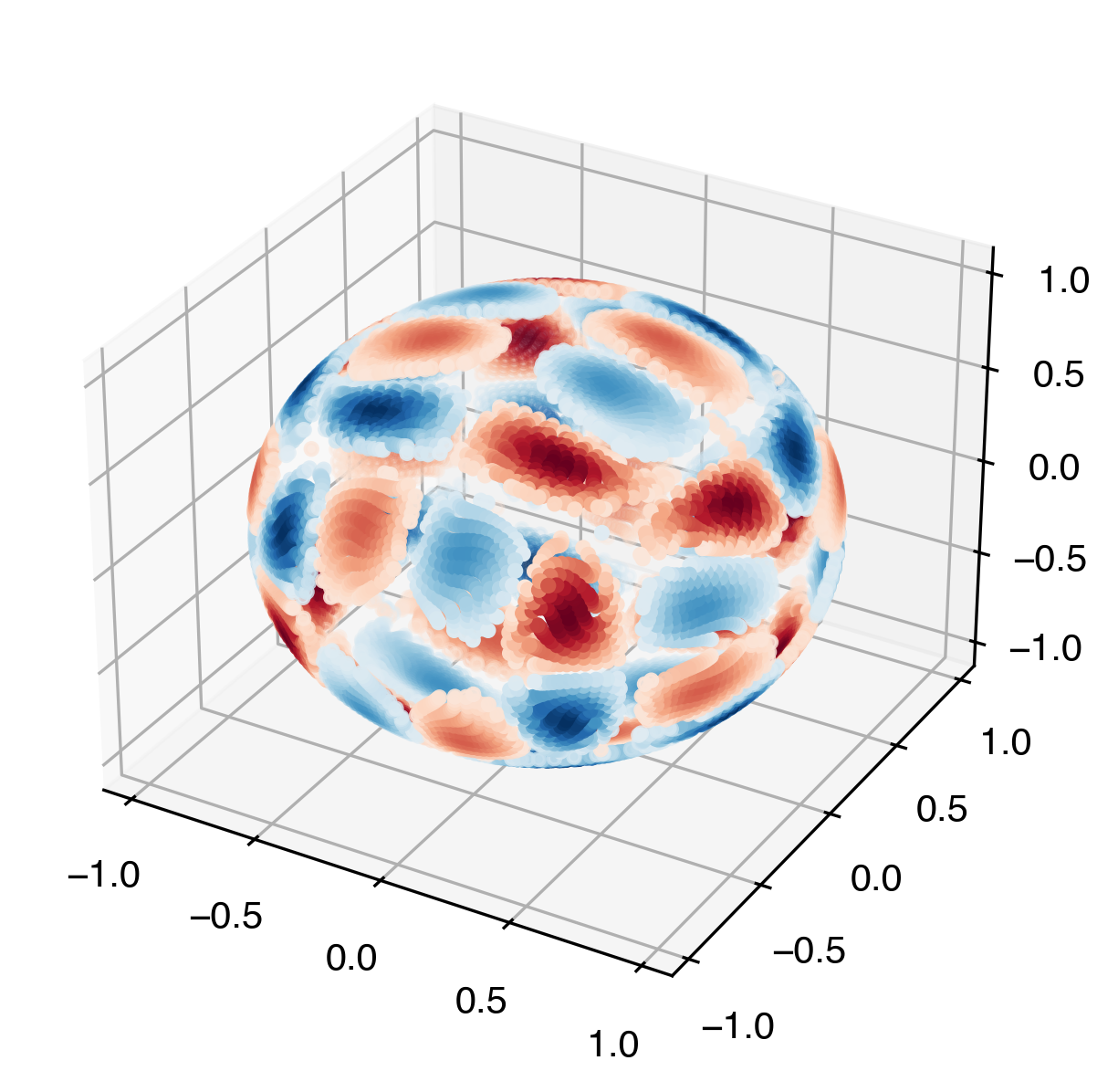}
\end{subfigure}
\begin{subfigure}[t]{0.295\textwidth}
\centering
\includegraphics[width=.95\textwidth]{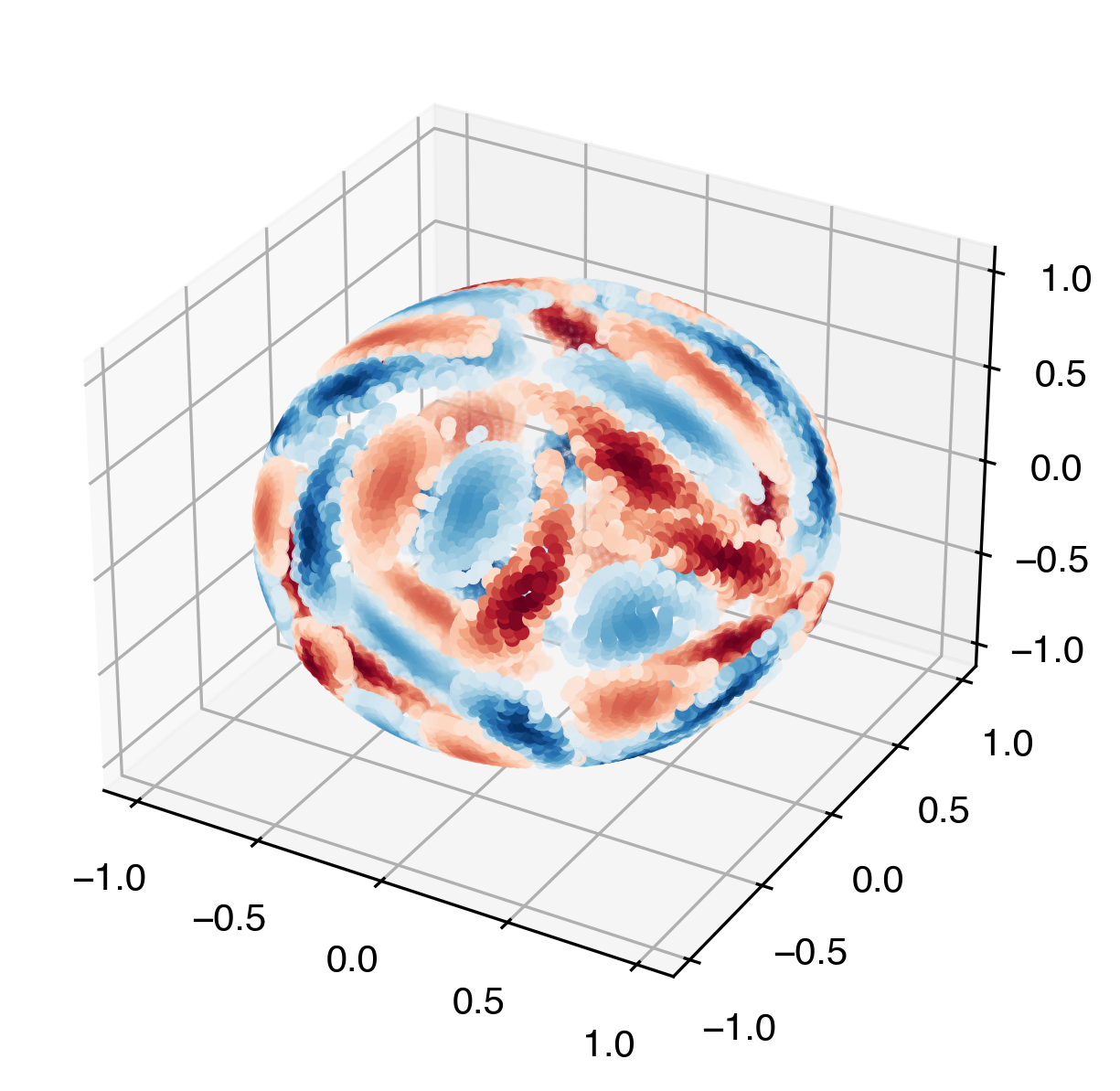}
\end{subfigure}
\caption{Shown is the solution of a regularised (vortex blob type) system constrained on the surface of the sphere. Red denotes positive potential vorticity, and blue denotes negative potential vorticity. The initial conditions are a projection of the initial conditions in \cite{SALTalgo}, where vortex blobs with little vorticity are removed. }
\label{fig:time zero many points}
\end{figure*}

It has been shown that both 2D deterministic Navier Stokes and stochastic Euler arise under different schemes as a mean field limit of an $N$ stochastic point vortex system \cite{fl2021mean, rosenzweig2020meanfield}. It has been shown numerically that for generic initial conditions the long time behaviour of the Euler equations on the sphere is governed by ``vortex blobs" that obey point vortex dynamics \cite{modin_viviani_2020}. One could conjecture if the Wong-Zakai anomaly for 2D spherical Euler can disrupt the behaviour in a similar vein to what was shown in this paper.

As a final concluding remark we note that the interplay of the two key ideas, composition of maps and homogenisation theory, can be independently generalised for the purposes of this paper. The composition of maps approach to geometric mechanics has been examined in recent works of \cite{HHS2023, HHS2023a} and is a classical approach to fluid dynamics with mean and fluctuating components \cite{Holm1996, andrews_mcintyre_1978} as well as in the context of the generalised Lagrangian mean (GLM) theory \cite{Holm2002}. In the context of this paper, we assume a fast fluctuating map that provided a means to produce noise via homogenisation theory and induce a sum of mean and fast terms in the velocity field. 

A natural extension to the ideas of this paper in terms of composition of maps is to establish an analogy to the tangled Poisson structures and momentum shifts seen in \cite{HHS2023a} caused by a fast flow. The interpretation of noisy velocity constraint in SALT as arising from a hybrid deterministic variational principle could provide additional understanding of  how SALT relates to the underlying deterministic physics one wishes to model at a coarser or unresolved scale. 
\newline 

In other directions, homogenisation extends to situations much more general than discussed in this paper. Most relevant to this work is the extension of the weak invariance principle (WIP) (discussed in \cref{sec:appendixh}) for deterministic homogenisation, which we utilised as a means to generate Brownian motion. More specifically, the WIP is what permitted us to pass to the limit and obtain a stochastic equation in \ref{sec:homogenisation} to derive the SALT velocity. Much like how stochastic homogenisation theory motivates the SALT ansatz, one could seek to motivate the rough variational principle seen in \cite{crisan2022variational} via an analogous \emph{rough} homogenisation, a natural first choice of example being fractional Brownian motion to incorporate memory effects. A key ingredient to this approach is to consider weak invariance principles that converge to rough paths more general than Brownian motion. The theory of rough geometric mechanics established in \cite{crisan2022variational} would then give a generalisation of the Euler-Poincaré equations in \eqref{abstractEP} with other rough paths. A selection of weak invariance principles more general than Brownian have been examined in \cite{Gehringer2022}, relating to Hermite processes and fractional Brownian motion. We conjecture a relevant extension in the appendix, motivating \cref{Method:FBM} from our numerics.

\paragraph{Acknowledgements.}

We remain gracious to our colleagues for their stimulating discussions and pleasant company. In this paper, we'd especially like to thank S. Takao,  R. Hu, O. Street, JM. Leahy, D.D. Holm, G. Pavliotis and the Friday afternoon Geometric Mechanics Research Group.
The work of TD is supported by an EPSRC PhD Scholarship. The work of JW is supported by the European Research Council (ERC) Synergy grant “Stochastic Transport in Upper Ocean Dynamics” (STUOD) – DLV-856408.

\newpage 
\bibliographystyle{plainurl}
\bibliography{bibliography.bib}
\newpage
\appendix
\section*{Appendices}
\addcontentsline{toc}{section}{Appendices}
\renewcommand{\thesubsection}{\Alph{subsection}}

\subsection{Stochastic Homogenization}\label{sec:appendixh}

In this appendix we shall overview some formulae and theorems applied in \cref{sec:homogenisation} from homogenisation theory, a mathematical mechanism from which diffusion processes arise as the limit of differential equations. In this section we outline the general theory and present a literature review. The systems we consider consist of slowly evolving dynamics coupled to an underlying ``fast" equation involving stochasticity. The scale separation of ``slow" verses ``fast" is governed by a parameter $\varepsilon$. The objective is elimination of the parameter and the fast dynamics through the derivation of an accurate limit dynamics when $\varepsilon \rightarrow 0$.  

The aforementioned fast equations typically fall into two cases: deterministic ODE with randomly distributed initial conditions or stochastic differential equations. In the first case each integral curve produces a random variable taking values in an appropriate function space depending on the realisation of the intial value. In the second case one formally views each realisation of the driving noise (Brownian motion in most, but not all situations) as giving a particular sample path of the solution of the SDE.

The results we discuss and state shall be kept simplified but sufficiently general enough to apply to our application of interest, which is to motivate the noise ansatz in \cref{sec:continuumeuler}, \ref{sec:SPV2D} with room for extension to the ideas presented in \cite{crisan2022variational}. 
Many results in this paper and the overarching framework of geometric mechanics are ``modular" in the sense that many general choices of fast dynamics can be coupled to the deterministic reconstruction equation in order to arrive at an applicable RDE for the Lagrangian paths that result in rough Euler Poincar\'e equations.

We describe the structure of a typical (finite dimensional) homogenization problem. Let $\mathcal{X}, \mathcal{Y}$ be spaces that we take to be either Euclidean or torii for some (possibly equal) dimensions $d, m$. Given smooth vector fields $f_i : \mathcal{X} \times \mathcal{Y} \rightarrow \mathcal{X}$ we are interested in the small $\varepsilon$ limit of the ODE system:

\begin{equation}\label{xeps}
\dot{x}^\varepsilon_t = \frac{1}{\varepsilon}f_0(x^\varepsilon_t,y^\varepsilon_t) + f_1(x^\varepsilon_t,y^\varepsilon_t), \qquad x^\varepsilon_0 = v. 
\end{equation}

Equation \eqref{xeps} can be thought of as a rescaling of the following ODE

\begin{equation}\label{originalxeps}
\dot{x}_{\varepsilon}(t) = \varepsilon f_0(x_{\varepsilon}(t),y(t)) + \varepsilon^2 f_1(x_{\varepsilon}(t),y(t)), \qquad x_{\varepsilon}(0) = v, 
\end{equation}

under the change of time $t \mapsto t /\varepsilon^2$ and defining $x^\varepsilon_t = x_{\varepsilon}(\varepsilon^{-2} t), \,\, y^\varepsilon_t = y(\varepsilon^{-2} t)$. We assume $y$ has dynamics defined by a flow map of a continuous time dynamical system $\phi_{t}: \mathcal{Y} \rightarrow \mathcal{Y}$. It is possible to allow full coupling via dependence on $x_{\varepsilon}$, i.e. $\phi^{x_{\varepsilon}}_{t}$ \cite{Engel2021} but for now we shall keep in mind dynamical systems of the form $\dot{y}(t) = h(y(t))$. We assume the $y$ dynamics are Ergodic and possess a unique invariant measure $\nu$. For an initial condition $y(0) \in \mathcal{Y}$ we wish to think of $\phi_{t}(y(0))$ as a stochastic process. With this assumption in mind it would be intuitive to rewrite the system \eqref{xeps} as an integral equation in the time $\tau = t / \varepsilon^2$,

\begin{equation}\label{inteqn}x^\varepsilon_t=v+\varepsilon \int_0^{\frac{t}{\varepsilon^2}} f_0(x_{\varepsilon}(\tau), y(\tau)) d \tau + \varepsilon^2 \int_0^{\frac{t}{\varepsilon^2}} f_1(x_{\varepsilon}(\tau), y(\tau)) d \tau.\end{equation}


The scaling in \eqref{inteqn} is such that in the $\varepsilon \rightarrow 0$ limit the second term ought to approach a $\nu$ average of $f_1$ via the Ergodic theorem and the first term should approximate Brownian motion through the central limit theorem. One postulates that for small values of $\varepsilon$ solutions of \eqref{xeps} are well approximated by the solutions of some ``homogenised" stochastic equation of the form:  

\begin{equation}\label{homogenisedguess}\mathrm{d}X = \overline{f}_1(X) \mathrm{d}t + \Sigma(X) * \mathrm{d} B_t.\end{equation}

Which has eliminated dependence on the fast variables but accounted for their effect during the averaging process. It is unclear what interpretation of the stochastic integral ought to be correct for \eqref{homogenisedguess}, thus we write $* \mathrm{d}B$ as a placeholder for Itô, Stratonovich or otherwise. The aim of homogenisation theory is understanding the correct form of \eqref{homogenisedguess}, via a suitable proof of convergence for the $y$ dynamics in \eqref{yeps} to a difussion process and interpreting the notion of the stochastic integral $* \mathrm{d} B_t$. 

There are two possibilities for introducing randomness into the fast flow (and thus \eqref{xeps}). Either the initial condition $y_0 \in \Omega$ is drawn randomly from a probability space and $\phi_{(\cdot)}$ is the flow of a sufficiently chaotic ODE, or $\phi_{(\cdot)}$ is the stochastic flow of an Ergodic diffusion process. 

The former case is known in the literature as deterministic homogenization \cite{kelly2017chaos}, while the latter has been referred as stochastic averaging or diffusive homogenisation \cite{Gehringer_2021}. Both fall under what has also been coined diffusion creation and homogenization theory \cite{Papanicolaou_1995}. With the introduction of the theory of rough paths by Lyons \cite{Lyons1998} \emph{rough} homogenisation has been introduced as a new tool to show convergence of the fast dynamics to noise \cite{kelly2017chaos, Gehringer_2021}. Both the stochastic averaging and deterministic homogenisation perspectives can represent the fast $\mathcal{Y}$ dynamics as being governed by an equation of the type:

\begin{equation}\label{yeps}
\dot{y}^\varepsilon_t = \frac{1}{\varepsilon^2}h(x^\varepsilon_t,y^\varepsilon_t) + \frac{1}{\varepsilon}b(x^\varepsilon_t,y^\varepsilon_t)\eta_t, \qquad y_\varepsilon(0) = y_0. 
\end{equation}

Where $h:  \mathcal{X} \times \mathcal{Y} \rightarrow \mathcal{Y}$ and $b:  \mathcal{X} \times \mathcal{Y} \rightarrow L(\mathcal{Y}; \mathcal{Y})$. $\eta_t(\omega) \in \mathcal{Y}$ will denote some form of noise. For non identically vanishing $b$ we are in the regime of stochastic averaging, whereas $b \equiv 0$ is deterministic homogenisation.


We may now understand equation \eqref{homogenisedguess} through three different perspectives. Firstly one can think of the slow $x_\varepsilon$ dynamics as a smooth approximation of a SDE, with the fast variable $y_\varepsilon$ treated as mollified white noise. Passing to the limit we return to a Brownian motion. This is the approach taken in \cite{b0b66424-1713-373e-87d8-fb6774b05aa3} and is deeply related to the approximation theory of stochastic differential equations \cite{ikeda2014stochastic}. In the context of geometric mechanics the approximation of stochastic Lagrangian trajectories by coloured noise was considered in \cite{DHP2023} and produced a neither Itô nor Stratonovich interpretation of the noise. In \cref{sec:homogenisation} homogenization produced these same correction terms and this is a consequence of the fact that the white noise limit in \cite{DHP2023} is a special case of stochastic averaging that introduces a shift to the L\'evy area, namely, \eqref{yeps} with $h(x^\varepsilon_t, y^\varepsilon_t) = Ax^\varepsilon_t, b(x^\varepsilon_t, y^\varepsilon_t) \equiv D, \eta_t = \dot{W}_t$ for constant positive definite matrices $A,D$. 

Homogenisation can also be understood via the method of multiple scales and perturbation expansions of the Feller (infinitesimal) generator \cite{Papanicolaou_1995}. An ansatz for the limiting SDE can be formally calculated by writing down the Kolmogorov equations for the coupled system \eqref{xeps}, \eqref{yeps} and deriving a generator for \eqref{homogenisedguess} which (weakly) determines the correct form of the equations. One can then prove convergence of solutions using the notion of two-scale convergence. This is the approach outlined in \cite{Pavliotis2008} \cite{Dror_2004}. The same neither Itô nor Stratonovich correction terms discussed in the context of the previous approach appear as additional leading order terms contributed by the $f_1$ component of the generator. A strength of this approach is that one can derive explicit representations of the drift and diffusion coefficients for \eqref{homogenisedguess} in terms of the functions $f_i, h, b$ and the measure $\nu$.

The third and most recent approach is to use rough path theory. The slow ODE is reformulated as a rough differential equation driven by smooth approximation to Brownian motion. We can demonstrate the idea of the rough paths approach in the simple case of $f_1(x^\varepsilon_t , y^\varepsilon_t) = f_1(x^\varepsilon_t)$ and $f_0(x^\varepsilon_t, y^\varepsilon_t) = w(x^\varepsilon_t)l(y^\varepsilon_t)$. We define the approximation $B_\varepsilon(t; y_0) = \varepsilon\int_0^\frac{t}{\varepsilon^2} l(\phi_s(y_0))\mathrm{d}s$, since $B_\varepsilon(t)$ is certainly $C^1$ (it is the integral of a continuous function) one can define the signature tensor $\mathbb{B}_\varepsilon(s,t) = \int_0^t (B_\varepsilon(r) - B_\varepsilon(s)) \otimes \mathrm{d}B_\varepsilon(r)$ as iterated integrals canonically in the Riemann sense. It suffices to prove convergence of $(B_\varepsilon, \mathbb{B}_\varepsilon)$ in the rough paths topology. The convergence of the driving signal $(B_\varepsilon, \mathbb{B}_\varepsilon)$ in this simplified case determines the form of \eqref{homogenisedguess} with $\overline{f}_1 = f_1$, $\Sigma = w$ and interpretation of $* \mathrm{d}B$ decided by the signature of the limit rough path. Convergence of the random element $B_\varepsilon(t;y_0)$ to a rough path lift of Brownian motion is known as an iterated weak invariance principle (WIP), as we are applying a form of the functional central limit theorem (the weak invariance principle) to $B_\varepsilon$ and it's iterated integrals. The iterated WIP is a precise formulation of the intuition seen from \eqref{inteqn}. The power of this technique is the ability to utilise the rough paths topology to show convergence of $(B_\varepsilon, \mathbb{B}_\varepsilon)$ and then continuity of the solution map for rough differential equations to conclude convergence of solution $x^\varepsilon \rightarrow x$. Note that the precise understanding the convergence of $(B_\varepsilon, \mathbb{B}_\varepsilon)$ resolves the ambiguity of what the placeholder $*\mathrm{d}B$ ought to resemble, as the notion of ``integration with respect to $\mathrm{d}B$" is encoded by what form the limit of the signature tensor $\mathbb{B}^\varepsilon$ takes. By reducing the proof of convergence to \eqref{homogenisedguess} into proving an iterated WIP, rough paths theory allows weaker assumptions than the central limit theorem. In some cases the rough driver in \eqref{homogenisedguess} and interpretation of noise $\eta_t$ in \eqref{yeps} can be generalised, for example fractional Brownian motion, provided an iterated WIP can be constructed, such as in \cite{Gehringer2022}. 

\subsubsection{Weak Invariance Principles}\label{sec:WIPs}

We shall now state a weak invariance principle that is sufficient to obtain Stratonovich noise and a Wong-Zakai anomaly, invoked in \cref{sec:homogenisation} and motivating the type of noise in section \ref{Example:Stratonovich Point Vortex System} and \ref{Example:Type I } in this paper. We shall also briefly discuss extensions for more general fractional Brownian noise allowing for memory effects, along with their own Wong-Zakai anomaly. 

We make use of the notation of \cite{kelly2017chaos} for $f: \mathcal{X} \times \mathcal{Y} \rightarrow \mathbb{R}$ depending on slow and fast dynamics. 
$$\|f\|_{C^{r, \kappa}(\mathcal{X} \times \mathcal{Y})}=\sum_{|k| \leq\lfloor r \rfloor} \sup _{x \in \mathcal{X}}\left\|D^k b(x, \cdot)\right\|_{C^\kappa}+\sum_{|k|=\lfloor r \rfloor} \sup _{x, x^\prime \in \mathcal{X}} \frac{\left\|D^k b(x, \cdot)-D^k b(x^\prime, \cdot)\right\|_{C^\kappa}}{|x-x^\prime|^{r -\lfloor r \rfloor}}$$

Where $r \geq 0, \kappa \in [0,1)$, $C^\kappa$ denotes the standard Hölder norm in the variable $y \in \mathcal{Y}$,  and the second sum is not included if $r \in \mathbb{N}$. Let us also denote $C^{r+, \kappa}(\mathcal{X} \times \mathcal{Y}) := C^{r + \delta, \kappa}(\mathcal{X} \times \mathcal{Y})$ for all $\delta^\prime > \delta > 0$ for some $\delta^\prime > 0$ and $C^{\kappa}(\mathcal{X} \times \mathcal{Y}) := C^{0,\kappa}(\mathcal{X} \times \mathcal{Y})$. For $\mathbb{R}^d$ valued $f$ we define $\| \cdot \|_{C^{r, \kappa}(\mathcal{X} \times \mathcal{Y}, \mathbb{R}^d)}$ by the sum of the scalar norms of the components of $f$. 

We shall be strict in writing $C^{\kappa}(\mathcal{X} \times \mathcal{Y}) \neq C^{\kappa}$ to prevent the mixed Hölder norm defined above being misinterpreted with the standard Hölder norm for single variable functions also used in this section, (such as for rough paths defined on $[0,T]$, or univariate functions on $\mathcal{Y}$) where the common definitions of Hölder continuity are used. 

\begin{theorem*}[Deterministic homogenisation theorem. Kelly, Melbourne 2017 \cite{kelly2017chaos}]\label{KM2017}
Let $\mathcal{Y} = M$ any smooth manifold of finite dimension with $\phi_t: M \rightarrow M$ a smooth flow. Let $\Lambda \subset M$ be a hyperbolic basic set of the dynamics \footnote{Recall a $\phi_t$ invariant set $\Lambda$ is called a hyperbolic basic set if $T_x M$ splits into stable and unstable spaces in the usual dynamical systems sense for all $x \in \Lambda$, $\phi_t$ is topologically transitive and periodic orbits are dense in $\Lambda$.} with $\nu$ its ergodic equilibrium measure. We denote $\kappa$-Hölder continuous real valued functions with vanishing $\nu$ mean by $C^{\kappa}_\nu(\Lambda)$. Define the following bilinear operator:

$$\mathfrak{B}(u,v) = \lim_{n \rightarrow \infty}\frac{1}{n}\int_\Lambda \int_0^n \int_0^t u \circ \phi_s  v \circ \phi_t \mathrm{d}s \mathrm{d}t, \quad u, v \in C^\kappa_\nu(\Lambda)$$

Then the following hold: 

The operator $\mathfrak{B}: C^\kappa_\nu(\Lambda) \times C^\kappa_\nu(\Lambda) \rightarrow \mathbb{R}$ is well defined as a convergent limit, bounded and positive semidefinite. Given $f_1 \in C^{1+,0}(\mathbb{R}^d \times M, \mathbb{R}^d), f_0 \in C^{2+, \kappa}_\nu(\mathbb{R}^d \times M, \mathbb{R}^d)$, there exists drift and diffusion coefficients $\overline{f}_1, \Sigma \Sigma^T$ defined by $\overline{f}_1(x)^i = \langle f_1^i(x, \cdot) \rangle + \mathfrak{B}(f_0^k(x, \cdot), \partial_k f_0^i(x, \cdot))$, $\left(\Sigma \Sigma^T\right)^{ij} = \mathfrak{B}(f_0^i(x, \cdot), f_0^j(x, \cdot)) + \mathfrak{B}(f_0^j(x, \cdot), f_0^i(x, \cdot))$ such that the family $x_\varepsilon$ converge in distribution, and in the sup norm topology to the unique solution $X$ of SDE \eqref{homogenisedguess} with $*\mathrm{d}B$ interpreted in the Itô sense. 
In addition, given the following integral converges:

$$\int_0^\infty \int_\Lambda u v \circ \phi_t \mathrm{d}\nu \mathrm{d} t < \infty$$

Then in addition we have the following formulae for $\mathfrak{B}(u,v)$

$$\mathfrak{B}(u,v) = \int_0^\infty \int_\Lambda u v \circ \phi_t \mathrm{d}\nu \mathrm{d} t$$
\end{theorem*}

From the above theorem, we have a representation of the coefficients of the homogenised SDE in terms of the the operator $\mathfrak{B}$. A key result to obtain the homogenisation theorem is the following abstract weak invariance principle. 

\begin{theorem*}[Iterated WIP. Kelly, Melbourne 2016 \cite{kelly2016}]\label{KM2014}
Let $M, \phi_t : M \rightarrow M, \nu, \Lambda, \mathfrak{B}$ be as in the deterministic homogenisation theorem. For $l \in L^\infty(\Lambda, \mathbb{R}^d)$ with $\nu$ mean zero define:

$$B_{n}(t)=n^{-1 / 2} \int_0^{t n} l \circ \phi_s \mathrm{d} s, \quad \quad \mathbb{B}_{n}(t)=n^{-1} \int_0^{t n} \int_0^s l \circ \phi_r \otimes l \circ \phi_s \mathrm{d} r \mathrm{d} s$$

Then as $n \rightarrow \infty$ we have convergence in distribution in the space $C([0,\infty), \mathbb{R}^m \oplus \mathbb{R}^m \otimes \mathbb{R}^m)$ with the $\kappa$-Hölder topology\footnote{The $\kappa$-Hölder seminorm for two parameter increments $\|X\|_{C^\kappa [0,T]} := \sup_{s \neq t \in [0,T]}\frac{|X_t - X_s|}{|t-s|^\kappa}$ induces a topology on the signature (hence rough paths) by defining $\|\mathbb{X}\|_{2 \kappa} = \sup_{s \neq t \in [0,T]}\frac{\|\mathbb{X}_{t,s}\|}{|t-s|^{2\kappa}}$ and $\|(X, \mathbb{X})\|_\kappa := \|X\|_{C^\kappa [0,T]} + \|\mathbb{X}\|_{2 \kappa}$.} for some $\kappa \in (\frac{1}{3},1)$.

$$(B_{n}, \mathbb{B}_{n}) \overset{d\,\,\,}{\rightarrow} (B, \widetilde{\mathbb{B}})$$

Where $B$ is an $m$ dimensional Brownian motion with covariance matrix $M$ and $\widetilde{\mathbb{B}}^{ij}(t) = \int_0^t B^i_s \mathrm{d}B^j_s + \mathfrak{B}(l^i, l^j)t$ is the Itô enhancement of Brownian motion with a perturbation to its signature.

Furthermore we have $M^{ij} = \operatorname{Sym}(\mathfrak{B}(l^i,l^j)), (\mathsf{s}^\prime)^{ij} = \operatorname{Alt}(\mathfrak{B}(l^i,l^j))$.
\end{theorem*}

\begin{remark}\label{rmkWIP} For a standard Brownian motion $W$ lifted as a Stratonovich rough path $\mathbb{W}^{\text{Strat}}_{s,t} -\frac{1}{2}I_{m \times m}(t-s) = \mathbb{W}^{\text{Itô}}_{s,t}$. It can be shown that if $DD^T = M$ is the Cholesky decomposition of the covariance matrix, then defining $B = DW, \mathsf{s}^\prime = D\mathsf{s}D^T$ allows writing the signature as: 
\begin{align*}
\widetilde{\mathbb{B}}_{t,s} &= D \otimes D (\mathbb{W}^{\text{Itô}}_{s,t}) + M(t-s) + \mathsf{s}^\prime(t-s) \\
&= D \otimes D (\mathbb{W}^{\text{Strat}}_{s,t}) - \frac{1}{2}D \otimes D (I_{m \times m})(t-s) + M(t-s) + \mathsf{s}^\prime = D\otimes D(\mathbb{W}^{\text{Strat}}_{s,t} + \mathsf{s}(t-s)),
\end{align*}

 where have made use of the tensor product identity $M \otimes N \left(X\right) = MXN^T$ holding for any $M \in \mathbb{R}^{n_1 \times n_1}, N \in \mathbb{R}^{n_2 \times n_2}, X \in \mathbb{R}^{n_1}\otimes \mathbb{R}^{n_2}, n_1, n_2 \in \mathbb{N}$. \end{remark}

The statement of the Iterated WIP above is written in the form assuming $f_0 = w(x)l(y)$ which is the case encountered in \ref{sec:homogenisation}. We refer the reader to \cite{kelly2017chaos} for the details involving the full proof for generic $f_0$ with no multiplicative structure.


It is natural to ask whether analogues of the WIP exist for more general stochastic processes. Several results are presented in \cite[Chapter 4]{ethier2005markov} on invariance principles for solutions of a stochastic process solving the martingale problem associated to an operator and measure. To our knowledge, similar generalisations relating to homogenisation are a newly developing field of research. Below we note an invariance principle of Gehringer, Li \cite{Gehringer2022} relating to fractional Brownian motion that motivates some of the extensions to this paper proposed in the conclusion.

Recall a fractional Brownian motion $B_t^H$ with Hurst parameter $H \in (\frac{1}{3}, 1)$ on a complete probability space $(\Omega, \mathcal{F}, \mathbb{P})$ satisfies
$$B_0^H = 0, \quad \mathbb{E}[(B_1^H)^2] = 1, \quad \mathbb{E}\left(B_t-B_s\right)\left(B_u-B_v\right)=\frac{1}{2}\left(|t-v|^{2 H}+|s-u|^{2 H}-|t-u|^{2 H}-|s-v|^{2 H}\right).$$

Define $y_t := \sigma \int_{-\infty}^t e^{-(t-s)}\mathrm{d}B_s^H$ to be the stationary fractional Ornstein-Uhlenbeck process. $y_t$ is well defined as a rough integral when $B_t^H$ is lifted to an a.s. $H$-Hölder rough path $(B_t^H, \mathbb{B}_t^H)$, with iterated integral canonically defined as $\mathbb{B}_t^H = \frac{1}{2}(B_t^H)^2$ as a consequence of $B_t^H$ being scalar. The constant $\sigma > 0$ is fixed in order for all times $t$ we have $y_t \sim N(0,1)$ and furthermore we have $y_t$ solving the Langevin equation (well defined as a RDE)

$$\mathrm{d}y_t = -y_t \mathrm{d}t + \sigma \mathrm{d}B^H_t, \quad y_0 = \sigma \int_{-\infty}^0 e^s \mathrm{d}B^H_s.$$

The fast fractional OU process $y^\varepsilon_t$ is likewise defined, now with the following $H$ dependent scaling and solving:

$$\mathrm{d}y^\varepsilon_t = -\frac{1}{\varepsilon}y^\varepsilon_t \mathrm{d}t + \frac{\sigma}{\varepsilon^H} \mathrm{d}B^H_t, \quad y^\varepsilon_0 = y_0, \quad y^\varepsilon_t = \frac{\sigma}{\varepsilon^H} \int_{-\infty}^t e^{-\frac{1}{\varepsilon}(t-s)} \mathrm{d}B^H_s$$

As a consequence of fractional scaling $B^H_{\frac{t}{\varepsilon}} = \varepsilon^H B_t^H$ we have $y_{\frac{t}{\varepsilon}}$ and $y^\varepsilon_t$ have the same distribution (and thus initial condition).

\begin{theorem*}[Scalar invariance principle for the fractional OU process. Gehringer, Li 2022 \cite{Gehringer2022}]\label{Li2022}

The stochastic process $B^H_\varepsilon(t) := \varepsilon^{H-1} \int_0^t y^\varepsilon_s \mathrm{d}s$ converges in $L^p$ (thus also in distribution) in the space of rough paths $C([0,T), \mathbb{R} \oplus \mathbb{R})$ with the $\kappa$-Hölder topology for any $\kappa \in (\frac{1}{3},H)$, as $\varepsilon \rightarrow 0$ to the path $(\sigma B^H, \frac{1}{2} (\sigma B^H)^2 )$. 

In addition to being stationary, the fractional Ornstein-Uhlenbeck process is ergodic. Let $\nu$ denote the invariant standard Gaussian measure for $y_t$, for any $f \in L^1(\mu)$ we denote $\overline{f} = \int_{\mathbb{R}} f(y) \mu(\mathrm{d}y)$ and have convergence in distribution of the ergodic average:

$$\int_0^t f(y^\varepsilon_s) \mathrm{d}s \overset{\,\,d\,\,\,\,}{\rightarrow} t \overline{f}$$

Let $C^r_b$ denote $r \in \mathbb{N}$ times continuously differentiable functions with derivatives up to order $k$ bounded ($k = 0$ interpreted as bounded continuous functions), for $w \in C^3_b(\mathbb{R}^d, \mathbb{R}^d), p \in C^2_b(\mathbb{R}^d, \mathbb{R}^d), q \in C_b(\mathbb{R}, \mathbb{R})$ the solution $x^\varepsilon$ of the ODE

$$\dot{x}^\varepsilon_t = \varepsilon^{H-1}w(x^\varepsilon_t)y^\varepsilon_t + p(x^\varepsilon_t) q(y^\varepsilon_t), \quad x^\varepsilon_0 = x_0$$

also converges in distribution and rough paths for a Hölder parameter $\kappa \in (\frac{1}{3},H)$ as $\varepsilon \rightarrow 0$ to the solution $X_t$ of the RDE

$$\mathrm{d}X_t = \sigma w(X_t) \mathrm{d}B^H_t + p(X_t)\overline{q}(X_t) \mathrm{d}t$$
\end{theorem*}

\begin{remark} One can view the fractional Ornstein Uhlenbeck process as \eqref{yeps} with $h = 1, b = \sigma, \eta = \dot{B}^H_t$. Taking $H = \frac{1}{2}$ returns us to the Brownian motion regime after performing a change in fast parameter to return to usual scaling (see \cite{Gehringer2022}[Rmk. 3.4] where the intersection with the results of other theorems is commented on). The ergodicity of the fractional OU process replaces the chaos assumptions needed in \eqref{KM2014} and provides the measure $\nu$, which is Gaussian. \end{remark}

The statement of the theorem in \cite{Gehringer2022} in fact much more general and can be phrased using multiplicative noise. Consider a the noise coefficient $\alpha(\varepsilon)w(x^\varepsilon_t)l(y^\varepsilon_t)$ for $l : \mathbb{R} \rightarrow \mathbb{R}$ continuous and centered. The centering condition can be, as in the other cases described in this paper, verified by computing a vanishing integral with the invariant Gaussian measure $\nu$ or can be verified by computing the Itô-Weiner chaos expansion of $l(x) = \sum_{k=0}^\infty c_k H_k(x)$. The chaos expansion of $G$ also must be considered in order to deduce the form of homogenised equations, Li and Gehringer define the Hermite rank $m_k$ to be the smallest $k \in \mathbb{N}$ such that $c_k \neq 0$ which affects the scaling and limit rough driver. It can be shown the scaling function $\alpha$ is picked according to 

$$H^*(m_k) := m_k(H-1) - 1 , \quad \alpha\left(\varepsilon, H^*\left(m_k\right)\right)=\left\{\begin{array}{cl}
\frac{1}{\sqrt{\varepsilon}} & \text { if } H^*\left(m_k\right)<\frac{1}{2}, \\
\frac{1}{\sqrt{\varepsilon|\log (\varepsilon)|}} & \text { if } H^*\left(m_k\right)=\frac{1}{2}, \\
\varepsilon^{H^*(m)-1} & \text { if } H^*\left(m_k\right)>\frac{1}{2} ,
\end{array}\right.$$

so that the noise coefficient converges to $cw(X_t) \mathrm{d}Z^{H^*(m_k), m_k}_t$ for some $c \in \mathbb{R}$ (see \cite{Gehringer2022}[Thm 3.7] for the precise statement of the convergence in each case and how the constant $c$ is determined from $H^*(m_k)$). The process $Z^{H^*(m_k), m_k}_t$ is known as a Hermite process and generalises fractional Brownian motion in the sense that $Z^{H, 1}_t = B^H_t$.

As remarked, a consequence of $\mathcal{Y} = \mathbb{R}$ is that an iterated weak invariance principle for the signature follows from any weak invariance principle for the path due to continuous dependence on the path (half the squaring function). This prevents the occurrence of a Wong-Zakai anomaly occurring, as the Lévy area of all scalar rough paths is identically $0$. The appearance of additional drift in the higher dimensional case can be conjectured and is exhibited in the related works \cite{Gehringer_2022}. In general, this phenomenon is non-specific to Brownian or fractional Brownian motion and solely relates to the regularity of the driving signal. It is shown in \cite{Friz_2015} that for distributional equations of the type

$$m \dot{y} = -M y + \dot{\boldsymbol{Z}}, \quad m > 0, \quad \boldsymbol{Z} \in C^\kappa([0,T], \mathbb{R}^m), \quad M \in \mathbb{R}^{m \times m} \, ,$$

the $m \rightarrow 0$ limit can never exhibit Wong-Zakai anomalies if the driver $\boldsymbol{Z}$ is sufficiently smooth with $\kappa > \frac{1}{2}$. Specialising to the case $\boldsymbol{Z}$ is a fBM with appropriate scaling, one defines the path $B^H_{\varepsilon}(t) := \varepsilon^{H-1}\int_0^ty^\varepsilon_s \mathrm{d}s \in \mathbb{R}^m$, with $y^\varepsilon_t$ satisfying,

$$\mathrm{d}y^\varepsilon_t = -\frac{1}{\varepsilon}Ay^\varepsilon_t \mathrm{d}t + \frac{1}{\varepsilon^H} D \mathrm{d}B^H_t, \quad y^\varepsilon_0 = y_0 = D \int_{-\infty}^0 e^{As} \mathrm{d}B^H_s, \quad y^\varepsilon_t = \frac{1}{\varepsilon^H} D \int_{-\infty}^t e^{-\frac{1}{\varepsilon}A(t-s)} \mathrm{d}B^H_s.$$

The matrices $A, D \in \mathbb{R}^{m \times m}$ are fixed and positive definite. A reasonable conjecture is whether $B^H_{\varepsilon}(t)$ satisfies an iterated weak invariance principle of $(B^H_{\varepsilon}(t) - B^H_{\varepsilon}(s), \mathbb{B}^H_{\varepsilon}(t,s) )\overset{\,\,d\,\,\,\,}{\rightarrow} (B^H_t - B^H_s, \mathbb{B}^H_{t,s} + \mathsf{s}(t-s)) $. As before, the signature $\mathbb{B}^H_{\varepsilon}(t,s)$ is canonically defined as iterated integrals of smooth paths, $\mathbb{B}^H_{t,s}$ defined as the geometric lift of a Gaussian process (see \cite[Appendix E.1]{crisan2022variational}) and $\mathsf{s} \in \mathfrak{s}\mathfrak{o}(m)$ is a deterministic perturbation to the Lévy area which appears conditionally on the antisymmetric part of $A$ akin to the physical Brownian motion case in \cite{DHP2023, Friz_2015}.

\subsubsection{Explicit representations of the drift and diffusion coefficients}\label{sec:homog_formulas}

In this section we recall derivations of the formulae for drift and diffusion coefficients of equation \eqref{homogenisedguess}, used in \cref{sec:homogenisation}. We outline the approach presented in \cite{Dror_2004} \cite{Pavliotis2008} who derive a form for the generator of the limit SDE. The system composed of \eqref{xeps} and \eqref{yeps} decompose over different scales in $\varepsilon$ their contributions to the Kolmogorov backward equation of the full system as follows:

\begin{equation}\label{fokker}\frac{\partial v}{ \partial t} = \frac{1}{\varepsilon^2}\mathcal{L}_1 v + \frac{1}{\varepsilon}\mathcal{L}_2 v + \mathcal{L}_3 v\end{equation}

Where $\mathcal{L}_1v = h \cdot \nabla_y v + \frac{1}{2} bb^T : \nabla_y \nabla_y v$, $\mathcal{L}_2v = f_0 \cdot \nabla_x v$ and $\mathcal{L}_3v = f_1 \cdot \nabla_x v$. We can then seek solutions of the form $v = v_0 + \varepsilon v_1 + \varepsilon^2 v_2 \ldots$ to substitute into \eqref{fokker} leading to the relations $\mathcal{L}_1v_0 = 0, \mathcal{L}_1 v_1 = -\mathcal{L}_2v_0, \mathcal{L}_1v_2 = \frac{\partial v_0}{\partial t} - \mathcal{L}_2v_1 - \mathcal{L}_3v_0$. 

The assumption that the fast $y$ dynamics is Ergodic implies $e^{\mathcal{L}_1t}f \rightarrow \langle f \rangle_\nu$ as $t \rightarrow \infty$ where the average $\langle \cdot \rangle_\nu$ is with respect to some unique invariant measure with a density $\nu = \rho_\infty(y) \mathrm{d}y$ satisfying $\mathcal{L}^*_1 \rho_\infty(y) = 0$.

$\mathcal{L}_1$ is the generator of an ergodic Markov process solely in the variable $y$, one such characterisation of this ergodicity is that its nullspace consists solely of constants forcing $v_0 = v_0(x,t)$. Furthermore, by the Fredholm alternative the ability to solve $\mathcal{L}_1 v_1 = -\mathcal{L}_2 v_0$ enforces an orthogonality condition of $-\mathcal{L}_2 v_0$ with respect to the kernel of $\mathcal{L}^*_1$, which by ergodicity solely consists of $\rho_\infty(y)$. The orthogonality condition becomes $\langle \mathcal{L}_2 v_0 \rangle_\nu = 0$ and since $v_0$ has no dependence on $y$ reduces to the sufficient condition that $\langle f_0 \rangle_\nu = 0$, this is the well known \emph{centering condition} that the slow dynamics \eqref{xeps} must satisfy to apply homogenisation. 

Provided all the required assumptions are met we are then able to write $v_1 = -\mathcal{L}_1^{-1}\mathcal{L}_2v_0$. This can be substituted into the $\mathcal{O}(1)$ relations in the expansion and apply the Fredholm alternative again when solving for $v_2$ to obtain the following condition

\begin{equation}\label{secretFP}\frac{\partial v_0}{\partial t}=-\left\langle\mathcal{L}_2 \mathcal{L}_1^{-1} \mathcal{L}_2 v_0\right\rangle_\nu+\left\langle\mathcal{L}_3 v_0\right\rangle_\nu.\end{equation}

Observe that \eqref{secretFP} is a new backward Kolmogorov equation for $v_0(x,t)$ consisting of terms involving a first order differential operator (in $x$) arising from $\mathcal{L}_3$ and a second order $x$ differential operator term from twice application of $\mathcal{L}_2$. The $y$ dependence of the coefficients of $\mathcal{L}_2, \mathcal{L}_3$ and the linear operator $\mathcal{L}^{-1}_1$ are averaged in a manner that preserved the leading order contributions from $\mathcal{L}_2$ and $\mathcal{L}_3$. 

It remains to extract a SDE associated to backward Kolmogorov equation \eqref{secretFP}. First note that we can express $\mathcal{L}_1^{-1} \mathcal{L}_2 v_0 = f \cdot \nabla_x v_0$, where $f$ is a solution of the cell problem $\mathcal{L}_1 f = f_0(x,y)$. One can then apply $\mathcal{L}_2$ using Leibniz's rule and match indices to write the composition of operators in Kolmogorov form.

\begin{equation}
\mathcal{L}_2 \mathcal{L}_1^{-1} \mathcal{L}_2 v_0 = \frac{1}{2} (2f_0f^T) : \nabla_x \nabla_x v_0 +  (\nabla_x f )f_0 \cdot \nabla_x v_0 \label{generatorcontributions}\end{equation}

Note that this formula contributes both drift and martingale components whilst the $\mathcal{L}_3v_0$ portion contributes as $f_1 \cdot \nabla_x v_0$. When taking averages the drift coefficient of the SDE is indeed given by the average of $f_1$, \emph{plus} an additional noise induced drift of $-\langle  \nabla_x f f_0 \rangle_\nu$ from \eqref{generatorcontributions}.

In order to extract the diffusion coefficient we solve $\Sigma \Sigma^T = -2\langle f_0 f^T\rangle_\nu$. Note that this is well defined if $f_0 f^T$ is negative definite which follows from the negative definiteness of $\mathcal{L}_1$ (see \cite{Pavliotis2008}[Thm 11.3]). 

\begin{remark} The derivation of the coefficients $\overline{f}_1, \Sigma$ via the generator are highly dependent on being able to solve the associated cell problem $\mathcal{L}_1 f = f_0$. Solving the cell problem implicitly requires $b \neq 0$. Equation \eqref{fokker} is valid as a backward equation for an ODE when no noise is present, however the operators $\mathcal{L}_i$ are not uniformly elliptic and in general the unique ergodic measure will not have a density with respect to the Lebesgue measure. One approach to extend this argument is considered in \cite{Engel2021} by introducing the regularised problem where $b = \sqrt{\delta} I$ and proceed as above. The $\delta \rightarrow 0$ limit requires the collection of generators to form a strongly continuous semigroup.\end{remark} 

The additional noise induced drift of $-\langle \nabla_x f f_0 \rangle_\nu$ can be thought of as governing the choice of stochastic integral $* \mathrm{d} W$. It has been shown in \cite{gottwald2022timereversibility, kelly2017chaos} that when the noise has a multiplicative structure $f_0(x_\varepsilon, y_\varepsilon) = w(x_\varepsilon)l(y_\varepsilon)$ this decomposes as the Itô-Stratonovich correction plus an additional drift depending on the Wong-Zakai anomaly. This term in the context of approximations of Brownian motion also appears in more general situations than WIP discussed here, see \cite{Gehringer_2022}. For general choices of $f_0$ the problem of what the additional drift means for the stochastic integral remains open, a sense of this behaviour can be illustrated in \cite{Pavliotis2008}[Sec 11.7.6]. 
Another approach for the case that $h(x,y) = h(y), b(x,y) = b(y)$ which bypasses solution of the cell problem is to represent the operator $\mathcal{L}^{-1}_1$ in integral form as $-\int_0^\infty e^{\mathcal{L}_1 t }\mathrm{d}t$. Note that if $g$ is centered, meaning 
 $\langle g \rangle = 0$ then we verify this expression gives a right inverse: 
 
 $$\mathcal{L}\left(-\int_0^\infty e^{\mathcal{L}_1 t}g \mathrm{d}t \right) = -\int_0^\infty \frac{\partial}{\partial t} e^{\mathcal{L}_1 t}g \mathrm{d}t = g - \lim\limits_{t \rightarrow \infty}e^{\mathcal{L}_1 t} g = g - \lim\limits_{t \rightarrow \infty} \int g \rho(y,t)\mathrm{d} y = g - \int g \rho_\infty(y)\mathrm{d} y = g.$$ 
 
 A similar argument shows the left inverse, since $\mathcal{L}_1$ commutes with it's one parameter subgroup of operators. We can then insert this integral form into the orthogonality condition \eqref{secretFP}
\begin{align*}
-\left\langle\mathcal{L}_2 \mathcal{L}_1^{-1} \mathcal{L}_2 v_0\right\rangle_\nu &= \left\langle f_0 \cdot \nabla \int_0^\infty e^{\mathcal{L}_1t}f_0 \cdot \nabla v_0 \mathrm{~d} t\right\rangle_\nu\\ 
&= \left\langle\int_0^{\infty} f_0 \mathrm{e}^{\mathcal{L}_1 t} f_0^T \mathrm{~d} t\right\rangle_\nu : \nabla_x\nabla_x v_0+\left\langle\int_0^{\infty} \nabla_x\left(\mathrm{e}^{\mathcal{L}_1 t} f_0\right) f_0 \mathrm{~d} t\right\rangle_\nu \cdot \nabla_x v_0.
\end{align*}

The above expression is now in terms of $\mathcal{L}_1$ and $f_0$ rather than $\mathcal{L}^{-1}_1$ or the cell problem solution $f$. Note that directly from the definition of the infinitesimal generator and the Feller semigroup we have that $e^{t\mathcal{L}_1}f_0(\cdot, y) = \mathbb{E}[f_0(\cdot, \phi_t(y))]$, where this expectation is respect to the probability measure $\mathbb{P}$ of the underyling Brownian motion\footnote{Compare this to the one-parameter semigroup $e^{t\mathcal{L}_1} f(\cdot,y) := f(\cdot,\phi_t(y)$ appearing in the deterministic homogenisation theorem \eqref{KM2014} defined by flowing along the fast dynamics of the random initial condition $y$.} and $\phi_t$ is the stochastic flow of the SDE \eqref{yeps}. We then arrive at one possible explicit representation of the SDE \eqref{homogenisedguess}:

\begin{align}\label{homogformula}
\overline{f}_1(x) &= \langle{f_1(x,y)}\rangle_{\nu(\mathrm{d}y)} + \left\langle \int_0^\infty \mathbb{E} \Big[\nabla_x f_0(x,\phi_t(y)) f_0(x,y)\Big] \mathrm{~d} t\right\rangle_{\nu(\mathrm{d}y)}\\
\Sigma\Sigma^T &= \frac{1}{2}(A + A^T)\\
A &= 2 \left\langle\int_0^{\infty} \mathbb{E}\left(f_0(x, y) \otimes f_0\left(x, \phi_t(y)\right)\right) \mathrm{d} t\right\rangle_{\nu(\mathrm{d}y)} 
\end{align}

The key assumption that the fast dynamics was $x$ independent allowed us to separate expectations and gradients when substituting the integral representation of $\mathcal{L}^{-1}_1$. By defining the product measure $\nu \otimes \mathbb{P}$ we can interchange the order of integrals in \eqref{homogformula}

\begin{align}\label{interchangedhomogformula}
\overline{f}_1(x) &= \langle{f_1(x,y)}\rangle_{\nu(\mathrm{d}y)} + \int_0^\infty \mathbb{E}^{\nu \otimes \mathbb{P}} \Big[\nabla_x f_0(x,\phi_t(y)) f_0(x,y)\Big] \mathrm{~d} t\\
\Sigma\Sigma^T &= \frac{1}{2}(A + A^T)\\
A &= 2 \int_0^{\infty} \mathbb{E}^{\nu \otimes \mathbb{P}}\left(f_0(x, y) \otimes f_0\left(x, \phi_t(y)\right)\right) \mathrm{d} t 
\end{align}

Computations with $\mathbb{E}^{\nu \otimes \mathbb{P}}$ are carried out by averaging stochastic integrals to zero under $\mathbb{E}$, and averaging $y$ with respect to the invariant measure as a double integral. The expectation with respect to the product measure in correction term and diffusion tensor formulae in \eqref{interchangedhomogformula} are analogous to the bilinear form $\mathfrak{B}$ seen in \cite{kelly2017chaos} discussed in the above subsection now from the viewpoint of the Feller generator of a stochastic differential equation.
\subsection{Lévy's Stochastic Area contribution}
\label{Sec:Lévy Area contribution in higher order Stratonovich Integrators}
The Stratonovich SDE, is effected by the commutator of the vector fields when higher order integrators are used. This occurs even in absence of a Wong-Zakai anomaly drift, when Lévy's stochastic area is approximated and used in a numerical method to capture and represent higher order terms in a (Wagner Platen/Taylor Stratonovich) expansion. In this section, we informally show how the Wong-Zakai anomaly drift term emerges naturally if one peturbs the Lévy area in a Stratonovich-Taylor expansion when viewed as a enhanced Brownian motion driven system. We also relate the point vortex system to a state valued SDE, and note how the terms in the Taylor Stratonovich expansion coincide with the point vortex approximation at the level of the continuum, for example the Itô-Stratonovich correction for the point vortex system is the usual Itô-Stratonovich correction for vector fields but evaluated and summed over the point vortices.

We are interested in converting the system of stochastically driven particles $\lbrace \b x_{\alpha}\rbrace_{\alpha\in 1,,...,n}\in \mathbb{R}^2$, 
\begin{align}
\b x_{\alpha}(t) &= \b x_{\alpha}(0) + \int_{0}^{t} \b u(\b x_{\alpha}(s),s)ds + \sum_{j=1}^{m}\int_{0}^{t}\b \xi_{j}(\b x_{\alpha}(s)) \circ \mathrm{d}W^{j}_t. 
\end{align}
into a single Stratonovich SDE system of state variables. To do so we let $\b q = (x_1,...,x_{\alpha},...,x_n ,y_1,...,y_{\alpha},...,y_n)^T\in \mathbb{R}^{2n}=\mathbb{R}^{d}$ denote the vector of state variables, the Stratonovich drift vector is made up of the stacked components of velocity evaluated at each point vortex 
$$\b f(\b q) = ( u_1(\b q),...,u_{\alpha}(\b q),..., u_n(\b q) ,v_1(\b q),...,v_{\alpha}(\b q),...,v_n(\b q))^T,$$
$G$ is a matrix whose $j$-th column $\b g_j$ is defined by the stacked components of $\b \xi_{j} = (\xi^{x}(\b q),\xi^{y}(\b q))^T$, where superscript $x,y$ denotes the $x,y$ component of $\b \xi$ in space, explicitly $$\b g_j (q) = (\xi^x_{j}(\b x_1;\b q),...,\xi^x_{j}(\b x_{\alpha};\b q),...,\xi^x_{j}(\b x_n;\b q) ,\xi^y_{j}(\b y_1;\b q),...,\xi^y_{j}(\b y_{\alpha};\b q),...,\xi^y_{j}(\b y_n;\b q))^T.$$ 
We denote $g^{k,j}$ the $j$th column, evaluated at position $k$ in state variable $\b q$, such that $g^{k,j} = (\xi^{x})^{k}_{j}$, for $k \in \lbrace 1,...,n\rbrace$, $g^{k,j} = (\xi^{y})^{k-n}_{j}$ for $k \in \lbrace 1+n,...,2n \rbrace$, recalling $q_k = x_{k}$ for $k \in \lbrace1,...,n\rbrace$ and $q_k = y_{k-n}$ for $k \in \lbrace 1+n,...,2n\rbrace$.
We have written the system as the Stratonovich system, 
\begin{align}
\mathrm{d} \b q_t = \b f (\b q_t)\mathrm{d}t + G(\b q_t)\circ \mathrm{d}W_t,
\end{align}

Where $\b q \in \mathbb{R}^{2n}$,$\b f(\b q):\mathbb{R}^{2n}\rightarrow \mathbb{R}^{2n}$, $G(\b q):\mathbb{R}^{2n} \rightarrow \mathbb{R}^{2n\times m}$, is integrated against the components of $W$ a $m$-dimensional signal. We define iterated integrals over an arbitrary time interval $[t^{n},t^{n+1}]$ over the index set $(j_1,...,j_{l})$, by
$$
J_{\left(j_1, \ldots, j_l\right)}=\int_{t^n}^{t^{n+1}} \cdots \int_{t^n}^{s_2} \circ d W_{s_{1}}^{j_1} \cdots \circ d W_{s_l}^{j_l}.
$$
and $J_{0} = \int_{t^n}^{t^{n+1}}\mathrm{d}t$. We now consider a Stratonovich Taylor expansion following \cite{kloeden1992stochastic} noting the following properties of the nested Stratonovich integrals 
$J_{\left(j_1, j_1\right)}=\frac{1}{2}\left(J_{j_1}\right)^2$ and $J_{\left(j_1, j_2\right)}= -J_{\left(j_2, j_1\right)} $ to seperate out the Itô-Stratonovich term from the Lévy area 
\begin{align}
q_{n+1}^k=q_n^k+f^k J_0+\sum_{j=1}^m g^{k, j} J_{j}+ \sum_{1\leq j_1,j_2\leq m} L^{j_1} g^{k, j_2} \operatorname{Sym}(J_{j_{1},j_{2}})  + \frac{1}{2}\sum_{1\leq j_1,j_2\leq m} (L^{j_1} g^{k,j_2} - L^{j_2} g^{k,j_1} )J_{\left(j_1, j_2\right)} + ... ,
\end{align}
where $k$ denotes the $k$th component of the system and $L$ is defined in $\cref{def:L}$. If one views the above expression as a enhanced Brownian motion driven system $( W, \mathbb{W})$ and includes a perturbation to the Lévy area term by the antisymmetric term $\Delta t\mathsf{s}^{j_1,j_2}\in \mathfrak{so}(m)$ where $\mathsf{s}^{j_1,j_2}$ has no temporal dependence, then one heuristically treats this as a drift in the following way,
\begin{align}
\begin{split}
q_{n+1}^k = q_n^k+\left( f^k + \frac{1}{2}\sum_{1\leq j_1,j_2\leq m} (L^{j_1} g^{k,j_{2}} - L^{j_2} g^{k,j_{1}} )\mathsf{s}^{j_1,j_2} \right)J_0 + \sum_{j=1}^m g^{k, j} J_{j}+ \frac{1}{2}\sum_{j_1=1}^m L^{j_1} g^{k, j_1} J_{j_{1}}^2   \\
+ \sum_{1\leq j_1\neq j_2\leq m} L^{j_1} g^{k, j_2} \operatorname{Sym}(J_{j_{1},j_{2}})+\frac{1}{2}\sum_{1\leq j_1,j_2\leq m} (L^{j_1} g^{k,j_{2}} - L^{j_2} g^{k,j_{1}} )J_{\left(j_1, j_2\right)}  + ... .
\end{split}\label{EQ: Stratonovich Taylor Expansion}
\end{align}
Where the last term can be written in terms of the Lévy area
$$
\frac{1}{2}\sum_{1\leq j_1,j_2\leq m} (L^{j_1} g^{k,j_{2}} - L^{j_2} g^{k,j_{1}} )J_{\left(j_1, j_2\right)} = \frac{1}{2}\sum_{1\leq j_1,j_2\leq m} (L^{j_1} g^{k,j_{2}} - L^{j_2} g^{k,j_{1}} )A_{\left(j_1, j_2\right)}
$$
where $A_{j_{1}, j_{2}} := J_{j_1, j_2}-\operatorname{Sym}(J_{j_1, j_2})$, denotes the skew symmetric part of the nested integration of multidimensional Brownian motion, where the symmetric operator is defined as $\operatorname{Sym}(M):=\frac{1}{2}(M+M^{T})$.

We first define the Stratonovich operator $L^{j}$ from \cite{kloeden1992stochastic}, and relate it to the point vortex approximation $\b q\in \mathbb{R}^{d}\mapsto \b x_{\alpha} \in \mathbb{R}^{N\times 2}$. 
\begin{align}
L^{j} &:= \sum_{k=1}^{2n} g^{k,j}\frac{\partial}{\partial q^{k}} \label{def:L} = \sum_{k=1}^{n} (\xi^x)^k_j \frac{\partial}{\partial x^{k}}+ \sum_{k=1}^{n} (\xi^y)_j^k\frac{\partial}{\partial y^{k}} =\sum_{\alpha = 1}^{n}\b \xi_{j}(\b x_{\alpha}) \cdot \nabla_{\alpha}
\end{align}
This allows us to calculate the $k$th component of the following expression 
\begin{align}
    L^{j_1} g^{k,j_2} &= \sum_{i=1}^{2n} g^{k,j_{1}}\frac{\partial}{\partial q^{i}} g^{k,j_2} = \sum_{\alpha=1}^{n} (\b \xi_{j_1}(\b x_{\alpha})\cdot \nabla_{\alpha} ) \b \xi^k_{j_2}(\b x_{\alpha}),
\end{align}
allowing the identification that the Itô-Stratonovich correction for the point vortex system is
\begin{align} \frac{1}{2}\sum_{j=1}^m \sum_{\alpha=1}^n (\b \xi_{j}\cdot \nabla_{\alpha} ) \b \xi_{j},
\end{align}
 the usual Itô-Stratonovich correction for vectorfields but evaluated and summed over the point vortices of the system. Similarly the $k$th component of the commutator can be computed as
\begin{align}
 L^{j_1} g^{k,j_2} - L^{j_2} g^{k,j_1} &=   \sum_{\alpha=1}^{n}(\b \xi^i_{j_1}\cdot \nabla_{\alpha}  ) \xi^k_{j_2} - (\b \xi^i_{j_2}\cdot \nabla_{\alpha} )\xi^k_{j_1}. \end{align}
 The commutator for the point vortex system is the regular commutator of stochastic vector fields but evaluated at the positions of the point particles and summed over the points in the system. Identifying the Lévy area and Wong Zakai anomaly drift contribution take the following form in the point vortex system
\begin{align}
\frac{1}{2}\sum_{1\leq j_1,j_2\leq m}\sum_{\alpha=1}^{n}[\b \xi_{j_1},  \b \xi_{j_2}] \left(A_{\left(j_1, j_2\right)}+ \mathsf{s}^{j_1, j_2}J_0\right).
\end{align}

Despite suggestive notion we haven't discretised in time, however the form \cref{EQ: Stratonovich Taylor Expansion} indicates that the drift and Lévy's stochastic area contribution should be treated differently in the context of a numerical scheme.


\subsection{Numerical approximation of the Lévy area}\label{sec:Numerical simulation of the Lévy area}
Following \cite{kloeden1992stochastic} (Chapter 5) and Kloeden, Platen and Wright \cite{kloeden1992approximation} the approximation of arbitrary higher order Stratonovich nested integrals can be considered, by using the relationship between different higher order terms found in the Stratonovich Taylor expansion, and the Fourier series expansion of a Brownian bridge. They consider the  truncated Fourier series of the Brownian bridge process
$\left\lbrace W^{j}_t-\frac{t}{\Delta t} W^j_{\Delta t}, 0 \leq t \leq \Delta t \right\rbrace
$ up to $K$ terms, $\forall j\in \lbrace 1,...m\rbrace$ where $W_t$ is the $m$-dimensional Wiener process $W_t=\left(W_t^1, \ldots, W_t^m\right)$ over $t \in [0, \Delta t]$. This component-wise takes the form
\begin{align}
W_t^j-\frac{t}{\Delta t} W_{\Delta t}^j=\frac{1}{2} a_{j, 0}+\sum_{r=1}^{K}\left(a_{j, r} \cos \left(\frac{2 r \pi t}{\Delta t}\right)+b_{j, r} \sin \left(\frac{2 r \pi t}{\Delta t}\right)\right),
\end{align}
where the random coefficients
\begin{align}
a_{j, r} &=\frac{2}{\Delta t} \int_0^{\Delta t}\left(W_s^j-\frac{s}{\Delta} W_{\Delta t}^j\right) \cos \left(\frac{2 r \pi s}{\Delta t}\right) d s \sim \sqrt{\frac{\Delta t}{2 \pi^2 r^2}} N\left(0 ; 1\right)\\
b_{j, r} &=\frac{2}{\Delta t} \int_0^{\Delta t}\left(W_s^j-\frac{s}{\Delta t} W_{\Delta t}^j\right) \sin \left(\frac{2 r \pi s}{\Delta t}\right) d s \sim \sqrt{\frac{\Delta t}{2 \pi^2 r^2}} N\left(0 ; 1\right)
\end{align}
for $j=1, \ldots, m$ and $r=0,1,2,...K$ are subsequently used in the approximation of a truncated Lévy area, defined below
\begin{align}
A^{K}_{j_1, j_2}  &=\frac{\pi}{\Delta t} \sum_{r=1}^{K} r\left(a_{j_1, r} b_{j_2, r}-b_{j_1, r} a_{j_2, r}\right), \\
&=\frac{1}{2 \pi} \sum_{r=1}^K \frac{1}{r}\left(\zeta_{j_1, r} \eta_{j_2, r}-\eta_{j_1, r} \zeta_{j_2, r}\right).
\end{align}
Where we have made use of the standard intermediary i.i.d Gaussian random variables
  $\zeta_{j, r}=\sqrt{\frac{2}{\Delta t}} \pi r a_{j, r}$, $ \eta_{j, r}=\sqrt{\frac{2}{\Delta t}} \pi r b_{j, r}$, in the approximation of the Lévy area. This approximation to the Lévy area has a mean square error bound found in \cite{kloeden1992approximation}, $
    \mathbb{E}(|A^{K} - A^{\infty}|^2) \leq \frac{\Delta t^2}{ 2\pi^2 K}$ such that the truncated Brownian bridge process should have terms inversely proportional to the timestep $K\approx \Delta t^{-1}$ in order to preserve the strong order of convergence for a strong order one scheme. Since our timestep is $1/500$, we use $K>500$ as the truncation of our Brownian bridge process. This is computed beforehand and used as $A_{j_1,j_2} = J_{j_1,j_2}$ for $j_1\neq j_2$, in the simulations requiring Lévy area approximation.

\subsection{Verification of stochastic commutation relations under translation}\label{sec:supplementary calcs}

We verify the vanishing of Poisson brackets used in \ref{Example:Stratonovich Point Vortex System} without having to appeal to symmetries or a stochastic variant of Noether's theorem. We remind the reader the general definition for $\Gamma \neq 0$

$$\b x_c(\lbrace \b x_{\beta}\rbrace_{\forall \beta}) = \left(\Gamma^{-1} \sum_{\beta = 1}^n \Gamma_{\beta}x_{\beta}, \,\Gamma^{-1} \sum_{\beta = 1}^n \Gamma_{\beta}y_{\beta}\right)^T, \quad \Gamma = \sum_{\beta = 1}^n \Gamma_\beta, \quad \nabla_\alpha \b x_c(\lbrace \b x_{\beta}\rbrace_{\forall \beta}) = 
 \frac{1}{\Gamma}\begin{pmatrix}
\Gamma_\alpha & 0 \\
0 & \Gamma_\alpha
\end{pmatrix}.
$$


For \cref{sec: case study} $n = 3, \Gamma_\alpha \equiv 1 \,\,\, \forall \alpha \in {1,2,3}$. We shall denote $\b x_c (\lbrace \b x_{\beta}\rbrace_{\forall \beta})$ as $\b x_c$ for brevity.

To derive $\{ H, \psi_1 \circ \mathrm{d} W^1_t\} = 0$ we require equal strength point vortices $\Gamma_\alpha \equiv K$ and initial conditions and a number of particles to ensure a stable circular configuration, shown to be $1 \leq n \leq 7$ in \cite{Montaldi2013}.
$$\{ H, \psi_1 \circ \mathrm{d} W^1_t\} = -\sum_{\alpha = 1}^n\nabla_\alpha H(\b x_\alpha) \cdot \b \xi_1(\b x_\alpha)\circ \mathrm{d}W^1_t = -\sum_{\alpha = 1}^n \frac{K^2}{2 \pi}\sum_{\beta = 1, \beta \neq \alpha}^n \frac{(\b x_\alpha - \b x_\beta)}{\|\b x_\alpha - \b x_\beta\|^2}  \cdot (\b x_\alpha - \b x_c)^\perp Ar(1-\frac{1}{n})e^{-\frac{r}{2}\|\b x_\alpha - \b x_c\|^2}$$

$$= -\frac{Ar(1 - \frac{1}{n}) K^2}{2 \pi } \left( \sum_{\alpha,\beta =1, \alpha\neq \beta}^n \frac{y_{\alpha} x_{\beta} - x_{\alpha} y_{\beta}}{\|\bx_{\alpha} - \bx_{\beta}\|^2} e^{-\frac{r}{2}\|\b x_\alpha - \b x_c\|^2} + \frac{y_{c} x_{\alpha} + y_\beta x_c - x_{c} y_{\alpha} - x_\beta y_c}{\|\bx_{\alpha} - \bx_{\beta}\|^2} e^{-\frac{r}{2}\|\b x_\alpha - \b x_c\|^2} \right) \circ \mathrm{d}W^1_t$$

The first term of this sum vanishes by a similar argument to the case $\b x_c = 0$, with no translation presented in \cref{Example:Stratonovich Point Vortex System}. The intervorticial and deviation from the center are constant and independent of indicies, it follows we can define the constant: \[\frac{Ar(1 - \frac{1}{n}) K^2}{2 \pi } \frac{e^{-\frac{r}{2}\|\b x_\alpha - \b x_c\|^2}}{\|\bx_{\alpha} - \bx_{\beta}\|^2} = C, \quad \text{if}, \quad \Gamma_\alpha = K \,\,\, \forall \alpha, \quad 1 \leq n \leq 7 \,,\]

and show the second term vanishes,

$$-C \left( \,\, \sum_{\alpha,\beta =1, \alpha\neq \beta}^n {y_{c} x_{\alpha} + y_\beta x_c - x_{c} y_{\alpha} - x_\beta y_c} \right) \circ \mathrm{d}W^1_t = -C(ny_c x_c + nx_c y_c - nx_c y_c - nx_c y_c) \circ \mathrm{d}W^1_t = 0 \, .$$

And we thus conclude $\{ H, \psi_1 \circ \mathrm{d} W^1_t\} = 0$ as required. To show $\{T_x, \psi_1 \circ \mathrm{d}W^1\} = 0$ note that $\psi_1$ is scalar function of distance from the center $\b x_c$. Let us denote with $\prime$ the derivative with respect to $\|\b x - \b x_c\|^2$, we may make use of the chain rule property of the Poisson bracket.
$$\{T_x, \psi_1 \circ \mathrm{d}W^1\} = \left\{T_x, A\exp \left( -\frac{r}{2} \| (\cdot) - \b x_c \|^2\right) \circ \mathrm{d}W^1_t\right\} = -\left\{T_x, \| (\cdot) - \b x_c \|^2 \right\} \frac{r}{2} \psi_1^\prime \circ \mathrm{d}W^1_t$$
It suffices to show that $T_x$ commutes with the distance from the center of vorticity. The assumption of equal strength vortices $\Gamma_\alpha \equiv K$ is required in order to conclude $\nabla_\alpha (\b x_\alpha - \b x_c) = I_{2 \times 2}(1 - \frac{1}{n})$ $\forall \alpha$. 
\begin{align*}-\left\{T_x, \| (\cdot) - \b x_c \|^2 \right\} &= \sum_{\alpha=1}^n \left(\nabla_\alpha \sum_{\beta=1}^n x_\beta\right) \cdot \nabla_\alpha^\perp \| \b x_\alpha - \b x_c \|^2 = -\sum_{\beta=1}^n \partial_{y_\beta}\| \b x_\beta - \b x_c \|^2 = -\sum_{\beta=1}^n 2(y_\beta - y_c)(1-\frac{1}{n})\\ 
&=  -2(1- \frac{1}{n}) \sum_{\beta=1}^n y_\beta - y_c = -2(1- \frac{1}{n}) (ny_c - n y_c) 
= 0 
\end{align*}
The same argument mutandis mutatis gives $\{T_y, \psi_1 \circ \mathrm{d}W^1\} = 0$. The proof of $\{R, \psi_1 \circ \mathrm{d}W^1\} = 0$ is also similar, also requiring $\nabla_\alpha (\b x_\alpha - \b x_c) = I_{2 \times 2}(1 - \frac{1}{n})$ to allow cancelling sums of coordinates of points:

$$\left\{R,  A\exp \left( -\frac{r}{2} \| (\cdot) - \b x_c \|^2_2\right) \right\} = -\left\{ \sum_{\beta=1}^n x^2_\beta + y^2_\beta , \| (\cdot) - \b x_c \|^2 \right\} \frac{r}{4} \psi_1^\prime \circ \mathrm{d}W^1_t$$

And so it suffices to show the angular impulse commutes with distance from the center of vorticity.
\begin{align*}
-\left\{ \sum_{\beta=1}^n x^2_\beta + y^2_\beta  , \| (\cdot) - \b x_c \|^2 \right\} &= \sum_{\beta=1}^n -2x_\beta\partial_{y_\beta}\| \b x_\beta - \b x_c \|^2 + 2y_\beta\partial_{x_\beta}\| \b x_\beta - \b x_c \|^2\\
&= \sum_{\beta=1}^n -4x_\beta (y_\beta - y_c)(1-\frac{1}{n}) + 4y_\beta (x_\beta - x_c)(1-\frac{1}{n})\\ 
&= 4(1-\frac{1}{n})\sum_{\beta=1}^n x_\beta y_\beta + x_\beta y_c- x_\beta y_\beta - y_\beta x_c\\
&= 4(1-\frac{1}{n})\left(y_c\sum_{\beta=1}^n x_\beta  -  x_c\sum_{\beta=1}^n y_\beta\right) = 4n(1-\frac{1}{n})(x_c y_c - x_c y_c) = 0.
\end{align*}


\end{document}